\documentclass[11pt]{article}
\usepackage{amssymb}

\usepackage[latin1]{inputenc}

\usepackage{times}

\title{Classical and Effective Descriptive Complexities of  $\omega$-Powers}

\author{Olivier Finkel  \\{\it   Equipe de Logique Mathématique}
  \\ CNRS et  Universit\'e Paris Diderot Paris 7
 \\ UFR de Mathématiques case 7012, site Chevaleret,\\75205 Paris Cedex 13, 
 France.\\ finkel@logique.jussieu.fr 
\and Dominique Lecomte
 \\{\it Equipe d'Analyse Fonctionnelle}
\\ Universit\' e Paris 6, tour 46-0, bo\^\i te 186,
\\  4, place Jussieu, 75 252 Paris Cedex 05, France.
\\ dominique.lecomte@upmc.fr
\\and 
\\ Universit\'e de Picardie, I.U.T. de l'Oise, site de Creil,
\\ 13, all\'ee de la fa\"\i encerie, 60 107 Creil, France.
}

\date{}

\def\ufootnote#1{\let\savedthfn\thefootnote\let\thefootnote\relax
\footnote{#1}\let\thefootnote\savedthfn\addtocounter{footnote}{-1}}

\newcommand{\bormxi}{{\bf\Pi}^{0}_{\xi}}
\newcommand{\bormlxi}{{\bf\Pi}^{0}_{<\xi}}

\newcommand{\bormone}{{\bf\Pi}^{0}_{1}}
\newcommand{\ca}{{\bf\Pi}^{1}_{1}}
\newcommand{\bormtwo}{{\bf\Pi}^{0}_{2}}

\newcommand{\bormeta}{{\bf\Pi}^{0}_{\eta}}

\newcommand{\bormxipo}{{\bf\Pi}^{0}_{\xi +1}}

\newcommand{\borel}{{\bf\Delta}^{1}_{1}}
\newcommand{\Borel}{{\it\Delta}^{1}_{1}}
\newcommand{\borone}{{\bf\Delta}^{0}_{1}}

\newcommand{\borthree}{{\bf\Delta}^{0}_{3}}
\newcommand{\boraone}{{\bf\Sigma}^{0}_{1}}
\newcommand{\boratwo}{{\bf\Sigma}^{0}_{2}}

\newcommand{\boraxi}{{\bf\Sigma}^{0}_{\xi}}
\newcommand{\ana}{{\bf\Sigma}^{1}_{1}}

\newcommand{\Ana}{{\it\Sigma}^{1}_{1}}
\newcommand{\Boraone}{{\it\Sigma}^{0}_{1}}
\newcommand{\Borone}{{\it\Delta}^{0}_{1}}
\newcommand{\Bormone}{{\it\Pi}^{0}_{1}}
\newcommand{\Bormtwo}{{\it\Pi}^{0}_{2}}
\newcommand{\Ca}{{\it\Pi}^{1}_{1}}

\newcommand{\Boraxi}{{\it\Sigma}^{0}_{\xi}}
\newcommand{\Bormxi}{{\it\Pi}^{0}_{\xi}}
\newcommand{\Bormxipo}{{\it\Pi}^{0}_{\xi +1}}

\newcommand{\Boran}{{\it\Sigma}^{0}_{n}}

\newcommand{\boran}{{\bf\Sigma}^{0}_{n}}
\newcommand{\boraxim}{{\bf\Sigma}^{0}_{\xi^-}}
\newcommand{\Boratwo}{{\it\Sigma}^{0}_{2}}

\newcommand{\Borthree}{{\it\Delta}^{0}_{3}}

\newcommand{\Boranpo}{{\it\Sigma}^{0}_{n+1}}

\newcommand{\Borxi}{{\it\Delta}^{0}_{\xi}}
\newcommand{\borxi}{{\bf\Delta}^{0}_{\xi}}
\newcommand{\Bormlxi}{{\it\Pi}^{0}_{<\xi}}
\newcommand{\borapeap}{{\bf\Sigma}^{0}_{1+\eta_{\alpha ,p}}}
\newcommand{\borapeapn}{{\bf\Sigma}^{0}_{1+\eta_{\alpha ,p,n}}}
\newcommand{\Borapeap}{{\it\Sigma}^{0}_{1+\eta_{\alpha ,p}}}

\newcommand{\Bormn}{{\it\Pi}^{0}_{n}}

\newcommand{\Born}{{\it\Delta}^{0}_{n}}

\newcommand{\boreta}{{\bf\Delta}^{0}_{\eta}}
\newcommand{\pn}{{\bf\Delta}^{1}_{n}}
\newcommand{\pan}{{\bf\Sigma}^{1}_{n}}
\newcommand{\pmn}{{\bf\Pi}^{1}_{n}}
\newcommand{\panpo}{{\bf\Sigma}^{1}_{n+1}}
\newcommand{\hs}{\hspace{12mm}

\noi}
\newcommand{\noi}{\noindent}
\newcommand{\ol}{ $\omega$-language}
\newcommand{\om}{\omega}
\newcommand{\Si}{\Sigma}

\newcommand{\Sio}{\Sigma^\omega}

\newcommand{\orl}{ $\omega$-regular language}

\newtheorem{thm} {Theorem} [section]
\newtheorem{defi} [thm] {Definition}
\newtheorem{cor} [thm] {Corollary}
\newtheorem{lem} [thm] {Lemma}
\newtheorem{prop} [thm] {Proposition}

\begin{document}

\maketitle

\ufootnote{This paper is an extended version of a conference paper which  appeared in the Proceedings
 of the 16th EACSL Annual Conference on Computer Science and Logic, CSL 07 \cite{csl07}.
Part of the results in this paper have been also presented at the International Conference Computability in Europe, CiE 07, Siena, Italy, June 2007.}

\begin{abstract} 
\noi  We prove that, for each countable ordinal $\xi\!\geq\! 1$, there exist some
$\boraxi$-complete $\omega$-powers, and some $\bormxi$-complete $\omega$-powers, extending previous works on the topological complexity 
of $\om$-powers \cite{Fin01a,Fin03a,Fin04,Lecomte01,Lecomte05,DF06}.
We prove effective versions of these results; in particular, for each recursive ordinal $\xi\! <\!\omega_1^{CK}$ there exist some recursive sets 
$A \subseteq 2^{< \om}$ such that   $A^\infty\!\in\!\Bormxi\!\setminus\!\boraxi$  (respectively, $A^\infty\!\in\!\Boraxi\!\setminus\!\bormxi$), where 
$\Bormxi$ and $\Boraxi$ denote classes of the hyperarithmetical hierarchy. 
To do this, we prove effective versions of a result by Kuratowski, describing a $\bormxi$ set as the range of a closed subset of the Baire space 
$\omega^\omega$ by a continuous bijection. This leads us to prove closure properties
 for the pointclasses $\Boraxi$ in arbitrary recursively presented Polish spaces. 
We apply our existence results to get better computations of the topological complexity of some sets of dictionaries considered in \cite{Lecomte05}.
\end{abstract} 

\noi {\bf  Keywords.}~$\omega$-power, Borel class, complete, 
effective descriptive set theory, hyperarithmetical hierarchy. 

\hs {\bf 2000 Mathematics Subject Classification.}~Primary: 03E15, 03B70, 
Secondary: 54H05, 68Q15, 68Q45, 68R15.\bigskip

\section{$\!\!\!\!\!\!$ Introduction.}

  We consider the finite alphabet $\Sigma\! =\!\{0,\ldots ,\Sigma\! -\! 1\}$, 
where $\Sigma\!\geq\! 2$ is an integer, and a language over this alphabet, i.e., 
a subset $A$ of the set $\Sigma^{<\omega}$ of finite words with letters in 
$\Sigma$. Notice that   a language of finite words will be also sometimes called   a dictionary, as in \cite{Lecomte05}. 
The set of infinite words over the alphabet $\Si$, i.e., of sequences of length $\om$ of letters of   $\Si$, is denoted $\Si^\om$. 

\begin{defi} The $\omega\! -\! power$ associated with $A$ is the set 
$A^\infty$ of the infinite sentences constructible with $A$ by concatenation. So 
we have ${A^\infty\! :=\!\{\ a_{0}a_{1}\ldots\!\in\!\Sigma^\omega\mid\forall i\! 
\in\!\omega~\ a_{i}\!\in\! A\ \}}$.\end{defi}

Notice that we denote here $A^\infty$ the $\omega$-power associated with $A$, as in \cite{Lecomte05}, while it is often denoted 
$A^\om$ in Theoretical Computer Science papers, as in \cite{sta,Fin01a,Fin03a,csl07}. Here we reserved the notation $A^\om$ to denote  the 
cartesian product of countably many copies of $A$ since this will be often used in this paper.\bigskip 

 In the theory of formal languages of infinite words, accepted by various kinds of automata, 
the $\omega$-powers appear very naturally  in the characterization of the class 
$REG_\om$  of  \orl s (respectively,  of the class $CF_\om$ of context free \ol s) 
 as the $\om$-Kleene closure 
of the family $REG$ of regular finitary languages (respectively,   of the 
family $CF$ of context free finitary languages),  see \cite{tho,lt,pp,stac,sta,mscs06,kms} for some references on this topic.\bigskip 

Since the set $\Sio$ of infinite words over a finite alphabet $\Si$ can be   equipped 
with the usual Cantor topology, the question of  the topological  complexity of  $\om$-powers of 
finitary  languages naturally arises and has  been posted by 
Niwinski \cite{Niwinski90},  Simonnet \cite{Simonnet92},  and  Staiger \cite{sta}.  

\hs 

\centerline{What are the possible levels of topological complexity for the 
$\omega$-powers?}\bigskip

As the concatenation map, from $A^\omega$ onto $A^\infty$, which associates $a_{0}a_{1}\ldots$ to $(a_{i})_{i\in\omega}$, 
is continuous, an $\omega$-power is always an analytic set.\bigskip 

  It has been recently  proved,   that 
for each integer $n\geq 1$, there exist some $\om$-powers of  (context-free) languages 
which are ${\bf \Pi}_n^0$-complete Borel sets,  \cite{Fin01a},  and   that there exists a 
(context-free) language $L$ such that $L^\om$ is analytic but not Borel, \cite{Fin03a}. Amazingly, 
the language $L$ is  very simple to describe and it is accepted by a simple $1$-counter automaton. Notice that 
Louveau has proved independently 
that analytic-complete $\omega$-powers exist, but the existence was proved in a non effective way.  
We refer the reader to \cite{HopcroftUllman79,ABB96} for basic notions about context-free languages.\bigskip

  The first author proved in \cite{Fin04} that  there exists a finitary language $V$ 
such that $V^\om$ is a Borel set of infinite rank. However  the only known fact 
on their complexity is that there is a (context-free) language $W$ such that $W^\om$ is Borel above 
${\bf \Delta_\omega^0}$, \cite{DF06}. In particular, it was still  unknown which could be  the possible infinite Borel ranks of $\om$-powers.\bigskip 

 The basic notions of descriptive set theory used in this paper will be 
recalled in the next section. We now state our results which extend the previous ones.

\vfill\eject

\begin{thm} (a) Let $3\!\leq\!\xi\! <\!\omega_1$, and ${\bf\Gamma}\!\not=\!\check {\bf\Gamma}$ be a Wadge class closed under finite unions satisfying the inclusions $\borxi ({\bf\Gamma})\!\subseteq\! {\bf\Gamma}\! =\!\borthree\hbox{\it -PU}({\bf\Gamma} )\!\subseteq\!\bormxipo$. Then there is 
$A\!\subseteq\! 2^{<\omega}$ such that $A^\infty$ is $\bf\Gamma$-complete.\smallskip 

\noindent (b) Let $1\!\leq\!\xi\! <\!\omega_1$. Then there is $A\!\subseteq\! 2^{<\omega}$ 
such that $A^\infty$ is $\boraxi$-complete.\smallskip

\noindent (c) Let $1\!\leq\!\xi\! <\!\omega_1$. Then there is $A\!\subseteq\! 2^{<\omega}$ such that 
$A^\infty$ is $\bormxi$-complete.\smallskip

\noindent (d) Let $1\!\leq\!\xi\! <\!\omega_1$. Then there is 
$A\!\subseteq\! 2^{<\omega}$ such that $A^\infty$ is $\check D_2 (\boraxi )$-complete.\smallskip

\noindent (e) Let $3\!\leq\!\xi\! <\!\omega_1$ and $\omega\!\leq\!\eta\! <\!\omega_1$ be an indecomposable ordinal. Then there is $A\!\subseteq\! 2^{<\omega}$ such that $A^\infty$ is 
$\check D_\eta (\boraxi )$-complete.\end{thm}

 So we get a complete knowledge of the Borel classes $\bf\Gamma$ for which there is 
$A\!\subseteq\! 2^{<\omega}$ such that $A^\infty$ is $\bf\Gamma$-complete. Indeed, the only 
class $\borxi$ admitting a complete set is $\borone$. And  
$A\! :=\!\{ s\!\in\! 2^{<\omega}\mid 0\!\prec\! s\ \ \hbox{\rm or}\ \ 1^2\!\prec\! s\}$ implies that 
$A^\infty\! =\! 2^\omega\!\setminus\! N_{10}$ is a $\borone$-complete set.\bigskip

 In this context coming from theoretical computer science, it is natural to wonder 
whether these examples are effective. We answer positively. The reader should 
see \cite{Moschovakis} for basic notions of effective descriptive set theory. It is known that 
$B\!\subseteq\! 2^\omega$ is $\boraxi$-complete if and only if 
$B\!\in\!\boraxi\!\setminus\!\bormxi$ (see 22.10 in \cite{Kechris}). The effective version 
of Theorem 1.2 is the following:

\begin{thm} (1) Let $1\!\leq\!\xi\! <\!\omega_1^{CK}$.\smallskip 

\noindent (a) There is $A\!\subseteq\! 2^{<\omega}$ such that 
$A^\infty\!\in\!\Boraxi\!\setminus\!\bormxi$.\smallskip

\noindent (b) There is $A\!\subseteq\! 2^{<\omega}$ such that $A^\infty\!\in\!\Bormxi\!\setminus\!\boraxi$.\smallskip 

 Moreover, $A$ can be coded by a $\Borone$ subset of $\omega$.\smallskip

\noindent (2) Similarly, let $\beta\!\in\! 2^\omega$ and $1\!\leq\!\xi\! <\!\omega_1^{\beta}$.\smallskip 

\noindent (a) There is $A\!\subseteq\! 2^{<\omega}$ such that 
$A^\infty\!\in\!\Boraxi (\beta )\!\setminus\!\bormxi$.\smallskip 

\noindent (b) There is $A\!\subseteq\! 2^{<\omega}$ such that 
$A^\infty\!\in\!\Bormxi (\beta )\!\setminus\!\boraxi$.\smallskip

 Moreover, $A$ can be coded by a $\Borone (\beta )$ subset of $\omega$.\end{thm}
 
 To prove Theorem 1.2,  we use a theorem of Kuratowski which is a level by level version of a theorem of Lusin and Souslin 
stating that every Borel set $B \subseteq 2^\om$ is the image of a closed subset of the Baire space 
$\om^\om$ by a continuous bijection. This theorem of Lusin and Souslin 
had already been used by Arnold in \cite{Arnold83} to prove that every Borel subset of 
$\Sio$, for a finite alphabet $\Si$,  is accepted by a non-ambiguous finitely branching  transition system with B\"uchi acceptance   
condition and our  first idea was to code the behaviour of such a transition system. This way, in the general case, 
we can manage to construct an $\om$-power of the same complexity as $B$. We now state Kuratowski's Theorem  \cite{Kuratowski} (see Corollary 33.II.1):
 
\begin {thm} 
Let $\xi\!\geq\! 1$ be a countable ordinal, $X$ a zero-dimensional Polish space, and 
$B\!\in\!\bormxipo (X)$. Then there is $C\!\in\!\bormone (\omega^\omega )$ and a continuous bijection 
$f\! :\! C\!\rightarrow\! B$ such that $f^{-1}$ is $\boraxi$-measurable (i.e., $f[U]$ is $\boraxi (B)$ for each open subset $U$ of $C$).
\end{thm}

 To prove Theorem 1.3, we first prove an effective version of Theorem 1.4. It has the following consequence.

\begin{thm} 
Let $\xi\!\geq\! 1$ be a countable ordinal, and $B\!\in\!\Bormxipo (2^\omega )$. Then there is 
$C\!\in\!\Bormone (\omega^\omega )$, a partial function 
$f\! :\!\omega^\omega\!\rightarrow\! 2^\omega$, recursive on $C$, 
and a partial function $g\! :\! 2^\omega\!\rightarrow\!\omega^\omega$, 
$\Boraxi$-recursive on $B$, such that $f$ defines a bijection from $C$ onto $B$ 
and $g$ coincides with $f^{-1}$.
\end{thm}

 To prove Theorems 1.3 and 1.5, we prove some results of effective 
descriptive set theory that cannot be found in \cite{Moschovakis}. We prove that the pointclasses $\boraxi$ are, 
uniformly and in the codes, closed under taking sections at points in spaces of type at most $1$, 
substitutions of partial recursive functions, finite intersections and unions, 
$\exists^\omega$, among other things.\bigskip

 In \cite{Lecomte05}, the following question is asked. What is the topological complexity of the set of dictionaries whose associated $\omega$-power is of a given level of complexity? More specifically, let 
 $1\!\leq\!\xi\! <\!\omega_1$. The following 
${\bf\Sigma}^1_2(2^{2^{<\omega}})\!\setminus\! D_2(\boraone )$ sets are introduced:
$$\begin{array}{ll}
{\bf\Sigma}_{\xi}\!\!\!\! & :=\!\{ A\!\subseteq\! 2^{<\omega}\mid A^\infty\!\in\!\boraxi\}
\hbox{\rm ,}\cr & \cr
{\bf\Pi}_{\xi}\!\!\!\! & :=\!\{ A\!\subseteq\! 2^{<\omega}\mid A^\infty\!\in\!\bormxi\}\hbox{\rm ,}\cr & \cr
\ {\bf\Delta}\!\!\!\! & :=\!\{ A\!\subseteq\! 2^{<\omega}\mid A^\infty\!\in\!\borel\}\!=\!
\{ A\!\subseteq\! 2^{<\omega}\mid A^\infty\!\in\!\ca\}.
\end{array}$$
The proof of Theorem 1.3 gives some more informations about the complexity of these sets. 
We will prove, using a result by J. Saint Raymond, that ${\bf\Sigma}_{\xi}$ and ${\bf\Pi}_{\xi}$ are $\ca$-hard if 
$\xi\!\geq\! 3$, which is a much better approximation of their complexity than the one in \cite{Lecomte05}. 
The proof of this fact has the following consequence. Theorem 1.2 shows that the $\omega$-powers are quite general objects. 
On the other hand, we will prove another result showing that they are not arbitrary.\bigskip

\noindent\bf Notation.\rm\ Let $\bf\Gamma$ be a class having a universal set 
${\cal U}^{2^\omega}_{\bf\Gamma}\!\subseteq\! (2^\omega )^2$, and ${\bf\Gamma'}$ another class. We set 

$${\cal U}({\bf\Gamma},{\bf\Gamma'})\! :=\!\{\alpha\!\in\! 2^\omega\mid ({\cal U}^{2^\omega}_{\bf\Gamma})_\alpha\!\in\! {\bf\Gamma'}\}.$$
Let $X,Y$ be zero-dimensional Polish spaces and $A\!\subseteq\! X$, $B\!\subseteq\! Y$. 
We will use the following notation to denote the Wadge quasi-order:
$$(X,A)\!\leq_W\! (Y,B)\ \Leftrightarrow\ \exists f\! :\! X\!\rightarrow\! Y\ \hbox{\rm continuous\ with}\ 
A\! =\! f^{-1}(B).$$
We write $(X,A)\! <_W\! (Y,B)$ if $(X,A)\!\leq_W\! (Y,B)$ and $(Y,B)\!\not\leq_W\! (X,A)$.\bigskip

 The consequence we mentioned is the following. If we choose suitable universal sets, then the following inequalities hold:
$$\begin{array}{ll}
{\cal U}(\bormxi ,\boraxi )\! & \leq_W {\bf\Sigma}_{\xi} <_W\ {\cal U}(\ana ,\boraxi )\cr & \cr 
{\cal U}(\boraxi ,\bormxi )\! & \leq_W {\bf\Pi}_{\xi} <_W\ {\cal U}(\ana ,\bormxi )\cr & \cr 
{\cal U}(\ana ,\ca )\! & \not\leq_W {\bf\Delta}\ <_W\ {\cal U}(\ana ,\borel )\! =\! {\cal U}(\ana ,\ca ).
\end{array}$$
This means that the $\omega$-powers are analytic sets that do not behave like arbitrary analytic sets. 
This also means that there is a strong difference between the Borel levels on one side, 
and the level of analytic sets on the other side. Actually, our method to prove Theorem 1.3 is 
a method that works for the Borel levels, and it cannot be extended to the level of analytic sets, 
even if Theorem 1.3 can be extended to the level of analytic sets (see \cite{Fin03a}). Note that we will prove that ${\cal U}(\ana ,\Borel )$ is ${\bf\Pi}^1_2$-complete.

\vfill\eject

 This paper is organized as follows:\bigskip

\noindent $\bullet$ In section 2 we prove Theorem 1.2.\bigskip
 
\noindent $\bullet$ In section 3 we recall a few basic facts of effective descriptive set 
theory, and fix some notation. Then we prove the results of effective 
descriptive set theory that we need for the sequel. This is where the closure 
properties for the pointclasses $\boraxi$ are proved.\bigskip

\noindent $\bullet$ In section 4 we prove Theorem 1.5.\bigskip

\noindent $\bullet$ In section 5 we prove Theorem 1.3.\bigskip

\noindent $\bullet$ In section 6 we study the complexity of some sets of dictionaries.

\section{$\!\!\!\!\!\!$ Proof of Theorem 1.2.}

$\underline{\bf{Basic\ facts\ and\ notation.}}$\bigskip

 In descriptive set theory, we study the topological complexity of definable 
subsets of Polish spaces, i.e., of separable and completely metrizable 
topological spaces.\bigskip 

\noindent $\bullet$ The notation for the $Borel\ classes$ in metrizable 
spaces is as follows: $\boraone$ is the class of open sets, and if $\xi\!\geq\! 1$ is a 
countable ordinal, then $\bormxi$ is the class of complements of $\boraxi$ sets, 
$\boraxi$ is the class of countable unions of sets in 
$\bigcup_{1\leq\eta <\xi}\ \bormeta$, and $\borxi$ is the class $\boraxi\cap\bormxi$. The class of 
$Borel\ sets$ is 
$$\borel\! :=\!\bigcup_{1\leq\xi <\omega_1}\ \boraxi\! =\!
\bigcup_{1\leq\xi <\omega_1}\ \bormxi .$$
$\bullet$ The class of $analytic\ sets$ is the class $\ana$ of subsets of Polish spaces that are continuous images of Polish spaces. One can prove that if $X$ is a Polish space, then 
$A\!\subseteq\! X$ is analytic if and only if $A$ is the projection on $X$ of a closed subset of 
$X\!\times\!\omega^\omega$ (see 14.3 in \cite{Kechris}). Then we can define the $projective\ classes$ in Polish spaces as follows; if $n\!\geq\! 1$ is an integer, then $\pmn$ is the class of complements of 
$\pan$ sets, $\panpo (X)$ is the class of projections on $X$ of sets in 
$\pmn (X\!\times\!\omega^\omega )$, and $\pn$ is the class $\pan\cap\pmn$.\bigskip

\noindent $\bullet$ If $\bf\Gamma$ is a class of sets in Polish spaces and $X$ is a Polish space, then a set ${\cal U}^{X}_{\bf\Gamma}\!\in\! {\bf\Gamma}(2^\omega\!\times\! X)$ is $universal$ for ${\bf\Gamma}(X)$ if 
${\bf\Gamma}(X)\! =\!\{ ({\cal U}^{X}_{\bf\Gamma})_\alpha\mid\alpha\!\in\! 2^\omega\}$ (where 
$({\cal U}^{X}_{\bf\Gamma})_\alpha\! :=\!\{ x\!\in\! X\mid (\alpha ,x)\!\in\! {\cal U}^{X}_{\bf\Gamma}\}$). For example, there are universal sets for $\boraxi (X)$, $\bormxi (X)$, $\ana (X)$, $\ca (X)$ for any Polish space $X$ (see 22.3 and 26.1 in \cite{Kechris}).\bigskip

\noindent $\bullet$ Recall that a Polish space is 
$zero\! -\! dimensional$ if it has a basis consisting of $\borone$ sets. 
Typically, let $K$ be a countable set. If $K$ is equipped with the discrete 
topology and $s\!\in\! K^{<\omega}$, then 
$N_{s}\! :=\!\{\alpha\!\in\! K^\omega\mid s\!\prec\!\alpha\}$ is a basic $\borone$  
set of $K^\omega$ ($s\!\prec\!\alpha$ means that $s$ is a beginning of $\alpha$). The 
length of $\gamma\!\in\! K^{\leq\omega}$ is denoted $\vert\gamma\vert$. If $\gamma\!\in\! K^{\leq\omega}$ 
and $k\!\in\!\omega$, then $\gamma\restriction k$ is the beginning of length $k$ of $\gamma$. If 
$s\!\prec\!\alpha\! =\!\alpha (0)\alpha (1)...$, then $\alpha\! -\! s$ is the sequence 
$\alpha (\vert s\vert )\alpha (\vert s\vert\! +\! 1)...$\bigskip 

\noindent $\bullet$ If $\bf\Gamma$ is a class of sets in zero-dimensional Polish 
spaces, closed under continuous preimages, then a subset $A$ of $X$ is 
${\bf\Gamma}\! -\! hard$ if for each $A'\!\in\! {\bf\Gamma}(X')$ there is a continuous map 
$f\! :\! X'\!\rightarrow\! X$ with $A'\! =\! f^{-1}(A)$.

 If $A\!\in\! {\bf\Gamma}(X)$ is $\bf\Gamma$-hard, then we say that $A$ is 
${\bf\Gamma}\! -\! complete$. We say that $\bf\Gamma$ 
is a $W\! adge\ class$ if there is a $\bf\Gamma$-complete set. We denote 
$\check {\bf\Gamma}\! :=\!\{\neg A\mid A\!\in\! {\bf\Gamma}\}$. If 
${\bf\Gamma}\!\not=\!\check {\bf\Gamma}\!\subseteq\!\borel$ is a Wadge class, then 
$A$ is ${\bf\Gamma}$-complete if and only if 
$A\!\in\! {\bf\Gamma}\!\setminus\!\check {\bf\Gamma}$.\bigskip

\noindent $\bullet$ If $\bf\Gamma$ is a Wadge class, then 
$\borxi (\bf\Gamma)\!\subseteq\! {\bf\Gamma }$ means that $E\!\in\! {\bf\Gamma }(X)$ if 
$E\!\in\!\borxi (A)$ and $A\!\in\! {\bf\Gamma }(X)$.\bigskip

\noindent $\bullet$ If $I$ is a set and ${\bf\Gamma}$ is a class or a set, then 
$(x_i)_{i\in I}\!\subseteq\! {\bf\Gamma}$ means that $x_i\!\in\! {\bf\Gamma}$ for each $i\!\in\! I$.\bigskip

\noindent $\bullet$ We set  
$\borxi\hbox{\rm -PU}({\bf\Gamma} )\! :=\!\{\bigcup_{n\in\omega}\ A_n\cap P_n\mid 
(A_n)_{n\in\omega}\!\subseteq\! {\bf\Gamma}\ \ \hbox{\rm and}\ \ (P_n)_{n\in\omega}\!\subseteq\!\borxi\ \ \hbox{\rm partition}\}$ if $1\!\leq\!\xi\! <\!\omega_1$. One can prove that if 
$\borone\!\subseteq\! {\bf\Gamma}\!\not=\!\check {\bf\Gamma}\!\subseteq\!\borel$ is a Wadge class, then there is a bigger $1\!\leq\!\xi\! <\!\omega_1$ (the $\underline{level}$ of ${\bf\Gamma}$) such that 
${\bf\Gamma}\! =\!\borxi\hbox{\rm -PU}({\bf\Gamma} )$ (see \cite{LouStRay}).\bigskip

 If $\eta\! <\!\omega_1$ and $(A_\theta )_{\theta <\eta}$ is an increasing sequence of subsets of some space $X$, then we set 
$$D_\eta [(A_\theta )_{\theta <\eta}]\! :=\!\{ x\!\in\! X\mid\exists\theta\! <\!\eta~\ \ 
x\!\in\! A_\theta\!\setminus\!\bigcup_{\theta'<\theta}\ A_{\theta'}\ \ \hbox{\rm and\ the\ parity\ of}\ \ \theta\ \ \hbox{\rm is opposite\ to\ that\ of}\ \ \eta\}.$$
If moreover $1\!\leq\!\xi\! <\!\omega_1$, then we set $D_\eta (\boraxi )\! :=\!\{ D_\eta 
[(A_\theta )_{\theta <\eta}]\mid (A_\theta )_{\theta <\eta}\!\subseteq\!\boraxi\}$. One can prove that 
$D_\eta (\boraxi )$ has level $\xi$ if $\eta\!\geq\! 1$ (see \cite{LouStRay}).\bigskip

\noindent $\bullet$ We say that $\omega\!\leq\!\eta\! <\!\omega_1$ is $indecomposable$ 
if $\eta$ cannot be represented as $\eta_1\! +\!\eta_2$ with $\eta_1,\eta_2\! <\!\eta$.
 It is known that the indecomposable ordinals are the $\omega^\theta$ with 
$1\!\leq\!\theta\! <\!\omega_1$ (see IV.2.16 in \cite{Levy}).\bigskip

\noindent $\underline{\bf{Proof \ of\ Theorem\ 1.2.}}$\bigskip

\noindent $\bullet$ We have already said  that the existence of the continuous bijection 
$f\! :\! C\!\rightarrow\! B$ 
given by Lusin and Souslin's Theorem had already  been used by Arnold in \cite{Arnold83} to prove that every Borel subset of 
$\Sio$, for a finite alphabet $\Si$,  is accepted by a non-ambiguous finitely branching  transition system with B\"uchi acceptance   
condition. We now recall the definition of these transition systems.\bigskip 

A $B\ddot{u}chi\ transition\ system$ is a tuple $\mathcal{T}=(\Si, Q, \delta, q_0, Q_f)$, where $\Si$ is a finite input alphabet, $Q$ is a countable set of states, 
$\delta \subseteq Q \times \Si  \times Q$ 
is the transition relation, $q_0 \in Q$ is the initial state, 
and $Q_f \subseteq Q$ is the set of final states. 
A run of $\mathcal{T}$ over an infinite word $\sigma \in \Sio$ is an infinite sequence 
of states $(t_i)_{i\geq 0}$, such that 
~ $t_0 = q_0$, ~and  for each $i\geq 0$,  ~ $(t_i, \sigma (i), t_{i+1}) \in \delta$. The run is said to be accepting iff there are infinitely 
many integers $i$ such that $t_i$ is in $Q_f$.\bigskip 

 The transition system is said to be $non\mbox{-}ambiguous$ if each infinite word $\sigma \in \Sio$ has at most 
one accepting run by $\mathcal{T}$.\bigskip

 The transition system is said to be $finitely\ branching$ if for each state $q\in Q$ and each $a \in \Si$, there are only finitely 
many states $q'$ such that $(q, a, q') \in \delta$.\bigskip 

 Our first idea was to code the behaviour of such a transition system. In fact this can be done on a part of infinite words of a special 
compact set $K_{0,0}$. However  we shall have also to consider more general sets $K_{N,j}$ and then we shall need the hypothesis of the 
$\boraxi$-measurability of the function $f$, which is given by Kuratowski's Theorem.
 
\vfill\eject 

\noindent $\bullet$ We now come to the proof of Theorem\ 1.2.\bigskip

\noindent (a) We may assume that $\borone\!\subseteq\! {\bf\Gamma}$, otherwise 
${\bf\Gamma}\! =\!\{\emptyset\}$ since $\borxi ({\bf\Gamma})\!\subseteq\! {\bf\Gamma}$, in which case $A\! :=\!\emptyset$ is suitable. This implies that $\borxi\!\subseteq\! {\bf\Gamma}$ since 
$\borxi ({\bf\Gamma})\!\subseteq\! {\bf\Gamma}$.\bigskip

\noindent $\bullet$ Let $B\!\in\! {\bf\Gamma}(2^\omega )\!\setminus\!\check {\bf\Gamma}$, and 
$P_\infty\! :=\!\{\alpha\!\in\! 2^\omega\mid\forall m\!\in\!\omega\ \exists n\!\geq\! m\ \ \alpha (n)\! =\! 1\}$, 
which is homeomorphic to $\omega^\omega$ (we associate $0^{\beta (0)}10^{\beta (1)}1...$ to 
$\beta\!\in\!\omega^\omega$). As 
$B\!\in\!\bormxipo$, Theorem 1.4 gives $C\!\in\!\bormone (P_\infty)$ and $f$. 
By Proposition 11 in \cite{Lecomte05}, it is enough to find $A\!\subseteq\! 4^{<\omega}$. 
The dictionary $A$ will be made of two pieces: we will 
have $A\! =\!\mu\cup\pi$. The set $\pi$ will code $f$, and $\pi^\infty$ will look like $B$ on some nice compact sets $K_{N,j}$. Outside a countable family of 
compact sets, we will hide $f$, so that $A^\infty$ will be the simple set $\mu^\infty$.\bigskip  

\noindent $\bullet$ We set $Q\! :=\!\{ (s,t)\!\in\! 2^{<\omega}\!\times\! 
2^{<\omega}\mid \vert s\vert\! =\!\vert t\vert\}$. We enumerate $Q$ as follows. We start 
with $q_{0}\! :=\! (\emptyset ,\emptyset )$. Then we put the sequences of length $1$ 
of elements of $2\!\times\! 2$, in the lexicographical ordering: 
$q_{1}\! :=\! (0,0)$, $q_{2}\! :=\! (0,1)$, $q_{3}\! :=\! (1,0)$, 
$q_{4}\! :=\! (1,1)$. Then we put the $16$ sequences of length $2$: 
$q_{5}\! :=\! (0^{2},0^{2})$, $q_{6}\! :=\! (0^{2},01)$, $\ldots$ And so on. We 
will sometimes use the coordinates of $q_{N}\! :=\! (q^{0}_{N},q^{1}_{N})$. We 
put $M_{j}\! :=\!\Sigma_{i<j}\ 4^{i+1}$. Note that the sequence 
$(M_j)_{j\in\omega}$ is strictly increasing, and that $q_{M_{j}}$ is the last 
sequence of length $j$ of elements of $2\!\times\! 2$.\bigskip

\noindent $\bullet$ If $l\!\in\!\omega$ and $(a_{i})_{i<l}\!\in\! (\omega^{<\omega} )^l$, 
then ${^\frown}_{i<l}\ a_{i}$ is the concatenation $a_{0}\ldots a_{l-1}$. Similarly, 
${^\frown}_{i\in\omega}\ a_{i}$ is the concatenation $a_0a_1\ldots$\bigskip
 
\noindent $\bullet$ Now we define the ``nice compact sets". We will sometimes view $2$ as 
an alphabet, and sometimes view it as a letter. To make this distinction clear, 
we will use the boldface notation $\bf 2$ for the letter, and the lightface notation $2$ otherwise.\bigskip

 We will have the same distinction with $3$ instead of $2$, so that $2=\{0, 1\}, 3= \{0, 1, {\bf 2}\}, 
4=\{0, 1, {\bf 2}, {\bf 3}\}$. Let $N,j$ be non-negative integers with $N\!\leq\! M_{j}$. We set 
$$K_{N,j}:=\{\ \gamma\! =\! {\bf 2}^{N}\ {^\frown}\  [\ {^\frown}_{i\in\omega}\ \ m_{i}\ {\bf 2}^{M_{j+i+1}}\ {\bf 3}\ 
{\bf 2}^{M_{j+i+1}}\ ]\!\in\! 4^\omega\mid\forall i\!\in\!\omega\ \ m_i\!\in\! 2\ \}.$$
As the map $\varphi_{N,j}\! :\! K_{N,j}\!\rightarrow\! 2^\omega$ defined by 
$\varphi_{N,j}(\gamma )\! :=\! (m_{i})_{i\in\omega}$ is a homeomorphism, $K_{N,j}$ is 
compact.\bigskip

\noindent $\bullet$ Now we will define the sets that ``look like $B$".\bigskip

\noindent - We define a function $c\! :\! B\!\times\!\omega\!\rightarrow\! Q$ by 
$c(\alpha ,l)\! :=\! [f^{-1}(\alpha ),\alpha ]\restriction l$. Note that $Q$ is countable, so that we equip it with the discrete topology. 
In these conditions, we prove that $c$ is $\boraxi$-measurable.\bigskip

 For any $q\in Q$, it holds that 
$c^{-1}(\{ q\})\! =\!\{ (\alpha ,l)\!\in\! B\!\times\!\omega\mid f^{-1}(\alpha )\restriction l\! =\! q^0\ \ \hbox{\rm and}\ \ \alpha\restriction l\! =\! q^1 \}$. 
But  $\alpha \restriction l = q^1$ means that ``$l\! =\!\vert q^1\vert$ and $\alpha$ belongs to the basic open set $N_{q^1}$". 
In the same fashion, $f^{-1}(\alpha ) \restriction l = q^0$ means that 
``$l\! =\!\vert q^0\vert$ and $f^{-1}(\alpha ) $ belongs to the basic open set $N_{q^0}$", or equivalently that 
``$l\! =\!\vert q^0\vert$ and $\alpha = f( f^{-1}(\alpha ) )$ belongs to 
$f[C\cap N_{q^0}]$". As $f[C\cap N_{q^0}]$ is a   $\boraxi$ subset of $B$, $c^{-1}(\{ q\})$ is a  
$\boraxi$ subset of $B\!\times\!\omega$ and $c$ is $\boraxi$-measurable.\bigskip 

\noindent - Let $N$ be an integer. We put 

$$E_N\! :=\!\{\ \alpha\!\in\! 2^\omega\mid q^{1}_{N}\alpha\!\in\! B\ \ 
\hbox{\rm and}\ \ c(q^{1}_{N}\alpha ,\vert q^{1}_{N}\vert )\! =\! q_{N}\ \}.$$ 
Note that $E_0\! =\!\{ \ \alpha\!\in\! 2^\omega\mid\alpha\!\in\! B~\ \hbox{\rm and}~\ 
c(\alpha ,0)\! =\!\emptyset\}\! =\! B$. Let us prove that $E_N\!\in\! {\bf\Gamma}(2^\omega )$ for each integer $N$.

\vfill\eject

 As $c$ is $\boraxi$-measurable and $\{ q_N\}\!\in\!\borone (Q)$, we get 
$c^{-1}(\{ q_N\})\!\in\!\borxi (B\!\times\!\omega)$. Note that the map 
$S\! :\!\{\alpha\!\in\! 2^\omega\mid q^{1}_{N}\alpha\!\in\! B\}\!\rightarrow\! B\!\times\!\omega$ defined by  
$S(\alpha )\! :=\! (q^{1}_{N}\alpha ,\vert q^{1}_{N}\vert )$ is continuous, so that $E_N\! =\! S^{-1}[c^{-1}(\{ q_N\})]$ 
is in $\borxi (\{\alpha\!\in\! 2^\omega\mid q^{1}_{N}\alpha\!\in\! B\})$. As 
$B\!\in\!\ {\bf\Gamma}(2^\omega )$ and the map $\alpha\!\mapsto\! q^{1}_{N}\alpha$ is continuous, 
$\{\alpha\!\in\! 2^\omega\mid q^{1}_{N}\alpha\!\in\! B\}$ is in ${\bf\Gamma}(2^\omega )$. Thus  
$E_N\!\in\! {\bf\Gamma}(2^\omega )$ since $\!\borxi ({\bf\Gamma} )\!\subseteq\! {\bf\Gamma}$.\bigskip

Now we define the transition system obtained from $f$.\bigskip  

\noindent - If $m\!\in\! 2$ and $n,p\!\in\!\omega$, then we write $n\buildrel m\over\rightarrow p$ if 
$q^{0}_{n}\!\prec\! q^{0}_{p}$ and $q^{1}_{p}\! =\! q^{1}_{n}m$.\bigskip  

\noindent  - As $f$ is continuous on $C$, the graph $\hbox{\rm Gr}(f)$ of $f$ is a closed subset of 
$C\!\times\! 2^\omega$. As $C$ is $\bormone (P_\infty )$, $\hbox{\rm Gr}(f)$ is also a closed subset of $P_\infty\!\times\! 2^\omega$. 
So there is a closed subset $F$ of $2^\omega\!\times\! 2^\omega$ such that $\hbox{\rm Gr}(f)\! =\! F\cap (P_\infty\!\times\! 2^\omega )$. 
We identify 
$2^\omega\!\times\! 2^\omega$ with $(2\!\times\! 2)^\omega$, i.e., we view 
$(\beta ,\alpha )$ as $[\beta (0),\alpha (0)],[\beta (1),\alpha (1)],...$ By Proposition 2.4 in \cite{Kechris}, there is 
$R\!\subseteq\! (2\!\times\! 2)^{<\omega}$, closed under initial segments, such that 
$F\! =\!\{ (\beta ,\alpha)\!\in\! 2^\omega\!\times\! 2^\omega\mid\forall k\!\in\!\omega\ \ 
(\beta ,\alpha )\restriction k\!\in\! R\}$. Notice that $R$ is a tree whose infinite branches form the set $F$. In particular, we get 

$$(\beta ,\alpha )\!\in\!\hbox{\rm Gr}(f)\ \Leftrightarrow\ \beta\!\in\! P_\infty\ \ \hbox{\rm and}\ \ \forall k\!\in\!\omega\ \ (\beta,\alpha)\restriction k\!\in\! R.$$
- We set $Q_{f}\! :=\{ (t,s)\!\in\! R\mid t\!\not=\!\emptyset\ \ \hbox{\rm and}\ \ t(\vert t\vert\! -\! 1)\! =\!1\}$. Notice that $Q_f$ is simply the set of pairs 
$(t,s)\!\in\! R$ such that the last letter of $t$ is a $1$.\bigskip  

 We have in fact already defined the transition system $\mathcal{T}$ obtained from $f$. This transition system has a countably  infinite set $Q$ of states and a set $Q_{f}$ of accepting states. The initial state is $q_{0}\! :=\! (\emptyset ,\emptyset )$. The input alphabet is $2=\{0, 1\}$ and the transition relation 
$\delta \subseteq Q \times 2   \times Q$ is given by: if $m\!\in\! 2$ and $n,p\!\in\!\omega$ then 
$(q_n, m, q_p) \in \delta$ iff     $n\buildrel m\over\rightarrow p$.\bigskip

 Recall that a run $(t_i)_{i\geq 0}$ of $\mathcal{T}$ is said to be B\"uchi accepting if there are infinitely 
many integers $i$ such that $t_i$ is in $Q_f$. 
Then the set of $\omega$-words over the alphabet $2$ which are accepted by the transition system
 $\mathcal{T}$ from the initial state $q_0$ with B\"uchi acceptance condition is exactly the Borel set $B$.\bigskip

\noindent $\bullet$ Now we define the finitary language $\pi$. We set
$$\pi\! :=\!\left\{ 
\begin{array}{ll}
& \!\!\!\!\!\! s\!\in\! 4^{<\omega}\mid\exists j,l\!\in\!\omega\ \ 
\exists (m_i)_{i\leq l}\!\in\! 2^{l+1}\ \ 
\exists (n_i)_{i\leq l}, (p_i)_{i\leq l}, 
(r_i)_{i\leq l}\!\in\!\omega^{l+1}\cr 
& \cr
& \ \ \ \ \ \ \ \ \ \ \ \ \ \ \ \ \ \ \ \ \ \ \ \ \ \ \ \ \ \ \ \ \ \ \ 
\ \ \ \ \ \ \ \ \ \ \ \ \ \ \ \ \ \ \ \ \ \ \ \ \ \ \ \ \ \ \ \ \ \ 
\begin{array}{ll}
& \!\!\!\!\!\! n_{0}\!\leq\! M_j\cr & \hbox{\rm and}\cr
& \!\!\!\!\!\!\forall i\!\leq\! l\ \ n_i\buildrel 
{m_i}\over\rightarrow p_i\ \ \hbox{\rm and}\ \ 
p_i\! +\! r_i = M_{j+i+1}\cr & \hbox{\rm and}\cr
& \!\!\!\!\!\!\forall i\! <\! l\ \ p_i = n_{i+1}\cr & \hbox{\rm and}\cr
& \!\!\!\!\!\! q_{p_l}\!\in\! Q_f\cr & \hbox{\rm and}\cr
& \!\!\!\!\!\! 
s = {^\frown}_{i\leq l}\ \ {\bf 2}^{n_i}\ m_i\ {\bf 2}^{p_i}\ {\bf 2}^{r_i}\ {\bf 3}\ {\bf 2}^{r_i}
\end{array}\!\!\!
\end{array}\right\}.$$
$\bullet$ Let us prove that $\varphi_{N,j}[\pi^\infty\cap K_{N,j}]\! =\! E_N$ if 
$N\!\leq\! M_j$.

\vfill\eject

 Let $\gamma\!\in\!\pi^\infty\cap K_{N,j}$, and 
$\alpha\! :=\!\varphi_{N,j}(\gamma )$. We can write 
$$\gamma = {^\frown }_{k\in\omega}\ [\ {^\frown }_{i\leq l_{k}}\ \ {\bf 2}^{n^k_i}\ 
m^k_i\ {\bf 2}^{p^k_i}\ {\bf 2}^{r^k_i}\ {\bf 3}\ {\bf 2}^{r^k_i}\ ]\hbox{\rm .}$$
As this decomposition of $\gamma$ is in $\pi$, we have 
$n^k_i\buildrel {m^k_i}\over\rightarrow p^k_i$ if $i\!\leq\! l_{k}$, 
$p^k_{i}\! =\! n^k_{i+1}$ if $i\! <\! l_{k}$, and $q_{p^k_{l_{k}}}\!\in\! Q_{f}$, 
for each $k\!\in\!\omega$. Moreover, $p^k_{l_{k}}\! =\! n^{k+1}_{0}$, for each 
$k\!\in\!\omega$, since $\gamma\!\in\! K_{N,j}$ implies that 
$p^k_{l_{k}}\! +\! r^k_{l_{k}}\! =\! r^k_{l_{k}}\! +\! n^{k+1}_{0}\! =\! M_{j+1+m}$ for some integer $m$. So we get 
$$N\buildrel {\alpha (0)}\over\rightarrow p^{0}_0
\buildrel {\alpha (1)}\over\rightarrow\ldots
\buildrel {\alpha (l_{0})}\over\rightarrow p^{0}_{l_{0}}
\buildrel {\alpha (l_{0}+1)}\over\rightarrow p^{1}_{0}
\buildrel {\alpha (l_{0}+2)}\over\rightarrow\ldots
\buildrel {\alpha (l_{0}+l_{1}+1)}\over\rightarrow p^{1}_{l_{1}}\ldots$$
In particular we have 
$$q^{0}_{N}\prec q^{0}_{p^{0}_0}\prec\ldots\prec q^{0}_{p^{0}_{l_{0}}}\prec 
q^{0}_{p^{1}_0}\prec\ldots\prec q^{0}_{p^{1}_{l_{1}}}\ldots$$
because $n\buildrel {m}\over\rightarrow p$ implies that $q_n^0\prec q_p^0$. Note that 
$\vert q^{1}_{p^{k}_{l_{k}}}\vert\! =\!\vert q^{1}_{N}\vert\! +\!
\Sigma_{j\leq k}\ (l_{j}\! +\! 1)$, so that the sequence 
$(\vert q^{0}_{p^{k}_{l_{k}}}\vert)_{k\in\omega}$ is strictly increasing since $|q_n^0|\! =\! |q_n^1|$ for each integer $n$. This implies the existence of $\beta\!\in\! P_{\infty}$ such that $q^{0}_{p^{k}_{l_{k}}}\prec\beta$ for each $k\!\in\!\omega$. Note that $\beta\!\in\! P_{\infty}$ because,  for each integer $k$, $q_{p^{k}_{l_{k}}}\!\in\! Q_f$. Note also that 
$(\beta ,q^1_N\alpha )\restriction k\!\in\! R$ for infinitely many $k$'s. As $R$ is closed under initial segments, $(\beta ,q^1_N\alpha )\restriction k\!\in\! R$ for every $k\!\in\!\omega$, so that 
$q^1_N\alpha\! =\! f(\beta )\!\in\! B$. Moreover, 
$$c(q^{1}_{N}\alpha ,\vert q^{1}_{N}\vert )\! =\! 
(\beta\restriction\vert q^{1}_{N}\vert ,q^{1}_{N})\! =\! 
(q^{0}_{N},q^{1}_{N})\! =\! q_{N}\hbox{\rm ,}$$ 
and $\alpha\!\in\! E_{N}$.\bigskip

 Conversely, let $\alpha\!\in\! E_{N}$. We have to see that 
$\gamma\! :=\!\varphi_{N,j}^{-1}(\alpha )\!\in\!\pi^\infty$. As $\gamma\!\in\! K_{N,j}$, 
we are allowed to write $\gamma = {\bf 2}^{N}\ {^\frown}\ [\ {^\frown}_{i\in\omega}\ \ \alpha (i)\ 
{\bf 2}^{M_{j+i+1}}\ {\bf 3}\ ^{M_{j+i+1}}\ ]$. We set 
$\beta\! :=\! f^{-1}(q^{1}_{N}\alpha )$. There is a sequence of integers $(k_{l})_{l\in\omega}$ such that 
${q_{k_{l}}\! =\! (\beta ,q^{1}_{N}\alpha )\restriction l}$. Note that 
$N\buildrel {\alpha (0)}\over\rightarrow k_{\vert q^{1}_{N}\vert +1}
\buildrel {\alpha (1)}\over\rightarrow k_{\vert q^{1}_{N}\vert +2}\ldots$ 
As $N\!\leq\! M_{j}$ we get $k_{\vert q^{1}_{N}\vert +i+1}\!\leq\! M_{j+i+1}$. So 
we can define $n_{0}\! :=\! N$, 
$p_{0}\! :=\! k_{\vert q^{1}_{N}\vert +1}$, 
$r_{0}\! :=\! M_{j+1}\! -\! p_{0}$, $n_{1}\! :=\! p_{0}$. Similarly, we can 
define $p_{1}\! :=\! k_{\vert q^{1}_{N}\vert +2}$, 
$r_{1}\! :=\! M_{j+2}\! -\! p_{1}$. We go on like this until we find some 
$q_{p_{i}}$ in $Q_{f}$. This clearly defines a word in $\pi$. And we can go 
on like this, so that $\gamma\!\in\!\pi^\infty$.\bigskip

 Thus $\pi^\infty\cap K_{N,j}$ is in 
${\bf\Gamma}(K_{N,j})\!\subseteq\! {\bf\Gamma}(4^\omega )$. Notice that we proved, among other things, the equality $\varphi_{0,0}[\pi^\infty\cap K_{0,0}]\! =\! B$. In 
particular, $\pi^\infty\cap K_{0,0}$ is ${\bf\Gamma}$-complete in $K_{0,0}$.\bigskip

 Notice that $\pi^\infty$ codes on $K_{0,0}$ the behaviour of the transition system accepting $B$. In a similar way, $\pi^\infty$ codes on $K_{N,j}$ the behaviour of the same transition system, but starting this time from the state $q_N$ instead of the initial state $q_0$. But some $\omega$-words in $\pi^\infty$ are not in $K_{0,0}$ and not even in any $K_{N,j}$ and we do not know what exactly the complexity of this set of $\omega$-words is. However we remark that all the words in $\pi$ have the same form 
${\bf 2}^N\ {^\frown}\ [\ {^\frown}_{i\leq l}\ \ m_i\ {\bf 2}^{P_i}\ {\bf 3}\ {\bf 2}^{R_i}\ ]$.\bigskip

\noindent $\bullet$ We are ready to define $\mu$. The idea is that an infinite sequence 
containing a word in $\mu$ cannot be in the union of the $K_{N,j}$'s.

\vfill\eject

 We set

$$\begin{array}{ll}
\mu^{0} & \!\!\!\! :=\!\left\{ 
\begin{array}{ll}
& \!\!\!\!\!\! s\!\in\! 4^{<\omega}\mid\ \exists l\!\in\!\omega\ \ 
\exists (m_i)_{i\leq l+1}\!\in\! 2^{l+2}\ \ \exists N\!\in\!\omega
\ \ \exists (P_i)_{i\leq l+1}, (R_i)_{i\leq l+1}\!\in\!\omega^{l+2}\cr 
& \cr
& \ \ \ \ \ \ \ \ \ \ \ \ \ \ \ \ \ \ \ \ \ \ \ \ \ \ \ \ \ \ \ \ \ \ \ 
\ \ \ \ \ \ \ \ \ \ \ \ \ \ \ \ \ \ \ \ \ \ \ \ \ \ \ \ \ \ \ \ \ \ 
\begin{array}{ll}
& \!\!\!\!\!\!\forall i\!\leq\! l\! +\! 1\ \ \exists j\!\in\!\omega\ \  P_i\! =\!  M_{j}\cr & \hbox{\rm and}\cr
& \!\!\!\!\!\!  P_l\!\not=\! R_{l}\cr & \hbox{\rm and}\cr
& \!\!\!\!\!\! s = {\bf 2}^N\ {^\frown}\ [\ {^\frown}_{i\leq l+1}\ \ m_i\ {\bf 2}^{P_i}\ {\bf 3}\ {\bf 2}^{R_i}\ ]
\end{array}\!\!\!
\end{array}\right\}\hbox{\rm ,}\cr & \cr 
\mu^{1} & \!\!\!\! :=\!\left\{ 
\begin{array}{ll}
& \!\!\!\!\!\! s\!\in\! 4^{<\omega}\mid\ \exists l\!\in\!\omega\ \ 
\exists (m_i)_{i\leq l+1}\!\in\! 2^{l+2}\ \ \exists N\!\in\!\omega
\ \ \exists (P_i)_{i\leq l+1}, (R_i)_{i\leq l+1}\!\in\!\omega^{l+2}\cr 
& \cr
& \ \ \ \ \ \ \ \ \ \ \ \ \ \ \ \ \ \ \ \ \ \ \ \ \ \ \ \ \ \ \ \ \ \ \ 
\ \ \ \ \ \ \ \ \ \ \ \ \ \ \ \ \ \ \ \ \ \ \ \ \ \ \ \ \ \ \ \ \ \ 
\begin{array}{ll}
& \!\!\!\!\!\!\forall i\!\leq\! l\! +\! 1\ \ \exists j\!\in\!\omega\ \  P_i\! =\!  M_{j}\cr & \hbox{\rm and}\cr
& \!\!\!\!\!\!\exists j\!\in\!\omega\ \  (P_l\!=\! M_{j}\ \ \hbox{\rm and}\ \ P_{l+1}\!\not=\! M_{j+1})\cr & \hbox{\rm and}\cr
& \!\!\!\!\!\! s = {\bf 2}^N\ {^\frown}\ [\ {^\frown}_{i\leq l+1}\ \ m_i\ {\bf 2}^{P_i}\ {\bf 3}\ {\bf 2}^{R_i}\ ]
\end{array}\!\!\!
\end{array}\right\}\hbox{\rm ,}\cr & \cr 
\mu & \!\!\!\! :=\!\mu^{0}\cup\mu^{1}.
\end{array}$$
All the words in $A$ have the same form 
${\bf 2}^N\ {^\frown}\ [\ {^\frown}_{i\leq l}\ \ m_i\ {\bf 2}^{P_i}\ {\bf 3}\ {\bf 2}^{R_i}\ ]$. Note that any finite concatenation of words of this form still has this form. Moreover, such a concatenation is in $\mu^i$ if its last word is in $\mu^i$.\bigskip

\noindent $\bullet$ Now we prove that $\mu^\infty$ is ``simple". The previous remarks show that 
$$\mu^\infty\! =\!\{\ \gamma\!\in\! 4^\omega\mid\exists i\!\in\! 2\ \ \forall j\!\in\!\omega\ \ 
\exists k,n\!\in\!\omega\ \ \exists t\!\in\! (\mu^i)^{n+1}\ \ n\!\geq\! j
\ \ \hbox{\rm and}\ \ \gamma\restriction k\! =\! {^\frown}_{l\leq n}\ t(l)\ \}.$$
This shows that $\mu^\infty\!\in\!\bormtwo (4^\omega )$.\bigskip

\noindent $\bullet$ Note again that all words in $A$ have the same form 
${\bf 2}^N\ {^\frown}\ [\ {^\frown}_{i\leq l}\ \ m_i\ {\bf 2}^{P_i}\ {\bf 3}\ {\bf 2}^{R_i}\ ]$. We set\bigskip 

\leftline{$P\! :=\!\{ {\bf 2}^N\ {^\frown}\ [\ {^\frown}_{i\in\omega}\ \ m_i\ {\bf 2}^{P_i}\ {\bf 3}\ {\bf 2}^{R_i}\ ]\!\in\! 4^\omega\mid (m_i)_{i\in\omega}\!\in\! 2^\omega\hbox{\rm ,}\ N\!\in\!\omega\hbox{\rm ,}\ 
(P_i)_{i\in\omega}\hbox{\rm ,}\ (R_i)_{i\in\omega}\!\in\!\omega^\omega\ \ \hbox{\rm and}$}\medskip
\rightline {$\forall i\!\in\!\omega\ \exists j\!\in\!\omega\ P_i\! =\! M_j\}.$}\bigskip

\noindent We define a map 
$F\! :\! P\!\setminus\!\mu^\infty\!\rightarrow (\{\emptyset\}\cup\mu )\!\times\!\omega^2$ as follows. Let 
$\gamma\! :=\! {\bf 2}^N\ {^\frown}\ [\ {^\frown}_{i\in\omega}\ \ m_i\ {\bf 2}^{P_i}\ {\bf 3}\ {\bf 2}^{R_i}\ ]\!\in\! P\!\setminus\!\mu^\infty$, and $j_0\!\in\!\omega$ with $P_0\! =\! M_{j_0}$. If $\gamma\!\in\! K_{N,j_0-1}$, then we put $F(\gamma )\! :=\! (\emptyset ,N,j_0)$.  If $\gamma\!\notin\! K_{N,j_0-1}$, then there is an integer $l$ maximal for which $P_l\!\not=\! R_l$ or there is $j\!\in\!\omega$ with $P_l\! =\! M_j$ and 
$P_{l+1}\!\not=\! M_{j+1}$. Let $j_1\!\in\!\omega$ with $P_{l+2}\! =\! M_{j_1}$. We put
$$F(\gamma )\! :=\! ({\bf 2}^N\ {^\frown}\ [\ {^\frown}_{i\leq l}\ \ m_i\ {\bf 2}^{P_i}\ {\bf 3}\ {\bf 2}^{R_i}\ ]\ {^\frown}\ m_{l+1}\ {\bf 2}^{P_{l+1}}\ {\bf 3},R_{l+1},j_1).$$
$\bullet$ Fix $\gamma\!\in\! A^\infty$. If $\gamma\!\notin\!\mu^\infty$, then 
$\gamma\!\in\! P\!\setminus\!\mu^\infty$, $F(\gamma )\! :=\! (t,S,j)$ is defined. Note that 
$t\ {\bf 2}^S\!\prec\!\gamma$, and that $j\! >\! 0$. Moreover, 
$\gamma\! -\! t\ {\bf 2}^S\!\in\! K_{0,j-1}$. Note also that 
$S\!\leq\! M_{j-1}$ if $t\! =\!\emptyset$, and that 
$t\ {\bf 2}^S\ \gamma (\vert t\vert\! +\! S)\ {\bf 2}^{M_{j}}\ {\bf 3}\!\notin\!\mu$. Moreover, there is an integer 
$N\!\leq\!\hbox{\rm min}(M_{j-1},S)$ ($N\! =\! S$ if $t\! =\!\emptyset$) such that 
$\gamma\! -\! t\ {\bf 2}^{S-N}\!\in\!\pi^\infty\cap K_{N,j-1}$, since the last word in $\mu$ in the decomposition of $\gamma$ (if it exists) ends before $t\ {\bf 2}^S$.\bigskip 
 
\noindent $\bullet$ In the sequel we will say that $(t,S,j)\!\in\! (\{\emptyset\}\cup\mu )\times\omega^2$ is 
$suitable$ if $S\!\leq\! M_{j}$ if $t\! =\!\emptyset$,  $t(\vert t\vert\! -\! 1)\! =\! {\bf 3}$ if 
$t\!\in\!\mu$, and $t\ {\bf 2}^S\ m\ {\bf 2}^{M_{j+1}}\ {\bf 3}\!\notin\!\mu$ if $m\!\in\! 2$. We set, for $(t,S,j)$ suitable, 
$$P_{t,S,j}:=\left\{\ \gamma\!\in\! 4^\omega\mid t\ {\bf 2}^S\!\prec\!\gamma\ \ \hbox{\rm and}\ \ 
\gamma\! -\! t\ {\bf 2}^S\!\in\! K_{0,j}\ \right\}.$$
Note that $P_{t,S,j}$ is a compact subset of $P\!\setminus\!\mu^\infty$, and that 
$F(\gamma )\! =\! (t,S,j\! +\! 1)$ if $\gamma\!\in\! P_{t,S,j}$. This shows that the 
$P_{t,S,j}$'s, for $(t,S,j)$ suitable, are pairwise disjoint. Note also that $\mu^\infty$ 
is disjoint from $\bigcup_{(t,S,j)\ \hbox{\rm suitable}}\ P_{t,S,j}$.\bigskip
   
\noindent $\bullet$ We set, for $(t,S,j)$ suitable and $N\!\leq\!\hbox{\rm min}(M_{j},S)$ 
($N\! =\! S$ if $t\! =\!\emptyset$),    
$$A_{t,S,j,N}:=\left\{\ \gamma\!\in\! P_{t,S,j}\mid\gamma\! -\! t\ {\bf 2}^{S-N}\!\in\!\pi^\infty\cap K_{N,j}\ \right\}.$$ 
Note that $A_{t,S,j,N}\!\in\! {\bf\Gamma}(4^\omega )$ since $N\!\leq\! M_j$.\bigskip
  
\noindent $\bullet$ The previous discussion shows that 
$$A^\infty\! =\!\mu^\infty\cup
\bigcup_{(t,S,j)\ \hbox{\rm suitable}}
\bigcup_{
\begin{array}{ll}
& N\leq\hbox{\rm min}(M_j,S)\cr 
& \ N=S\ \hbox{\rm if}\ t=\emptyset
\end{array}}
\ A_{t,S,j,N}.$$ 
As $\bf\Gamma$ is closed under finite unions, the set 
$$A_{t,S,j}:=\!\bigcup_{
\begin{array}{ll}
& N\leq\hbox{\rm min}(M_j,S)\cr 
& \ N=S\ \hbox{\rm if}\ t=\emptyset
\end{array}}
\ A_{t,S,j,N}$$ 
is in ${\bf\Gamma}(4^\omega )$.\bigskip

\noindent $\bullet$ We can write
$$A^\infty\! =\!\mu^\infty\!\setminus\!\left(\bigcup_{(t,S,j)\ \hbox{\rm suitable}}\ P_{t,S,j}
\right)\ \cup\bigcup_{(t,S,j)\ \hbox{\rm suitable}}\ A_{t,S,j}\cap P_{t,S,j}.$$
Note that the $P_{t,S,j}$'s and $\bigcup_{(t,S,j)\ \hbox{\rm suitable}}\ P_{t,S,j}$ are $\borthree$ subsets 
of $4^\omega$ since $(P_{t,S,j})_{(t,S,j)\ \hbox{\rm suitable}}$ is a countable family of closed sets. 
Moreover, $\mu^\infty$ is a $\bormtwo\!\subseteq\! {\bf\Gamma}$ subset of $4^\omega$. 
This implies that $A^\infty$ is in $\borthree\hbox{\rm -PU}({\bf\Gamma})\! =\! {\bf\Gamma}$. Moreover, the set 
$A^\infty\cap P_{\emptyset ,0,0}\! =\!\pi^\infty\cap P_{\emptyset ,0,0}\! =\!
\pi^\infty\cap K_{0,0}$ is $\bf\Gamma$-complete. 
This shows that $A^\infty$ is $\bf\Gamma$-hard (any reduction with values in $K_{0,0}$ is also a reduction with values in $4^\omega$). 
Thus $A^\infty$ is $\bf\Gamma$-complete.\bigskip

 We can now end the proof of Theorem 1.2.\bigskip

\noindent (b) If $\xi\! =\! 1$, then we can take 
$A\! :=\!\{ s\!\in\! 2^{<\omega}\mid 0\!\prec\! s\ \ \hbox{\rm or}\ \ \exists k\!\in\!\omega\ \ 10^k1\!\prec\! s\}$ and $A^\infty\! =\! 2^\omega\!\setminus\!\{ 10^\infty\}$ is $\boraone$-complete.\bigskip

\noindent $\bullet$ If $\xi\! =\! 2$, then we will see in Theorem 2 the existence of $A\!\subseteq\! 2^{<\omega}$ such that $A^\infty$ is $\boratwo$-complete.\bigskip

\noindent $\bullet$  So we may assume that $\xi\!\geq\! 3$, and we just have to apply (a) to 
${\bf\Gamma}\! :=\!\boraxi$.

\vfill\eject 

\noindent (c) If $\xi\! =\! 1$, then we can take $A\! :=\!\{ 0\}$ and $A^\infty\! =\!\{ 0^\infty\}$ is 
$\bormone$-complete.\bigskip

\noindent $\bullet$ If $\xi\! =\! 2$, then we can take $A\! :=\!\{ 0^k1\mid k\!\in\!\omega\}$ and 
$A^\infty\! =\! P_\infty$ is $\bormtwo$-complete.\bigskip

\noindent $\bullet$ So we may assume that $\xi\!\geq\! 3$, and we just have to apply (a) to ${\bf\Gamma}\! :=\!\bormxi$.\bigskip 

\noindent (d) First notice that $D_2(\boraxi )\! =\!\{ B\!\setminus\! C\mid B,C\!\in\!\boraxi\}$. Indeed, 
$\subseteq$ is clear, and $\supseteq$ comes from the fact that 
$B\!\setminus\! C\! =\! (B\cup C)\!\setminus\! C$. This implies that 
$\check D_2(\boraxi )\! =\!\{ B\cup C\mid B\!\in\!\boraxi\ \hbox{\rm and}\ C\!\in\!\bormxi\}$. A consequence of this is the closure of $\check D_2(\boraxi )$ under finite unions. Another consequence is  
$\borxi [\check D_2(\boraxi )]\!\subseteq\!\check D_2(\boraxi )$. Indeed, if 
$D\! :=\! B\cup C\!\in\!\check D_2(\boraxi )(X)$ and $E\!\in\!\borxi (D)$, then choose 
${\bf\Sigma}\!\in\!\boraxi (X)$ and ${\bf\Pi}\!\in\!\bormxi (X)$ such that 
$E\! =\! {\bf\Sigma}\cap D\! =\! {\bf\Pi}\cap D$. We get 
$E\! =\! ({\bf\Sigma}\cap B)\cup ({\bf\Pi}\cap C)\!\in\!\check D_2(\boraxi )(X)$.\bigskip

\noindent $\bullet$ If $\xi\! =\! 1$, then we can take $A\! :=\!\{s\in 2^{<\omega}\mid 0\!\prec\! s~\ \hbox{\rm or}~\ \exists q\!\in\!\omega ~{(101)^q}1^3\!\prec\! s~\ \hbox{\rm or}~\ s\! =\! 10^2\}$ and $A^\infty\! =\!\bigcup_{p\in\omega}~[N_{{(10^2)^p} 0}\cup (\bigcup_{q\in\omega}~N_{{(10^2)^p} {(101)^q} 1^3})]\cup 
\{(10^2)^\infty\}$ is $\check D_2(\boraone )$-complete (see $\S$7 in \cite{Lecomte05}, and also example 9 in \cite{Staiger97b}).\bigskip

\noindent $\bullet$ If $\xi\! =\! 2$, then we can take $A\! :=\!\{ s\!\in\! 2^{<\omega}\mid 1^2\!\prec\! s~\ \hbox{\rm or}~\ s\! =\! 0\}$ and 
$$A^\infty\! =\! \left(\{0^\infty\}\cup\bigcup_{p\in\omega}\ N_{0^p1^2}\right)\cap [(2^\omega\!\setminus\! P_{\infty})
\cup\{\alpha\!\in\! 2^\omega\mid\forall m\!\in\!\omega~\exists n\!\geq\! m~\ \alpha (n)\! =\! \alpha 
(n\! +\! 1)\! =\! 1\}]$$ 
is $\check D_2(\boratwo )$-complete (see $\S$7 in \cite{Lecomte05}).\bigskip

\noindent $\bullet$ So we may assume that $\xi\!\geq\! 3$, and we just have to apply (a) to ${\bf\Gamma}\! :=\!\check D_2(\boraxi )$.\bigskip

\noindent (e) Let $X$ be a zero-dimensional Polish space, and $E,F\!\in\! D_\eta(\boraxi )(X)$. By Lemma 4.2 in \cite{eng}, $E\!\times\! F$ is $D_\eta(\boraxi )$. Now let ${\cal C}\!\subseteq\! 2^\omega$ be $D_\eta(\boraxi )$-complete, $h\! :\! 2^\omega\!\times\! 2^\omega\!\rightarrow\! 2^\omega$ continuous with ${\cal C}\!\times\! {\cal C}\! =\! h^{-1}({\cal C})$, and 
$f,g\! :\! X\!\rightarrow\! 2^\omega$ continuous with $E\! =\! f^{-1}({\cal C})$ and 
$F\! =\! g^{-1}({\cal C})$. It is clear that the map $c\! :\! X\!\rightarrow\! 2^\omega$ defined by 
$c(x)\! :=\! h[f(x),g(x)]$ satisfies $E\cap F\! =\! c^{-1}({\cal C})$. This shows that $D_\eta(\boraxi )$ is closed under finite intersections. Thus $\check D_\eta(\boraxi )$ is closed under finite unions.\bigskip

 Note also that if $D\!\in\! D_\eta(\boraxi )$ and $B\!\in\!\boraxi$, then $B\cup D\!\in\! D_\eta(\boraxi )$. Indeed, let $(A_\theta )_{\theta <\eta}\!\subseteq\!\boraxi$ be an increasing sequence with 
$D\! =\! D[(A_\theta )_{\theta <\eta}]$. We set $B_0\! :=\!\emptyset$, $B_1\! :=\! B$, and 
$B_{2+\theta}\! :=\! A_\theta\cup B$ if $\theta\! <\!\eta$. Then 
$(B_\theta )_{\theta <\eta}\!\subseteq\!\boraxi$ is increasing, and 
$D[(B_\theta )_{\theta <\eta}]\! =\! B\cup\bigcup_{2\rho +1<\eta}\ (A_{2\rho +1}\cup B)\!\setminus\! 
(A_{2\rho}\cup B)\! =\! B\cup D$ since $\eta$ is even. This shows that if $D\!\in\!\check D_\eta(\boraxi )$ and $B\!\in\!\bormxi$, then $B\cap D\!\in\!\check D_\eta(\boraxi )$. This implies the inclusion 
$\borxi [\check D_\eta(\boraxi )]\!\subseteq\!\check D_\eta(\boraxi )$.\bigskip

 Now we can apply (a) to ${\bf\Gamma}\! :=\!\check D_\eta(\boraxi )$.\hfill{$\square$}\bigskip
 
 As we already said, a Borel class remains for which we did not provide a complete $\omega$-power yet: the class $\boratwo$. Note that it is easy to see that  the classical example of a $\boratwo$-complete set, the set $2^\omega\!\setminus\! P_\infty$,  is not an $\omega$-power. However we are going to prove the following result.\bigskip 

\noindent\bf Theorem 2\it\ \ There is a  recursive (and even context-free) language 
$A\!\subseteq\! 2^{<\omega}$ such that $A^\infty\!\in\!\Boratwo\!\setminus\!\bormtwo$.\rm\bigskip

\noindent {\bf Proof.} By Proposition 11 in \cite{Lecomte05}, it is enough to find $A\!\subseteq\! 3^{<\omega}$. We set, for $j\! <\! 3$ and $s\!\in\! 3^{<\omega}$,
$$\begin{array}{ll}
n_j(s)\!\! & :=\ \hbox{\rm Card}\{ i\! <\!\vert s\vert\mid s(i)\! =\! j\}\hbox{\rm ,}\cr & \cr
\ \ \ \ \ \ T\!\! & :=\ \{\alpha\!\in\! 3^{\leq\omega}\mid\forall l\! <\! 1\! +\!\vert\alpha\vert\ \  
n_2(\alpha\restriction l)\!\leq\! n_1(\alpha\restriction l)\}.
\end{array}$$
So $T$ is the tree of sequences for which any initial segment contains more coordinates equal to 1 than  coordinates equal to 2.\bigskip

\noindent $\bullet$ We inductively define, for $s\!\in\! T\cap 3^{<\omega}$, $s^{\hookleftarrow}\!\in\! 2^{<\omega}$ as follows:
$$s^{\hookleftarrow}\! :=\!\left\{\!\!\!\!\!\!\begin{array}{ll}
& \emptyset\ \ \hbox{\rm if}\ \ s\! =\!\emptyset\hbox{\rm ,}\cr & \cr
& t^{\hookleftarrow}\varepsilon\ \ \hbox{\rm if}\ \ s\! =\! t\varepsilon\ \ \hbox{\rm and}\ \ 
\varepsilon\! <\! 2\hbox{\rm ,}\cr & \cr
& t^{\hookleftarrow}\hbox{\rm ,\ except\ that\ its\ last\ 1\ is\ replaced\ with\ 0,\ if}\ s\! =\! t{\bf 2}.
\end{array}\right.$$
$\bullet$ We will extend this definition to infinite sequences. To do this, we introduce a notion of limit. Fix $(s_n)_{n\in\omega}\!\subseteq\! 2^{<\omega}$. We define 
${\displaystyle\lim_{n\rightarrow\infty}{s_n}}\!\in\! 2^{\leq\omega}$ as follows. For each 
$t\!\in\! 2^{<\omega}$,
$$t\!\prec\! {\displaystyle\lim_{n\rightarrow\infty}{s_n}}\ \Leftrightarrow\ \exists n_0\!\in\!\omega\ \ 
\forall n\!\geq\! n_0\ \ t\!\prec\! s_n.$$
$\bullet$ If $\alpha\!\in\! T\cap 3^{\omega}$, then we set 
$\alpha^{\hookleftarrow}\! :=\! {\displaystyle\lim_{n\rightarrow\infty}{(\alpha\restriction n)^{\hookleftarrow}}}$. We define $e\! :\! T\cap 3^{\omega}\!\rightarrow\! 2^\omega$ by 
$e(\alpha )\! :=\!\alpha^{\hookleftarrow}$. Note that $T\cap 3^{\omega}\!\in\!\Bormone (3^{\omega})$, and $e$ is a $\Boratwo$-recursive partial function on $T\cap 3^{\omega}$, since for 
$t\!\in\! 2^{<\omega}$ we have 
$$t\!\prec\! e(\alpha )\ \Leftrightarrow\ \exists n_0\!\in\!\omega\ \ \forall n\!\geq\! n_0\ \ t\!\prec\! 
(\alpha\restriction n)^{\hookleftarrow}.$$
$\bullet$ We set $E\! :=\!\{s\!\in\! T\cap 3^{<\omega}\mid n_2(s)\! =\! n_1(s)\ \ \hbox{\rm and}\ \ 
s\!\not=\!\emptyset\ \ \hbox{\rm and}\ \ 1\!\prec\! [s\restriction (\vert s\vert\! -\! 1)]^{\hookleftarrow}\}$. Note that 
$\emptyset\!\not=\! s^{\hookleftarrow}\!\prec\! 0^\infty$, and that $s(\vert s\vert\! -\! 1)\! =\! {\bf 2}$ changes $s(0)\! =\! [s\restriction (\vert s\vert\! -\! 1)]^{\hookleftarrow}(0)\! =\! 1$ into $0$ if $s\!\in\! E$.\bigskip

\noindent $\bullet$ If $S\!\subseteq\! 3^{<\omega}$, then $S^*\! :=\!\{ {^\frown}_{i<l}\ s_i\!\in\! 3^{<\omega}\mid l\!\in\!\omega\ \ \hbox{\rm and}\ \ (s_i)_{i<l}\!\subseteq\! S\}$. We put 
$$A\! :=\!\{ 0\}\cup E\cup\{ {^\frown}_{j\leq k}\ (c_j1)\!\in\! 3^{<\omega}\mid
[\forall j\!\leq\! k\ \ c_j\!\in\! (\{ 0\}\cup E)^*]\ \ 
\hbox{\rm and}\ \ [k\! >\! 0\ \ \hbox{\rm or}\ \ (k\! =\! 0\ \ \hbox{\rm and}\ \ c_0\!\not=\!\emptyset )]\}.$$ 
Note that $A$ is recursive.\bigskip

\noindent $\bullet$ In the proof of Theorem 1.2.(b) we met the set 
$\{ s\!\in\! 2^{<\omega}\mid 0\!\prec\! s\ \ \hbox{\rm or}\ \ \exists k\!\in\!\omega\ \ 10^k1\!\prec\! s\}$. We will call this set $B$, and $B^\infty\! =\! 2^\omega\!\setminus\!\{ 10^\infty\}$ is $\boraone$-complete (and even $\Boraone$). Let us show that $A^\infty\! =\! e^{-1}(B^\infty )$.\bigskip

\noindent - By induction on $\vert t\vert$, we get $(st)^{\hookleftarrow }={s^{\hookleftarrow }}
{t^{\hookleftarrow }}$ if $s,t\!\in\! T\cap 3^{<\omega}$. Let us show that 
$(s\beta)^{\hookleftarrow }\! =\! {s^{\hookleftarrow }}{\beta ^{\hookleftarrow }}$ if moreover 
$\beta\!\in\! T\cap 3^{\omega }$.\bigskip 

 Assume that $t\!\prec\! (s\beta )^{\hookleftarrow }$. Then there is $m_{0}\!\geq\!\vert s\vert$ such that, for 
$m\geq m_{0}$, 
$$t\!\prec\! [(s\beta )\restriction m]^{\hookleftarrow }\! =\! [s\beta\restriction (m\! -\!\vert s\vert )]^{\hookleftarrow }\! =\! 
{s^{\hookleftarrow }}[\beta\restriction (m\! -\!\vert s\vert )]^{\hookleftarrow }.$$ 
This implies that $t\prec {s^{\hookleftarrow }}{\beta ^{\hookleftarrow }}$ if 
$\vert t\vert\! <\!\vert s^{\hookleftarrow }\vert $. If $\vert t\vert\!\geq\!\vert s^{\hookleftarrow }\vert$, then there is $m_{1}\!\in\!\omega$ such that, for $m\!\geq\! m_{1}$, $\beta ^{\hookleftarrow }\restriction 
(\vert t\vert\! -\!\vert s^{\hookleftarrow }\vert )\!\prec\! [\beta\restriction (m\! -\!\vert s\vert )]^{\hookleftarrow }$. Here again, we get $t\!\prec\! {s^{\hookleftarrow }} {\beta ^{\hookleftarrow }}$. Thus 
$(s\beta)^{\hookleftarrow }\! =\! {s^{\hookleftarrow }}{\beta ^{\hookleftarrow }}$.\bigskip 

 Let $(s_{i})_{i\in\omega}\!\subseteq\! T\cap 3^{<\omega}$. Then ${^\frown}_{i\in\omega}\ s_{i}\!\in\! T$, and 
$({^\frown}_{i\in\omega}\ s_{i})^{\hookleftarrow }\! =\! {^\frown}_{i\in\omega}\ s_{i}^{\hookleftarrow }$, by the previous facts.\bigskip

\noindent - Let $(a_{i})_{i\in\omega}\!\in\! (A\!\setminus\!\{\emptyset\})^\omega$ and 
$\alpha\! :=\! {^\frown}_{i\in\omega}\ a_i$. As $A\!\subseteq\! T$, $e(\alpha )\! =\! ({^\frown}_{i\in\omega}\ a_{i})^{\hookleftarrow }\! =\! {^\frown}_{i\in\omega}\ a_{i}^{\hookleftarrow }$.\bigskip

 If $a_0\!\in\!\{ 0\}\cup E$, then $\emptyset\!\not=\! a_0^{\hookleftarrow}\!\prec\!0^\infty$, thus 
$e(\alpha )\!\in\! N_0\!\subseteq\! 2^\omega\!\setminus\!\{ 10^\infty\}\! =\! B^\infty$.\bigskip

 If $a_0\!\notin\!\{ 0\}\cup E$, then $a_0\! =\! {^\frown}_{j\leq k}\ (c_j1)$, thus 
$a_0^{\hookleftarrow}\! =\! {^\frown}_{j\leq k}\ (c_j^{\hookleftarrow}1)$.\bigskip

\ \ \ If $c_0\!\not=\!\emptyset$, then $e(\alpha )\!\in\! B^\infty$ as before.\bigskip

\ \ \ If $c_0\! =\!\emptyset$, then $k\! >\! 0$, so that 
$e(\alpha )\!\not=\! 10^\infty$ since $e(\alpha )$ has at least two coordinates equal to $1$.\bigskip

 We proved that $A^\infty\!\subseteq\! e^{-1}(B^\infty )$.\bigskip

\noindent - Assume that $e(\alpha )\!\in\! B^\infty$. We have to find 
$(a_i)_{i\in\omega}\!\subseteq\! A\!\setminus\!\{\emptyset\}$ with 
$\alpha\! =\! {^\frown}_{i\in\omega}\ a_i$. We split into cases:\bigskip

\noindent 1. $e(\alpha )\! =\! 0^\infty$.\bigskip

\noindent 1.1. $\alpha (0)\! =\! 0$.\bigskip

 In this case $\alpha\! -\! 0\!\in\! T$ and $e(\alpha\! -\! 0)\! =\! 0^\infty$. Moreover, $0\!\in\! A$. We put 
$a_0\! :=\! 0$.\bigskip

\noindent 1.2. $\alpha (0)\! =\! 1$.\bigskip

 In this case there is a coordinate $j_0$ of $\alpha$ equal to ${\bf 2}$ ensuring that $\alpha (0)$ is replaced with a $0$ in $e(\alpha )$. We put $a_0\! :=\!\alpha\restriction (j_0\! +\! 1)$, so that  
 $a_0\!\in\! E\!\subseteq\! A$, $\alpha\! -\! a_0\!\in\! T$ and $e(\alpha\! -\! a_0)\! =\! 0^\infty$.\bigskip
 
 Now the iteration of the cases 1.1 and 1.2 shows that $\alpha\!\in\! A^\infty$.\bigskip
 
\noindent 2. $e(\alpha )\! =\! 0^{k+1}10^\infty$ for some $k\!\in\!\omega$.\bigskip

 As in case 1, there is $c_0\!\in\! (\{ 0\}\cup E)^*$ such that 
$c_0\!\prec\!\alpha$, $c_0^{\hookleftarrow}\! =\! 0^{k+1}$, 
$\alpha\! -\! c_0\!\in\! T$ and $e(\alpha\! -\! c_0)\! =\! 10^\infty$. Note that 
$\alpha (\vert c_0\vert )\! =\! 1$, $\alpha\! -\! (c_01)\!\in\! T$ and $e[\alpha\! -\! (c_01)]\! =\! 0^\infty$. We put 
$a_0\! :=\! c_01$, and argue as in case 1.\bigskip
 
\noindent 3. $e(\alpha )\! =\! ({^\frown}_{j\leq l+1}\ 0^{k_j}1)0^\infty$ for some $l\!\in\!\omega$.\bigskip

 The previous cases show the existence of $(c_j)_{j\leq l+1}\!\subseteq\! (\{ 0\}\cup E)^*$ such that 
$a_0\! :=\! {^\frown}_{j\leq l+1}\ c_j1\!\prec\!\alpha$, $\alpha\! -\! a_0\!\in\! T$ and 
$e(\alpha\! -\! a_0)\! =\! 0^\infty$. We are done since $a_0\!\in\! A$.\bigskip
 
\noindent 4. $e(\alpha )\! =\! {^\frown}_{j\in\omega}\ 0^{k_j}1$.\bigskip

 An iteration of the discussion of case 3 shows that we can take $a_i$ of the form 
${^\frown}_{j\leq l+1}\ c_j1$.\bigskip

\noindent $\bullet$ The previous discussion shows that $A^\infty\! =\! e^{-1}(B^\infty )$. As $e$ is $\Boratwo$-recursive, $e^{-1}(B^\infty)\!\in\!\Boratwo (3^{\omega})$.\bigskip 

 It remains to see that $e^{-1}(B^\infty)\!\notin\!\bormtwo$. We argue by contradiction. We know that  $B^\infty\! =\! 2^\omega\!\setminus\!\{ 10^\infty\}$, so 
$e^{-1}(\{10^\infty\})\! =\! (T\cap 3^\omega)\!\setminus\! e^{-1}(B^\infty )$ is a $\boratwo$ subset of $ 3^\omega$ since $T\cap 3^{\omega}$ is closed in $3^{\omega}$. Thus $e^{-1}(\{10^\infty\})$ is  a countable union of compact subsets of $3^{\omega}$.\bigskip

 Consider now the cartesian product $(\{ 0\}\cup E)^{\omega}$ of countably many copies of 
$\{ 0\}\cup E$. The set $\{ 0\}\cup E$ is countable and it can be equipped with the discrete 
topology. The product  $(\{ 0\}\cup E)^{\omega}$ is equipped with the product topology 
of the discrete topology on $\{ 0\}\cup E$. In these conditions, the topological space 
$(\{ 0\}\cup E)^{\omega}$ is homeomorphic to the Baire space $\omega^\omega$.\bigskip 

 Consider now the map  $h\! :\! (\{ 0\}\cup E)^{\omega}\!\rightarrow\! e^{-1}(\{10^\infty\})$ defined by 
$h(\gamma )\! :=\! 1[{^\frown}_{i\in\omega}\ \gamma_i]$ for each sequence 
$\gamma\! =\! (\gamma_0, \gamma_1, \ldots )\!\in\! (\{ 0\}\cup E)^{\omega}$. 
The map $h$  is a homeomorphism by the previous discussion. As 
$(\{ 0\}\cup E)^{\omega}$ is homeomorphic to $\omega^\omega$, 
the Baire space $\omega^\omega$ is also homeomorphic to $e^{-1}(\{10^\infty\})$. This implies that $\omega^\omega$ is a countable union of compact sets. But this is absurd, by Theorem 7.10 in \cite{Kechris}.\bigskip

\noindent $\bullet$ It remains to see that $A$ is context-free. We assume here that the reader is familiar with the theory 
of formal languages and of context-free languages; basic notions  may be found in the Handbook Chapter \cite{ABB96}.\bigskip 

 It is easy to see that the language $E$ is in fact accepted by a $1$-counter automaton: it is the set of words 
$s\!\in\!  3^{<\omega}$ such that 
$$\forall 1\!\leq\! l\! <\!\vert s \vert\ \  
n_2(s \restriction l)\!< \! n_1(s \restriction l)\ \ \hbox{\rm and}\ \ n_2(s)\! =\! n_1(s)\ \ \hbox{\rm and}\ \ s(0)\! =\! 1\ \ 
\hbox{\rm and}\ \ s(\vert s\vert\! -\! 1)\! =\! {\bf 2}.$$ 
 This implies that $A$ is also accepted by a $1$-counter automaton because the class of $1$-counter languages is closed under concatenation and 
star operation. In particular $A$ is a context-free language because the class of languages accepted by  $1$-counter automata form a strict subclass 
of the class of context-free languages.\hfill{$\square$}\bigskip

\noindent\bf Remark.\rm\ The operation $\alpha\!\rightarrow\!\alpha^{\hookleftarrow}$ we have defined 
 is very close to the erasing operation defined by J. Duparc in his study of the Wadge hierarchy (see \cite{Duparc01}). 
However we have  modified this operation in such a way that $\alpha^{\hookleftarrow}$ is always  infinite 
when $\alpha$ is infinite, and that it has the good property with regard to $\omega$-powers and topological complexity.\bigskip
 
\noindent\bf Question.\rm\ What are the Wadge classes $\bf\Gamma$ for which there is 
$A\!\subseteq\! 2^{<\omega}$ such that $A^\infty$ is $\bf\Gamma$-complete? 
We have seen that Theorem 1.2 solves completely the case where $\bf\Gamma$ is a Borel class, and it also solves the problem for some other Wadge classes.
 The problem is solved for a few other Wadge classes in \cite{Lecomte01,Lecomte05}. 
We do not know (yet?) any Wadge class for which this problem cannot be solved.

\section{$\!\!\!\!\!\!$ Effective descriptive set theory background.}

$\underline{\bf{Basic\ facts\ and\ notation.}}$\bigskip

\noindent $\bullet$ In \cite{Moschovakis}, the classical $arithmetical\ hierarchy$ is defined 
as follows (see 3E). Let $X$ be a recursively presented Polish space,  
$[N(X,k)]_{k\in\omega}$ an effective enumeration of a neighborhood basis for the 
topology of $X$, and $B\!\subseteq\! X$. We say that $B\!\in\!\Boraone (X)$ if 
there is a recursive map $\varepsilon\! :\!\omega\!\rightarrow\!\omega$ such that 
$B\! =\!\bigcup_{i\in\omega}\ N[X,\varepsilon (i)]$. If $n\!\geq\! 1$ is an integer, then 
$\Bormn$ is the class of complements of $\Boran$ sets. We say that 
$B\!\in\!\Boranpo$ if there is $C\!\in\!\Bormn (\omega\!\times\! X)$ such that 
$B\! =\!\exists^\omega C\! :=\!\{ x\!\in\! X\mid\exists i\!\in\!\omega\ (i,x)\!\in\! C\}$. We also set 
$\Born\! :=\!\Boran\cap\Bormn$.\bigskip

\noindent $\bullet$ We say that $\gamma\!\in\!\Boraone$ if 
$\{ k\!\in\!\omega\mid\gamma\!\in\! N(\omega^\omega ,k)\}\!\in\!\Boraone (\omega)$. Let 
$\beta\!\in\! 2^\omega$. The $relativization$  
$\Boraone (\beta )$ of $\Boraone$ to $\beta$ is defined as follows. A set $P\!\subseteq\! X$ is in 
$\Boraone (\beta )$ if there is $Q\!\in\!\Boraone (2^\omega\!\times\! X)$ such that $P\! =\! Q_\beta$. As  before we say that $\gamma\!\in\!\Boraone (\beta )$ if 
$\{ k\!\in\!\omega\mid\gamma\!\in\! N(\omega^\omega ,k)\}\!\in\!\Boraone (\beta )(\omega)$.\bigskip

\noindent $\bullet$ Recall the existence of a $good\ parametrization$ 
in $\Boran$ for $\boran$ (see 3E.2, 3F.6 and 3H.1 in \cite{Moschovakis}). This means that 
there is a system of sets 
$G^{\Boran ,X}\!\in\!\Boran (\omega^\omega\!\times\! X)$ such that for each 
recursively presented Polish space $X$ and for each $P\!\subseteq\! X$, 
$$\begin{array}{ll}
P\!\in\!\boran & \!\!\Leftrightarrow\  \ 
\exists\gamma\!\in\!\omega^\omega\ \ P\! =\! G^{\Boran ,X}_{\gamma}
\hbox{\it ,}\cr & \cr
P\!\in\!\Boran & \!\!\Leftrightarrow \ \ 
\exists\gamma\!\in\!\Boraone\ \ P\! =\! G^{\Boran ,X}_{\gamma}.
\end{array}$$
Moreover, if $X$ is a recursively presented Polish space of type at most 1 
(i.e., a finite product of spaces equal to $\omega$, $\omega^\omega$ or 
$2^\omega$), and $Y$ is a recursively presented Polish space, then there is 
$S^{X,Y}_{\Boran}\! :\!\omega^\omega\!\times\! X\!\rightarrow\!\omega^\omega$ recursive so 
that 
$$(\gamma ,x,y)\!\in\! G^{\Boran ,X\times Y}\ \ \Leftrightarrow\ \ \  
[S^{X,Y}_{\Boran}(\gamma ,x),y]\!\in\! G^{\Boran ,Y}.$$
Note that $G^{\Boran ,X}$ is universal for $\boran (X)$ (with $\omega^\omega$ instead of $2^\omega$).\bigskip

\noindent $\bullet$ Let $f\! :\! X\!\rightarrow\! Y$ be a partial function, $D\!\subseteq\!\hbox{Domain}(f)$ and $P\!\subseteq\! X\!\times\!\omega$. Then $P\ computes\ f\ on\ D$ if 
$$x\!\in\! D\ \ \Rightarrow\ \ \forall k\!\in\!\omega\ \ [f(x)\!\in\! N(Y,k)\ \Leftrightarrow\ (x,k)\!\in\! P].$$
If $P$ is in some pointclass $\it\Gamma$ and computes $f$ on $D$, then we say that $f$ is 
$\it\Gamma\! -\! recursive\ on\ D$. This means that $f^{-1}[N(Y,k)]$ is in $\it\Gamma$, uniformly in $k$. We also say $recursive\ on\ D$ for $\Boraone$-recursive on $D$.\bigskip

\noindent $\bullet$ We also recall the notation for the coding of partial 
recursive functions from $X$ into $Y$ introduced in \cite{Moschovakis} (see 7A). We first define 
a partial function $U\! :\!\omega^\omega\!\times\! X\!\rightarrow\! Y$ by 
$$\begin{array}{ll}
U(\gamma ,x)\!\downarrow 
& \!\Leftrightarrow\ \ U(\gamma ,x)\ \hbox{\rm is\ defined}\ 
\Leftrightarrow\ \exists y\!\in\! Y\ \forall 
k\!\in\!\omega\ \ [y\!\in\! N(Y,k)\Leftrightarrow 
(\gamma ,x,k)\!\in\!G^{\Boraone ,X\times\omega}]\hbox{\rm ,}\cr & \cr U(\gamma ,x) 
& \!\!\! :=\ \ \hbox{\rm the\ unique}\ y\!\in\! Y\ \hbox{\rm such that}\ \forall k\!\in\!\omega\ \ 
[y\!\in\! N(Y,k)\Leftrightarrow (\gamma ,x,k)\!\in\! G^{\Boraone ,X\times\omega}].
\end{array}$$
Now let $\gamma\!\in\!\omega^\omega$. The function 
$\{\gamma\}^{X,Y}\! :\! X\!\rightarrow\! Y$ is defined by $\{\gamma\}^{X,Y}(x)\! :=\! U(\gamma ,x)$. 
Then a partial function $f\! :\! X\!\rightarrow\! Y$ is recursive on its domain if and only if there is 
$\gamma\!\in\!\Boraone$ such that $f(x)\! =\!\{\gamma\}^{X,Y}(x)$ when $f(x)$ is defined. More generally, the functions of the form $\{\gamma\}^{X,Y}$ are the partial continuous functions from a subset of $X$ into $Y$. We will write $\{\gamma\}$ instead of $\{\gamma\}^{X,Y}$ when $Y\! =\!\omega^\omega$, in order to simplify the notation.\bigskip

 If $X$ is of type at most $1$ and $Z$ is a recursively presented Polish space, then there is 
a recursive map $S^{X,Y,Z}_{\Boraone}\! :\!\omega^\omega\!\times\! X\!\rightarrow\!\omega^\omega$ 
such that $\{\gamma\}^{X\times Y,Z}(x,y)\! =\!\{ S^{X,Y,Z}_{\Boraone}(\gamma ,x)\}^{Y,Z}(y)$ if 
$(\gamma ,x)\!\in\!\omega^\omega\!\times\! X$.\bigskip

 Kleene's Recursion Theorem asserts that if $f\! :\!\omega^\omega\!\times\! X\!\rightarrow\! Y$ is 
recursive on its domain, then there is $\varepsilon^*\!\in\!\Boraone$ such that 
$f(\varepsilon^*,x)\! =\!\{\varepsilon^*\}^{X,Y}(x)$ when $f(\varepsilon^*,x)$ is defined (see 7A.2 in \cite{Moschovakis}). This will be the fundamental tool in the sequel. It is very useful to prove effective versions of classical results.\bigskip

\noindent $\bullet$ We will use the following basic maps:\bigskip

\noindent - We first define a one-to-one map 
$<.>:\!\omega^{<\omega}\!\rightarrow\!\omega$. Let $(p_n)_{n\in\omega}$ be the 
sequence of prime numbers. We set $<\!\emptyset\! >:=\! 1$, and, if 
$t\! :=\! (t_0,...,t_{l})\!\in\!\omega^{l+1}$, then we set 
$\overline{t}\! :=<t_0,...,t_{l}>:=\! p_0^{t_0+1}...p_{l}^{t_{l}+1}$.\bigskip

\noindent - If $k\!\in\!\omega$, then we say that ``$\hbox{\rm Seq}(k)$" (i.e., ``$k$ is a sequence") if 
$k\! =<t_0,...,t_{l-1}>$ for some $t_0,...,t_{l-1}$.\bigskip

\noindent - The length $\hbox{\rm lh}(k)$ of $k\!\in\!\omega$ is $l$ if $\hbox{\rm Seq}(k)$ and 
$k\! =<t_0,...,t_{l-1}>$, $0$ otherwise.\bigskip

\noindent - If $k,i\!\in\!\omega$, then we define $(k)_i\! :=\! t_i$ if 
$\hbox{\rm Seq}(k)$, $k\! =<t_0,...,t_{l-1}>$ and $i\! <\! l$, $0$ otherwise.\bigskip

\noindent - If $\gamma\!\in\!\omega^\omega$ and $i\!\in\!\omega$, then we define 
$(\gamma )_i\!\in\!\omega^\omega$ by $(\gamma )_i(j)\! :=\!\gamma ({\bf <i,j>})$. But 
here we do not use the injection $(i,j)\!\mapsto\ <i,j>$ above, since we want a 
bijection from $\omega^{2}$ into $\omega$. So we use the notation ${\bf <i,j>}$ 
for $2^i\cdot (2j\! +\! 1)\! -\! 1$, when $(\gamma )_i$ is concerned. 
The inverse bijection is denoted $s\!\mapsto\! {\bf [(s)_{0},(s)_{1}]}$.\bigskip

\noindent $\underline{\bf{Borel\ codes\ and\ closure\ properties.}}$\bigskip

\noindent\bf Notation.\rm\ We give a coding of Borel sets slightly different from 
the one given in \cite{Moschovakis} (see 7B), since there is a problem for $\boraone$. It can be found in some unpublished 
notes written by  Louveau, \cite{Louveau}. We define by induction on the countable ordinal $\xi\!\geq\! 1$ the set $BC_{\xi}$ of Borel codes for $\boraxi$ as follows. 
If $\gamma\!\in\!\omega^\omega$, then we define 
$\gamma^*\!\in\!\omega^\omega$ by $\gamma^*(i)\! :=\!\gamma (i\! +\! 1)$. We set 
$$\begin{array}{ll}
BC_{1} & \!\!\! :=\{\ \gamma\!\in\!\omega^\omega\mid\gamma (0)\! =\! 0\ \}
\hbox{\rm ,}\cr & \cr
BC_{\xi} & \!\!\! :=\Big\{\ \gamma\!\in\!\omega^\omega\mid\gamma (0)\! =\! 1\ \ 
\hbox{\rm and}\ \ \forall i\!\in\!\omega\ \ \{\gamma^{*}\}(i)\!\downarrow\ 
\hbox{\rm and}\ \ \{\gamma^{*}\}(i)\!\in\!\bigcup_{1\leq\eta <\xi}\ BC_{\eta}\ \Big
\}\ \hbox{\rm if}\ \xi\! \geq\! 2.\end{array}$$
The set of Borel codes is $BC\! :=\!\bigcup_{1\leq\xi <\omega_{1}}\ BC_{\xi}$. We 
also set $BC^{*}\! :=\!\bigcup_{2\leq\xi <\omega_{1}}\uparrow\ BC_{\xi}$.\bigskip

 Now let $X$ be a recursively presented Polish space. We define 
$\rho^{X}\! :\! BC\!\rightarrow\!\borel (X)$ by induction:
$$\rho^{X}(\gamma ):=\left\{\!\!\!\!\!\begin{array}{ll}
& \!\!\!\bigcup_{i\in\omega}\ N[X,\gamma^{*}(i)]\ \ \hbox{\rm if}\ \ 
\gamma\!\in\! BC_{1}\hbox{\rm ,}\cr & \cr
& \!\!\!\bigcup_{i\in\omega}\ X\!\setminus\rho^{X}[\{\gamma^{*}\}(i)]\ \ 
\hbox{\rm if}\ \ \gamma\!\in\! BC^{*}.
\end{array}\right.$$ 
Clearly, $\rho^{X}[BC_{\xi}]\! =\!\boraxi (X)$, by induction on $\xi$. The following is a consequence of 7B.1.(ii).(a) in \cite{Moschovakis}. It expresses the fact that the class of Borel sets is uniformly closed under complementation. 

\begin {lem} There is a recursive map $u_{\neg}\! :\!\omega^\omega\!\rightarrow\!\omega^\omega$ 
such that for each $1\!\leq\!\xi\! <\!\omega_{1}$ and for each $\gamma\!\in\! BC_{\xi}$, 
$u_{\neg}(\gamma )\!\in\! BC_{\xi +1}$, and 
$\rho^{X}[u_{\neg}(\gamma )]\! =\!\neg\rho^{X}(\gamma )$ for each 
recursively presented Polish space $X$.\end{lem}

\noindent\bf Proof.\rm\ Just copy the proof of 7B.1.(ii).(a) in \cite{Moschovakis}: it gives more than the statement in \cite{Moschovakis}.
\hfill{$\square$}

\bigskip

 In the sequel we will need a refinement of 7B.1.(iii) in \cite{Moschovakis}:
 
\begin {lem} Let $X$ be a recursively presented Polish space of type at most 1. 
Then there is a recursive map 
$u^X_{s}\! :\!\omega^\omega\!\times\! X\!\rightarrow\!\omega^\omega$ such that for 
each $1\!\leq\!\xi\! <\!\omega_{1}$, for each $\gamma\!\in\! BC_{\xi}$ and for each 
$x\!\in\! X$, $u^X_{s}(\gamma ,x)\!\in\! BC_{\xi}$, and 
$\rho^{Y}[u^X_{s}(\gamma ,x)]\! =\!\rho^{X\times Y}(\gamma )_{x}$ for each 
recursively presented Polish space $Y$.\end{lem}

 Some of the ideas of the proof are contained in 7A.3 in \cite{Moschovakis}.\bigskip
  
\noindent\bf Proof.\rm\ For $\xi\! =\! 1$, using the description of basic clopen sets in products (see 3B.1 in \cite{Moschovakis}), we define a subset of $\omega^\omega\!\times\! X\!\times\!\omega$ by 
$$(\gamma ,x,k)\!\in\! P\ \ \Leftrightarrow\ \ \exists i\!\in\!\omega\ \ \bigg(\ k = 
\Big< 0,\Big(\!\gamma^{*}(i)\!\Big)_{2}\Big>\ \ \hbox{\rm and}\ \ x\!\in\! 
N\Big[ X,\Big< 0,\Big(\!\gamma^{*}(i)\!\Big)_{1}\Big>\Big]\ \bigg).$$ 

 By 3D.5 in \cite{Moschovakis}, $\Boraone$ is closed under recursive substitutions, so that $P\!\in\!\Boraone$. By 3C.4 in \cite{Moschovakis}, there is $P^{*}\!\in\!\Borone (\omega^\omega\!\times\! X\!\times\!\omega^{2})$ with 
$$(\gamma ,x,k)\!\in\! P\ \Leftrightarrow\ \exists n\!\in\!\omega\ \ 
(\gamma ,x,k,n)\!\in\! P^*$$
(the idea is that in a space of type at most 1, an open set is a countable union of clopen sets). We define a map $g\! :\!\omega^\omega\!\times\! X\!\rightarrow\!\omega^\omega$ by  
$$g(\gamma ,x)(j)\! :=\!\left\{\!\!\!\!\!\!\!\!\begin{array}{ll} 
& (j\! -\! 1)_{0}\ \ \hbox{\rm if}\ \ j\! >\! 0\ \ \hbox{\rm and}\ \ 
[\gamma ,x,(j\! -\! 1)_{0},(j\! -\! 1)_{1}]\!\in\! P^{*}\hbox{\rm ,}\cr & \cr
& 0\ \ \hbox{\rm otherwise.}
\end{array}\right.$$
Clearly, $g$ is recursive and $g(\gamma ,x)\!\in\! BC_{1}$.\bigskip

\noindent $\bullet$ For the general case, we define a partial function $\psi\! :\! 
(\omega^\omega )^{2}\!\times\! X\!\times\!\omega\!\rightarrow\!
\omega^\omega$ by 
$$\psi (\varepsilon ,\gamma ,x,i)\! :=\!\left\{\!\!\!\!\!\!\!\!
\begin{array}{ll}
& g[\{\gamma^{*}\}(i),x]\ \ \hbox{\rm if}\ \ \{\gamma^{*}\}(i)(0)\! =\! 0
\hbox{\rm ,}\cr & \cr
& \{\varepsilon \}[\{\gamma^{*}\}(i),x]\ \ \hbox{\rm if}\ \ \{\gamma^{*}\}(i)(0)
\! =\! 1.
\end{array}\right.$$
The idea is that we want to build a recursive map $u^X_s$, that will have a recursive code $\varepsilon^{*}$. The function $\psi$ describes the properties that we want for $u^X_s$, and Kleene's Recursion Theorem will give the recursive code. By 3G.1 and 3G.2 in \cite{Moschovakis}, the collection of partial functions which are recursive on their domain is closed under composition, so that $\psi$ is recursive on its domain. Let $\nu\!\in\!\Boraone$ such that 
$$\psi (\varepsilon ,\gamma ,x,i)\! =\!\{\nu\}(\varepsilon ,\gamma ,x,i)$$ 
if $\psi (\varepsilon ,\gamma ,x,i)$ is defined. Note that $\{\nu\}(\varepsilon ,\gamma ,x,i)\! =\! 
\{ S^{(\omega^\omega )^{2}\times X,\omega ,\omega^\omega}_{\Boraone}
(\nu ,\varepsilon ,\gamma ,x)\}(i)$ when it is defined. We define a recursive map 
$\varphi\! :\! (\omega^\omega )^{2}\!\times\! X\!\rightarrow\!\omega^\omega$ by
$$\varphi (\varepsilon ,\gamma ,x)\! :=\!\left\{\!\!\!\!\!\!\!\!\begin{array}{ll}
& g(\gamma ,x)\ \ \hbox{\rm if}\ \ \gamma (0)\! =\! 0\hbox{\rm ,}\cr & \cr
& 1^\frown S^{(\omega^\omega )^{2}\times X,\omega ,\omega^\omega}_{\Boraone}
(\nu ,\varepsilon ,\gamma ,x)\ \ \hbox{\rm if}\ \ \gamma (0)\!\not=\! 0.
\end{array}\right.$$
By Kleene's Recursion Theorem, there is $\varepsilon^{*}\!\in\!\Boraone$ such that 
$\varphi (\varepsilon^{*},\gamma ,x)\! =\!\{\varepsilon^{*}\}(\gamma ,x)$ for 
each $(\gamma ,x)$ in $\omega^\omega\!\times\! X$. We put 
$u^X_{s}(\gamma ,x)\! :=\!\{\varepsilon^{*}\}(\gamma ,x)$. Note that the map $u^X_{s}$ is a total recursive map. We prove that $u^X_{s}(\gamma ,x)$ satisfies the required properties by induction on $\xi$.\bigskip

\noindent $\bullet$ Let $(\gamma ,x)\!\in\! BC_{1}\!\times\! X$. We have 
$u^X_{s}(\gamma ,x)\! =\!\{\varepsilon^{*}\}(\gamma ,x)\! =\!
\varphi (\varepsilon^{*} ,\gamma ,x)\! =\! g(\gamma ,x)$. So $u^X_{s}(\gamma ,x)$ 
is in $BC_{1}$, by the previous discussion. If moreover $Y$ is a recursively 
presented Polish space, then using the proof of 3B.1 in \cite{Moschovakis} we get
$$\begin{array}{ll}
y\!\in\!\rho^{X\times Y}(\gamma )_{x}\!\!\!
& \Leftrightarrow\ \exists i\!\in\!\omega\ \ (x,y)\!\in\! 
N[X\!\times\! Y,\gamma^{*}(i)]\cr 
& \Leftrightarrow\ \exists i\!\in\!\omega\ \ \bigg( 
y\!\in\! N\Big[Y,\Big< 0,\Big(\gamma^{*}(i)\Big)_{2}\Big>\Big]\ \ \hbox{\rm and}
~\ \ x\!\in\! N\Big[X,\Big< 0,\Big(\gamma^{*}(i)\Big)_{1}\Big>\Big]\bigg)\cr 
& \Leftrightarrow\ \exists k\!\in\!\omega\ \ [y\!\in\! N(Y,k)~\ \hbox{\rm and}~\ 
(\gamma ,x,k)\!\in\! P]\cr
& \Leftrightarrow\ \exists i\!\in\!\omega\ \ \Big( y\!\in\! N[Y,(i)_{0}]~\ \hbox{\rm and}
~\ [\gamma ,x,(i)_{0},(i)_{1}]\!\in\! P^*\Big)\cr
& \Leftrightarrow\ \exists i\!\in\!\omega\ \ y\!\in\! N\Big( Y,[g(\gamma ,x)]^{*}(i)\Big)\cr 
& \Leftrightarrow\  y\!\in\!\rho^{Y}[g(\gamma ,x)].
\end{array}$$
$\bullet$ Now let $(\gamma ,x)\!\in\! BC_{\xi}\!\times\! X$, with $\xi\! \geq\! 2$. We have 
$$u^X_{s}(\gamma ,x)\! =\!\{\varepsilon^{*}\}(\gamma ,x)\! =\!
\varphi (\varepsilon^{*} ,\gamma ,x)\! =\! 1^\frown 
S^{(\omega^\omega )^{2}\times X,\omega ,\omega^\omega}_{\Boraone}
(\nu ,\varepsilon^{*} ,\gamma ,x).$$
As $\gamma\!\in\! BC_{\xi}$, $\{\gamma^{*}\}(i)$ is defined for each 
integer $i$. In particular, $\psi (\varepsilon^{*} ,\gamma ,x,i)$ is defined 
for each $(\gamma ,x,i)$ in $\omega^\omega\!\times\! X\!\times\!\omega$ 
since $\{\gamma^{*}\}(i)(0)\!\in\! 2$, and equal to 
$$\{\nu\}(\varepsilon^{*} ,\gamma ,x,i)\! =\! 
\{ S^{(\omega^\omega )^{2}\times X,\omega ,\omega^\omega}_{\Boraone}
(\nu ,\varepsilon^{*} ,\gamma ,x)\}(i).$$ This shows that 
$\{ u^X_{s}(\gamma ,x)^{*}\}(i)$ is defined for each integer $i$. If $\{\gamma^{*}\}(i)(0)\! =\! 0$, then 
$$\{ u^X_{s}(\gamma ,x)^{*}\}(i)\! =\!g[\{\gamma^{*}\}(i),x]\! =\! u^X_{s}[\{\gamma^{*}\}(i),x].$$ 
As $\{\gamma^{*}\}(i)\!\in\! BC_1$, $u^X_{s}[\{\gamma^{*}\}(i),x]$ is in $BC_{1}$ too. Similarly, if $\{\gamma^{*}\}(i)(0)\! =\! 1$, then 
$$\{ u^X_{s}(\gamma ,x)^{*}\}(i)\! =\!
\{\varepsilon^{*}\}[\{\gamma^{*}\}(i),x]\! =\! u^X_{s}[\{\gamma^{*}\}(i),x].$$ 
Then $u^X_{s}[\{\gamma^{*}\}(i),x]\!\in\! BC_{\eta}$ for some $1\!\leq\!\eta\! <\!\xi$, by induction assumption. This shows that $u^X_{s}(\gamma ,x)$ is in $BC_{\xi}$. If $Y$ is a recursively 
presented Polish space, then $\rho^{Y}\Big( u^X_{s}[\{\gamma^{*}\}(i),x]\Big)\! =\!
\rho^{X\times Y}[\{\gamma^{*}\}(i)]_{x}$, by 
induction assumption. This shows that 
$\rho^{Y}[u^X_{s}(\gamma ,x)]\! =\!\rho^{X\times Y}(\gamma )_{x}$.\hfill{$\square$}\bigskip

 Lemma 3.2 expresses, among other things, the fact that the pointclasses $\boraxi$ are uniformly closed under taking sections at points in spaces of type at most $1$. Similarly, we now prove another lemma stating, among other things, that the pointclasses $\boraxi$ are uniformly closed under substitutions of partial recursive functions (when $\delta$ below is recursive).

\begin {lem} Let $X,Y$ be recursively presented Polish spaces. Then there is  
${u^{X,Y}_{r}\! :\! (\omega^\omega )^{2}\!\rightarrow\!\omega^\omega}$ recursive such 
that for each $1\!\leq\!\xi\! <\!\omega_{1}$, for each $\gamma\!\in\! BC_{\xi}$ and for 
each $\delta\!\in\!\omega^\omega$, $u^{X,Y}_{r}(\gamma ,\delta )\!\in\! BC_{\xi}$. 
Moreover, we have 
$x\!\in\!\rho^{X}[u^{X,Y}_{r}(\gamma ,\delta )]\Leftrightarrow\{\delta\}^{X,Y}(x)\!
\in\!\rho^{Y}(\gamma )$ if $\{\delta\}^{X,Y}(x)$ is defined.\end{lem}

\noindent\bf Proof.\rm\ The scheme of the proof is quite similar to that of Lemma 3.2. Indeed, this is again an application of Kleene's Recursion Theorem. For $\xi\! =\! 1$, we choose 
$P\!\in\!\Boraone (\omega^\omega\!\times\! X\!\times\!\omega )$ such that 
$$U(\delta ,x)\!\downarrow\ \ \ \Rightarrow\ \ \forall k\!\in\!\omega\ \ 
[\ U(\delta ,x)\!\in\! N(Y,k)\Leftrightarrow (\delta ,x,k)\!\in\! P\ ].$$
(this is possible since $U$ is recursive on its domain; see 7A.1 in \cite{Moschovakis}). By 3C.4 and 3C.5 in \cite{Moschovakis}, there is 
$P^{*}\!\in\!\Borone (\omega^\omega\!\times\!\omega^{3})$ with 
$$(\delta ,x,k)\!\in\! P\ \Leftrightarrow\ \exists i\!\in\!\omega\ \ \Big( x\!\in\! N[X,(i)_{0}]\ \ 
\hbox{\rm and}\ \ [\delta ,k,(i)_{0},(i)_{1}]\!\in\! P^*\Big).$$ 
We define a map $g\! :\! (\omega^\omega )^2\!\rightarrow\!\omega^\omega$ by  
$$g(\gamma ,\delta )(j)\! :=\!\left\{\!\!\!\!\!\!\!\!\begin{array}{ll} 
& \Big( (j\! -\! 1)_{0}\Big)_{0}\ \ \hbox{\rm if}\ \ j\! >\! 0\ \ 
\hbox{\rm and}\ \ \Big[\delta ,\gamma^{*}[(j\! -\! 1)_{1}],
\Big( (j\! -\! 1)_{0}\Big)_{0},\Big( (j\! -\! 1)_{0}\Big)_{1}\Big]\!\in\! P^*\hbox{\rm ,}\cr & \cr 
& 0\ \ \hbox{\rm otherwise.}
\end{array}\right.$$
Clearly, $g$ is recursive and $g(\gamma ,\delta )\!\in\! BC_{1}$.

\vfill\eject

\noindent $\bullet$ For the general case, we define a partial function 
$\psi\! :\! (\omega^\omega )^{3}\!\times\!\omega\!\rightarrow\!\omega^\omega$ by 
$$\psi (\varepsilon ,\gamma ,\delta ,i)\! :=\!\left\{\!\!\!\!\!\!\!\!
\begin{array}{ll}
& g[\{\gamma^{*}\}(i),\delta ]\ \ \hbox{\rm if}\ \ \{\gamma^{*}\}(i)(0)\! =\! 0
\hbox{\rm ,}\cr & \cr
& \{\varepsilon \}[\{\gamma^{*}\}(i),\delta ]\ \ \hbox{\rm if}\ \ \{\gamma^{*}\}
(i)(0)\! =\! 1.
\end{array}\right.$$
We argue as in the proof of Lemma 3.2 to define 
$\varphi\! :\! (\omega^\omega )^{3}\!\rightarrow\!\omega^\omega$, and we put 
$u^{X,Y}_{r}(\gamma ,\delta )\! :=\!\{\varepsilon^{*}\}(\gamma ,\delta )$. 
The map $u^{X,Y}_{r}$ is a total recursive map. We show that 
$u^{X,Y}_{r}(\gamma ,\delta )$ satisfies the required properties by induction on $\xi$.\bigskip 

\noindent $\bullet$ If $(\gamma ,\delta )\!\in\! BC_{1}\!\times\!\omega^\omega$ and 
$\{\delta\}^{X,Y}(x)$ is defined, then
$$\begin{array}{ll}
x\!\in\!\rho^{X}[u^{X,Y}_{r}(\gamma ,\delta )]\!\!\!
& \Leftrightarrow\ \exists k\!\in\!\omega\ \ x\!\in\! 
N[X,g(\gamma ,\delta )^{*}(k)]\cr 
& \Leftrightarrow\ \exists k\!\in\!\omega\ \ x\!\in\! 
N\Big[ X,\Big( (k)_{0}\Big)_{0}\Big]\ \hbox{\rm and}\ 
\Big[\delta ,\gamma^{*}[(k)_{1}],\Big( (k)_{0}\Big)_{0},
\Big( (k)_{0}\Big)_{1}\Big]\!\in\! P^*\cr 
& \Leftrightarrow\ \exists j\!\in\!\omega\ \ \exists i\!\in\!\omega\ \ 
\Big( x\!\in\! N[X,(i)_{0}]~\ \hbox{\rm and}~\ 
[\delta ,\gamma^{*}(j),(i)_{0},(i)_{1}]\!\in\! P^*\Big)\cr 
& \Leftrightarrow\ \exists j\!\in\!\omega\ \ [\delta ,x,\gamma^{*}(j)]\!\in\! P\cr 
& \Leftrightarrow\ \exists j\!\in\!\omega\ \ \{\delta\}^{X,Y}(x)\!\in\! 
N[Y,\gamma^{*}(j)]\cr 
& \Leftrightarrow\ \{\delta\}^{X,Y}(x)\!\in\!\rho^{Y}(\gamma ).
\end{array}$$
$\bullet$ Now let $\gamma\!\in\! BC_{\xi}$ with $\xi\! \geq\! 2$, and $\delta\!\in\!\omega^\omega$. As in the proof of Lemma 3.2, $u^{X,Y}_{r}(\gamma ,\delta )\!\in\! BC_{\xi}$. If $\{\delta\}^{X,Y}(x)$ is defined, then 
$$x\!\in\!\rho^{X}\Big( u^{X,Y}_{r}[\{\gamma^{*}\}(i),\delta ]\Big)\ 
\Leftrightarrow\ \{\delta\}^{X,Y}(x)\!\in\!\rho^{Y}[\{\gamma^{*}\}(i)]\hbox{\rm ,}$$
by induction assumption. This shows that $x\!\in\!\rho^{X}[u^{X,Y}_{r}(\gamma ,\delta )]
\Leftrightarrow\{\delta\}^{X,Y}(x)\!\in\!\rho^{Y}(\gamma)$.\hfill{$\square$}\bigskip

 As a corollary, one can prove  the uniform closure of the pointclasses $\boraxi$ under fixations of recursive arguments. It is sometimes convenient to ``view a code in $BC_1$ as an element of 
 $BC_2$", even if it is not formally correct. The next lemma expresses this:
 
\begin {lem} Let $X$ be a recursively presented Polish space. Then there is 
$u^X_{*}\! :\!\omega^\omega\!\rightarrow\!\omega^\omega$ recursive such that for each 
$\gamma\!\in\! BC_1$ (resp., $BC^*$), $u^X_{*}(\gamma )\!\in\! BC_{2}$ (resp.,  
$u^X_{*}(\gamma )\! =\!\gamma$),  and $\rho^{X}[u^X_{*}(\gamma )]\! =\!\rho^{X}(\gamma )$.\end{lem}
 
\noindent\bf Proof.\rm\ We define $R\!\in\!\Boraone (\omega^\omega\!\times\! X)$ by 
$(\gamma ,x)\!\in\! R\ \Leftrightarrow\ \exists i\!\in\!\omega\ \ x\!\in\! N[X,\gamma^*(i)]$. As 
$R\!\in\!\Boratwo$ there is $C\!\in\!\Bormone (\omega\!\times\!\omega^\omega\!\times\! X)$ 
such that $R\! =\!\exists^\omega C$. Let $\varepsilon_0\!\in\!\Boraone$ such that 
$\neg C\! =\!\bigcup_{i\in\omega}\ N[\omega\!\times\!\omega^\omega\!\times\! X,\varepsilon_0(i)]$. Note that $0^\frown\varepsilon_0\!\in\!\Boraone\cap BC_1$ and 
$\neg C\! =\!\rho^{\omega\times\omega^\omega\times X}(0^\frown\varepsilon_0)$. Using Lemma 3.2, 
we see the existence of $\gamma_0\!\in\!\Boraone$ such that 
$\{\gamma_0\}(\gamma ,i)\! =\! u_s^{\omega\times\omega^\omega}(0^\frown\varepsilon_0,i,\gamma )$ for each $(\gamma ,i)\!\in\!\omega^\omega\!\times\!\omega$. Then we define $u^X_{*}(\gamma )\! :=\! 1^\frown 
S^{\omega^\omega,\omega ,\omega^\omega}_{\Boraone}(\gamma_0,\gamma )$ if $\gamma\!\in\! BC_1$, $\gamma$ otherwise.\hfill{$\square$}\bigskip 
 
 We now prove another lemma stating, among other things, that the pointclasses 
$\boraxi$ are uniformly closed under finite intersections and unions:

\begin {lem} Let $X$ be a recursively presented Polish space. There is 
$u^X_{f}\! :\! 2\!\times\!\omega\!\times\! \omega^\omega\!\rightarrow\!\omega^\omega$ recursive such that for each $(\xi ,a,n,\gamma )\!\in\! (\omega_{1}\!\setminus\!\{ 0\})\!\times\! 2\!\times\!\omega\!\times\! \omega^\omega$,\smallskip 

\noindent (a) If $(\gamma )_i\!\in\! BC_1\cup BC_{\xi}$ for each $i\!\leq\! n$, then 
$u^X_{f}(a,n,\gamma )$ is in $BC_1\cup BC_{\xi}$. Moreover, the equalities 
$\rho^{X}[u^X_{f}(0,n,\gamma )]\! =\!\bigcap_{i\leq n}\rho^{X}[(\gamma )_{i}]$ and 
$\rho^{X}[u^X_{f}(1,n,\gamma )]\! =\!\bigcup_{i\leq n}\rho^{X}[(\gamma )_{i}]$ hold.\smallskip

\noindent (b) If moreover $\xi\! \geq\! 2$ and $(\gamma )_i\!\in\! BC_{\xi}$ for some $i\!\leq\! n$, then 
$u^X_{f}(a,n,\gamma )$ is in $BC_{\xi}$.\end{lem}
 
\noindent\bf Proof.\rm\ Once again, this is an application of Kleene's Recursion Theorem. For $\xi\! =\! 1$, by 3B.2 in \cite{Moschovakis} there is $f \! :\!\omega^3\!\rightarrow\!\omega$ recursive such that,  for $(u,n)\!\in\!\omega^2$, 
$$\bigcap_{i\leq n}\ N[X,(u)_i]\! =\!\bigcup_{m\in\omega}\ N[X,f(u,n,m)].$$
We set $g(a,n,\gamma )(0)\! :=\! 0$ and 
$$g(0,n,\gamma )(i\! +\! 1)\! :=\! f\bigg(\!\! <[(\gamma )_0]^*\Big[\Big( (i)_0\Big)_0\Big],... ,
[(\gamma )_n]^*\Big[\Big( (i)_0\Big)_n\Big]>,n,(i)_1\bigg)\hbox{\rm ,}$$
$$g(1,n,\gamma )(i\! +\! 1)\! :=\!\left\{\!\!\!\!\!\!\!\!
\begin{array}{ll}
& [(\gamma )_{(i)_0}]^*[(i)_1]~\ \hbox{\rm if}~\ (i)_0\!\leq\! n\hbox{\rm ,}\cr & \cr 
& 0~\ \hbox{\rm otherwise.}
\end{array}
\right.$$
Note that $g(a,n,\gamma )\!\in\! BC_1$. If $(\gamma )_i\!\in\! BC_1$ for each $i\!\leq\! n$, then we get 
$$\begin{array}{ll}\rho^{X}[g(0,n,\gamma )]\!\!\!\!
& =\bigcup_{i\in\omega}\ N[X,g(0,n,\gamma )(i\! +\! 1)]\cr
& =\bigcup_{(j,m)\in\omega^2}\ N\Big[X,f\Big(\!\! <[(\gamma )_0]^*[( j)_0],... ,
[(\gamma )_n]^*[( j)_n]>,n,m\Big)\Big]\cr
& =\!\bigcup_{j\in\omega}\ \bigcap_{i\leq n}\ N\Big(X,[(\gamma )_i]^*[( j)_i]\Big)\cr
& =\bigcap_{i\leq n}\ \bigcup_{j\in\omega}\ N\Big(X,[(\gamma )_{i}]^*(j)\Big)\cr
& =\bigcap_{i\leq n}\ \rho^{X}[(\gamma )_i].
\end{array}$$
Moreover, 
$$\rho^{X}[g(1,n,\gamma )]
\! =\!\bigcup_{i\in\omega}\ N[X,g(1,n,\gamma )(i\! +\! 1)]
\! =\!\bigcup_{i\leq n}\ \bigcup_{j\in\omega}\ N\Big(X,[(\gamma )_{i}]^*(j)\Big)
\! =\!\bigcup_{i\leq n}\ \rho^{X}[(\gamma )_i].$$
$\bullet$ For the general case, using Lemma 3.4 we define a partial function 
$h\! :\! \omega^\omega\!\times\!\omega^2\!\rightarrow\!\omega^\omega$ by 
$$\Big( h(\gamma ,n,j)\Big)_i\! :=\!\left\{\!\!\!\!\!\!\!\!
\begin{array}{ll}
& \Big\{\Big( u^X_{*}[(\gamma )_i]\Big)^*\Big\}[(j)_i]~\ \hbox{\rm if}~\ i\!\leq\! n\hbox{\rm ,}\cr & \cr
& 0^\infty\ \ \hbox{\rm otherwise.}
\end{array}\right.$$ 
It allows us to define another partial function $\psi\! :\! \omega^\omega\!\times\! 2\!\times\! \omega\!\times\!\omega^\omega\!\times\!\omega\!\rightarrow\!\omega^\omega$ by 
$$\psi (\varepsilon ,a,n,\gamma ,j)\! :=\!\{\varepsilon \}[1\! -\! a,n,h(\gamma ,n,j)].$$
We argue as in the proof of Lemma 3.2 to define $\nu$ and a recursive map 
$\varphi\! :\!\omega^\omega\!\times\! 2\!\times\! \omega\!\times\!\omega^\omega\!\rightarrow\!\omega^\omega$ by
$$\varphi (\varepsilon ,a,n,\gamma )\! :=\!\left\{\!\!\!\!\!\!
\begin{array}{ll}
& g(a,n,\gamma )\ \ \hbox{\rm if}\ \ (\gamma )_i(0)\! =\! 0\ \hbox{\rm for\ each}\ i\!\leq\! n\hbox{\rm ,}\cr & \cr
& 1^\frown S^{\omega^\omega\times 2\times\omega\times\omega^\omega ,\omega ,\omega^\omega}_{\Boraone}
(\nu ,\varepsilon ,a,n,\gamma )\ \ \hbox{\rm if}\ \ (\gamma )_i(0)\!\not=\! 0\ \hbox{\rm for\ some}\ i\!\leq\! n.
\end{array}\right.$$
By Kleene's Recursion Theorem, there is $\varepsilon^{*}\!\in\!\Boraone$ such that 
$\varphi (\varepsilon^{*},a,n,\gamma )\! =\!\{\varepsilon^{*}\}
(a,n,\gamma )$ for each $(a,n,\gamma )$ in $2\!\times\!\omega\!\times\!\omega^\omega$. We put 
$u^{X}_{f}(a,n,\gamma )\! :=\!\{\varepsilon^{*}\}(a,n,\gamma )$. The map $u^{X}_{f}$ is a total recursive map. We show that $u^{X}_{f}(a,n,\gamma )$ satisfies the required properties by induction on $\xi$.\bigskip 

\noindent $\bullet$ Assume that $(\gamma )_i\!\in\! BC_{1}$ holds for each $i\!\leq\! n$. We have 
$u^{X}_{f}(a,n,\gamma )\! =\! g(a,n,\gamma )$, so we are done, by the previous discussion. Assume now that $\xi\! \geq\! 2$, and that $(\gamma )_i\!\in\! BC_1\cup BC_{\xi}$ for each $i\!\leq\! n$. We may assume that $(\gamma )_i\!\in\! BC_{\xi}$ holds for some $i\!\leq\! n$. Then 
$\Big\{\Big( u^X_{*}[(\gamma )_i]\Big)^*\Big\}(k)$ is defined for each integer $k$. In particular, 
$h(\gamma ,n,j)$ and $\psi (\varepsilon^{*} ,a,n,\gamma ,j)$ are defined for each 
$(a,j)$ in $2\!\times\!\omega$. Thus  
$$\{ u^{X}_{f}(a,n,\gamma )^{*}\}(j)\! =\! u^X_f[1\! -\! a,n,h(\gamma ,n,j)]$$ 
is defined for each integer $j$. As 
$\Big\{\Big( u^X_{*}[(\gamma )_i]\Big)^*\Big\}(k)$ is in some $BC_{\eta_k}$ with 
$1\!\leq\!\eta_k\! <\!\xi$ for each integer $k$, there is $1\!\leq\!\eta\! <\!\xi$ such that $\Big( h(\gamma ,n,j)\Big)_i$ is in $BC_1\cup BC_{\eta}$ for each $i\!\leq\! n$. By induction assumption, we get 
$u^X_f[1\! -\! a,n,h(\gamma ,n,j)]\!\in\! BC_1\cup BC_{\eta}$. This shows that 
$u^{X}_{f}(a,n,\gamma )\!\in\! BC_{\xi}$. Moreover, by induction assumption we get 
$$\begin{array}{ll}\rho^X[u^X_f(0,n,\gamma )]\!\!
& \! =\ \bigcup_{j\in\omega}\ \neg\bigcup_{i\leq n}\ \rho^X\Big[\Big(h(\gamma ,n,j)\Big)_i\Big]\cr & \cr
& \! =\ \bigcup_{j\in\omega}\ \bigcap_{i\leq n}\ \neg\rho^X
\Big[\Big\{\Big( u^X_{*}[(\gamma )_i]\Big)^*\Big\}[(j)_i]\Big]\cr & \cr
& \! =\ \bigcap_{i\leq n}\ \bigcup_{j\in\omega}\ \neg\rho^X
\Big[\Big\{\Big( u^X_{*}[(\gamma )_i]\Big)^*\Big\}(j)\Big]\cr & \cr
& \! =\ \bigcap_{i\leq n}\ \rho^X\Big( u^X_{*}[(\gamma )_i]\Big)\cr & \cr
& \! =\ \bigcap_{i\leq n}\ \rho^X[(\gamma )_i].
\end{array}$$
Similarly, we get $\rho^X[u^X_f(1,n,\gamma )]\! =\!\bigcup_{i\leq n}\ 
\rho^X[(\gamma )_i]$.\hfill{$\square$}\bigskip

 In the sequel we will need a last closure property, asserting, among other things, that the pointclasses $\boraxi$ are uniformly closed under $\exists^\omega$:
 
\begin {lem} (a) There is a recursive map 
$u_{\exists}\! :\!\omega^\omega\!\rightarrow\!\omega^\omega$ such that for each 
$1\!\leq\!\xi\! <\!\omega_{1}$ and for each $\gamma\!\in\! BC_{\xi}$, 
$u_{\exists}(\gamma )\!\in\! BC_{\xi}$, and  
$x\!\in\!\rho^{X}[u_{\exists}(\gamma )]\ \Leftrightarrow\ \exists n\!\in\!\omega
\ \ (n,x)\!\in\!\rho^{\omega\times X}(\gamma )$, for each recursively 
presented Polish space $X$ and for each $x\!\in\! X$.

\noindent (b) There is a recursive map $u_{(.)}
\! :\!\omega^\omega\!\rightarrow\!\omega^\omega$ such that for each $1\!\leq\!\xi\! <\!\omega_{1}$, $(\gamma )_n\!\in\! BC_{\xi}$ for each $n\!\in\!\omega$ implies that 
$u_{(.)}(\gamma )\!\in\! BC_{\xi}$, and 
$x\!\in\!\rho^{X}[u_{(.)}(\gamma )]\Leftrightarrow\exists n\!\in\!\omega ~\ x\!\in\!
\rho^{X}[(\gamma )_n]$ for each recursively presented Polish space $X$ and for 
each $x\!\in\! X$.\end{lem}

\noindent\bf Proof.\rm\ Once again we code the properties that we want. So a look at the end of the proofs of (a) and (b) can give an idea of the intuition behind them.\bigskip

\noindent (a) By 3B.1 in \cite{Moschovakis}, there are $g$ and $h$ recursive such that 
$N(\omega\!\times\! X,k)\! =\! N[\omega ,g(k)]\!\times\! N[X,h(k)]$ for each integer $k$. If $\gamma (0)\! =\! 0$, then we put
$$u_{\exists}(\gamma )(j)\! :=\!\left\{\!\!\!\!\!\!\!\!
\begin{array}{ll}
& h\Big(\gamma^*[(j\! -\! 1)_0]\Big)\ \ \hbox{\rm if}\ \ j\! >\! 0\ \ \hbox{\rm and}\ \ (j\! -\! 1)_1\!\in\! N
\Big[\omega ,g\Big(\gamma^*[(j\! -\! 1)_0]\Big)\Big]\hbox{\rm ,}\cr & \cr 
& 0\ \hbox{\rm otherwise.}
\end{array}\right.$$

 Using Lemma 3.2, we define  a partial function 
$f\! :\!\omega^\omega\!\times\!\omega\!\rightarrow\!\omega^\omega$ by 
$$f(\gamma ,i)\! :=\! u^\omega_s\Big(\{\gamma^*\}[(i)_1],(i)_0\Big).$$
As $f$ is recursive on its domain, there is $\varepsilon_0\!\in\!\Boraone$ such 
that $f(\gamma ,i)\! =\!\{\varepsilon_0\}(\gamma ,i)$ if  $f(\gamma ,i)$ is 
defined. If $\gamma (0)\!\not=\! 0$, then we put 
$u_{\exists}(\gamma )\! :=\! 1^\frown S_{\Boraone}^{\omega^\omega ,\omega ,
\omega^\omega}(\varepsilon_0,\gamma )$. This defines a recursive map $u_{\exists}$. 
If $\gamma (0)\!\not=\! 0$ and $\gamma\!\in\! BC_{\xi}$, then 
$\{[u_{\exists}(\gamma )]^*\}(i)\! =\!\{S_{\Boraone}^{\omega^\omega ,\omega ,
\omega^\omega}(\varepsilon_0,\gamma )\}(i)\! =\!\{\varepsilon_0\}(\gamma ,i)\! =
\! f(\gamma ,i)$. Thus 
$u_{\exists}(\gamma )\!\in\! BC_{\xi}$, even if $\gamma (0)\! =\! 0$. Let $x\!\in\! X$. If 
$\gamma (0)\! =\! 0$, then 
$$\begin{array}{ll}\exists n\!\in\!\omega\ \ (n,x)\!\in\!\rho^{\omega\times X}(\gamma )\!\!\!\!
& \Leftrightarrow\exists n\!\in\!\omega\ \exists p\!\in\!\omega\ \ (n,x)\!\in\! N[\omega\!\times\! X,\gamma^*(p)]\cr 
& \Leftrightarrow\exists n\!\in\!\omega\ \exists p\!\in\!\omega\ \ n\!\in\! N\Big(\omega ,g[\gamma^*(p)]\Big)~\ 
\hbox{\rm and}~\ x\!\in\! N\Big(X,h[\gamma^*(p)]\Big)\cr
& \Leftrightarrow\exists i\!\in\!\omega\ \ (i)_1\!\in\! N\Big[\omega ,g\Big(\gamma^*[(i)_0]\Big)\Big]~\ 
\hbox{\rm and}~\ x\!\in\! N\Big[X,h\Big(\gamma^*[(i)_0]\Big)\Big]\cr 
& \Leftrightarrow x\!\in\!\rho^{X}[u_{\exists}(\gamma )].\end{array}$$
If $\gamma (0)\!\not=\! 0$, then 
$$\begin{array}{ll}\exists n\!\in\!\omega\ \ (n,x)\!\in\!\rho^{\omega\times X}(\gamma )\!\!\!\!
& \Leftrightarrow\exists n\!\in\!\omega\ \exists p\!\in\!\omega\ \ (n,x)\!\notin\!\rho^{\omega\times X}
[\{\gamma^*\}(p)]\cr 
& \Leftrightarrow\exists i\!\in\!\omega\ \ \neg\ [(i)_0,x]\!\in\!\rho^{\omega\times X}
\Big(\{\gamma^*\}[(i)_1]\Big)\cr 
& \Leftrightarrow\exists i\!\in\!\omega\ \ \neg\ x\!\in\!\rho^{X}[f(\gamma ,i)]\cr
& \Leftrightarrow x\!\in\!\rho^{X}[u_{\exists}(\gamma )].\end{array}$$
(b) If $(\gamma )_0(0)\! =\! 0$, then we put
$$u_{(.)}(\gamma )(j)\! :=\!\left\{\!\!\!\!\!\!\!\!
\begin{array}{ll}
& 0\ \ \hbox{\rm if}\ \ j\! =\! 0\hbox{\rm ,}\cr & \cr
& [(\gamma )_{(j-1)_0}]^*[(j\! -\! 1)_1]\ \hbox{\rm otherwise.}
\end{array}\right.$$
We define  a partial function $f'\! :\!\omega^\omega\!\times\!\omega\!\rightarrow\!\omega^\omega$ by 
$f'(\gamma ,i)\! :=\!\{ [(\gamma )_{(i)_0}]^*\}[(i)_1]$. As $f'$ is recursive on 
its domain, there is $\gamma_0\!\in\!\Boraone$ such that 
$f'(\gamma ,i)\! =\!\{\gamma_0\}(\gamma ,i)$ if $f'(\gamma ,i)$ is defined. If 
$(\gamma )_0(0)\!\not=\! 0$,  then we put $u_{(.)}(\gamma )\! :=\! 1^\frown 
S_{\Boraone}^{\omega^\omega ,\omega ,\omega^\omega}(\gamma_0,\gamma )$. This 
defines a recursive map $u_{(.)}$. If $\xi\! \geq\! 2$ and 
$(\gamma )_n\!\in\! BC_{\xi}$ for each integer $n$, then 
$\{[u_{(.)}(\gamma )]^*\}(i)\! =\!\{S_{\Boraone}^{\omega^\omega ,\omega ,
\omega^\omega}(\gamma_0,\gamma )\}(i)\! =\!\{\gamma_0\}(\gamma ,i)\! =\! 
f'(\gamma ,i)$. Thus $u_{(.)}(\gamma )\!\in\! BC_{\xi}$, even if $\xi\! =\! 1$. Let $x\!\in\! X$. If $(\gamma )_0(0)\! =\! 0$, then 
$$\begin{array}{ll}\exists n\!\in\!\omega ~\ x\!\in\!\rho^{X}[(\gamma )_n]\!\!\!\!
& \Leftrightarrow\exists n\!\in\!\omega\ \exists p\!\in\!\omega\ \ x\!\in\! N\Big( X,[(\gamma )_n]^*(p)\Big)\cr 
& \Leftrightarrow\exists i\!\in\!\omega\ \ x\!\in\! N\Big( X,[(\gamma )_{(i)_0}]^*[(i)_1]\Big)\cr 
& \Leftrightarrow x\!\in\!\rho^{X}[u_{(.)}(\gamma )].\end{array}$$
If $(\gamma )_0(0)\!\not=\! 0$, then 
$$\begin{array}{ll}\exists n\!\in\!\omega ~\ x\!\in\!\rho^{X}[(\gamma )_n]\!\!\!\!
& \Leftrightarrow\exists n\!\in\!\omega\ \exists p\!\in\!\omega\ \ x\!\notin\!\rho^{X}
\Big(\{ [(\gamma )_n]^*\}(p)\Big)\cr 
& \Leftrightarrow\exists i\!\in\!\omega\ \ \neg\ x\!\in\!\rho^{X}
\Big(\{ [(\gamma )_{(i)_0}]^*\}[(i)_1]\Big)\cr 
& \Leftrightarrow\exists i\!\in\!\omega\ \ \neg\ x\!\in\!\rho^{X}[f'(\gamma ,i)]\cr
& \Leftrightarrow x\!\in\!\rho^{X}[u_{(.)}(\gamma )].\end{array}$$
This finishes the proof.\hfill{$\square$}

\vfill\eject

\noindent $\underline{\bf{The\ hyperarithmetical\ hierarchy.}}$\bigskip

 The notion of a hyperarithmetical set is defined in \cite{Moschovakis} (see 7B): a subset of a recursively 
presented Polish space is $hyperarithmetical$ if it is Borel and has 
a recursive Borel code. We can define a hyperarithmetical hierarchy, extending 
the arithmetical hierarchy. The following characterization of the arithmetical 
pointclasses $\Boran$ can be found in  Louveau's notes \cite{Louveau}:
 
\begin {thm} Let $X$ be a recursively presented Polish space, and $n\!\geq\! 1$ an integer. Then 
$$\Boran (X)\! =\!\{\rho^{X}(\gamma )\mid\gamma\!\in\!\Boraone\cap BC_{n}\}.$$
\end{thm}

 Actually, we will use only a small part of it. More specifically, we will only use the fact that if 
$P$ is $\Boraone (X)$, then there is $\gamma\!\in\!\Boraone\cap BC_{1}$ with 
$P\! =\!\rho^{X}(\gamma )$. It is very simple: there is $\varepsilon\!\in\!\Boraone$ 
such that $P\! =\!\bigcup_{i\in\omega}\ N[X,\varepsilon (i)]\! =\!
\rho^{X}(0^\frown\varepsilon )$. Thus $\gamma\! :=\! 0^\frown\varepsilon$ is 
suitable. The following definition comes naturally after Theorem 3.7, and can be found in  \cite{Louveau}:
 
\begin{defi} Let $X$ be a recursively presented Polish space, and 
$1\!\leq\!\xi\! <\!\omega_{1}$. Then we set 
$$\begin{array}{ll}
\Boraxi (X) & \!\!\!\! =\{\rho^{X}(\gamma )\mid\gamma\!\in\!\Boraone\cap BC_{\xi}\}
\hbox{\rm ,}\cr &\cr 
\Bormxi (X) & \!\!\!\! =\check\Boraxi (X)\hbox{\rm ,}\cr & \cr
\Borxi (X) & \!\!\!\! =\Boraxi (X)\cap\Bormxi (X).
\end{array}$$
This defines the $hyperarithmetical\ hierarchy$.
\end{defi}

 Note that Lemma 3.3 (resp., 3.5, 3.6) implies that the hyperarithmetical pointclasses are closed under recursive substitutions (resp., finite intersections and unions, $\exists^\omega$). Now we construct recursive maps giving codes for the basic neighborhoods and their complements in spaces of type at most 1.
 
\begin {lem} Let $X$ be a recursively presented Polish space of type at most 1.

\noindent (a) There is a recursive map $u_{N}\! :\!\omega\!\rightarrow\!\omega^\omega$ such that 
$u_{N}(k)\!\in\! BC_{1}$, and $\rho^{X}[u_{N}(k)]\! =\!\ N(X,k)$. 

\noindent (b)  There is a recursive map $u^X_{\neg N}\! :\!\omega\!\rightarrow\!\omega^\omega$ such that $u^X_{\neg N}(k)\!\in\! BC_{1}$, and $\rho^{X}[u^X_{\neg N}(k)]\! =\!
\neg N(X,k)$.\end{lem}

\noindent\bf Proof.\rm\ (a) Put $u_{N}(k)\! :=\! 0k0^\infty$.\bigskip

\noindent (b)  By 3C.3 in \cite{Moschovakis}, the equivalence $(x,k)\!\in\! R\Leftrightarrow x\!\notin\! N(X,k)$ defines 
$R\!\in\!\Boraone (X\!\times\!\omega )$. By Theorem 3.7 there is $\gamma_0\!\in\!\Boraone\cap BC_1$ with $R\! =\!\rho^{X\times\omega}(\gamma_0)$. Using Lemma 3.2 we set 
$u^X_{\neg N}(k)\! :=\! u_s^{\omega}(\gamma_0,k)$.\hfill{$\square$}\bigskip

 Now we use this to prove that, uniformly in $\xi\!\geq\! 2$, a set in the pointclass $\boraxi (X)$ (resp., 
$\Boraxi (X)$) is the disjoint union of sets in $\bormlxi$ (resp., $\Bormlxi$), if $X$ is a space of 
type at most 1. We will use the notation 
$$E\! =\!\displaystyle{\bigcup^{\bullet}_{i\in\omega}}\ E_i$$ 
to express the fact that $E$ is the disjoint union of the $E_i$'s.
  
\begin {lem} Let $X$ be a recursively presented Polish space of type at most 1. 
Then there is a recursive map $u^X_{d}\! :\!\omega^\omega\!\rightarrow\!\omega^\omega$ 
such that $u^X_{d}(\gamma )\!\in\! BC_{\xi}$ if $\gamma\!\in\! BC_{\xi}$, for each 
$1\!\leq\!\xi\! <\!\omega_1$. Moreover,\bigskip 

\noindent (a) There is a recursive map 
$u^X_{c}\! :\!\omega^\omega\!\times\!\omega\!\rightarrow\!\omega^\omega$ such that\bigskip 

\noindent (1) $u^X_c(\gamma ,i)\!\in\! BC_1$ for each $(\gamma ,i)\!\in\!\omega^\omega\!\times\!\omega$.\smallskip

\noindent (2) $\{[u^X_d(\gamma )]^*\}(i)$ is defined, in $BC_1$, 
$\rho^{X}\Big(\{[u^X_d(\gamma )]^*\}(i)\Big)\!\in\!\borone$ and 
$\neg\rho^{X}\Big(\{[u^X_d(\gamma )]^*\}(i)\Big)\! =\!\rho^{X}[u^X_c(\gamma ,i)]$ for each 
$(\gamma ,i)\!\in\! BC_{1}\!\times\!\omega$.

\noindent (b) If $1\!\leq\!\xi\! <\!\omega_{1}$ and $\gamma\!\in\! BC_{\xi}$, 
then $\rho^{X}(\gamma )\! =\!\displaystyle{\bigcup^{\bullet}_{i\in\omega}}\ 
\neg\rho^{X}\Big(\{[u^X_d(\gamma )]^*\}(i)\Big)$.\end{lem}

\noindent\bf Proof.\rm\ For $\xi\! =\! 1$, a look at the computation of  $\rho^X(\gamma )$   at the end of this point can help to understand what is going on. We first define a map 
$\tilde f\! :\!\omega^\omega\!\times\!\omega\!\rightarrow\!\omega^\omega$, using Lemma 3.9, as follows:
$$\Big(\tilde f(\gamma ,i)\Big)_j\! :=\!\left\{\!\!\!\!\!\!\!\!\begin{array}{ll}
& u^X_{\neg N}[\gamma^*(j)]~\ \hbox{\rm if}~\ j\! <\! i\hbox{\rm ,}\cr & \cr 
& u_{N}[\gamma^*(i)]~\ \hbox{\rm if}~\ j\!\geq\! i.
\end{array}\right.$$
As $\tilde f$ is recursive, the formula $u^X_{c}(\gamma ,i)\! :=\! u^X_f[0,i,\tilde f(\gamma ,i)]$ defines 
$u^X_{c}$ recursive such that $u^X_c(\gamma ,i)\!\in\! BC_1$ for each $(\gamma ,i)\!\in\!\omega^\omega\!\times\!\omega$ (see Lemma 3.5). Then, using Lemma 3.9, we define a map $f\! :\!\omega^\omega\!\times\!\omega\!\rightarrow\!\omega^\omega$:
$$\Big( f(\gamma ,i)\Big)_j\! :=\!\left\{\!\!\!\!\!\!\!\!\begin{array}{ll}
& u_{N}[\gamma^*(j)]~\ \hbox{\rm if}~\ j\! <\! i\hbox{\rm ,}\cr & \cr 
& u^X_{\neg N}[\gamma^*(i)]~\ \hbox{\rm if}~\ j\!\geq\! i.
\end{array}\right.$$
As $f$ is recursive, and using Lemma 3.5, there is $\varepsilon_0\!\in\!\Boraone$ such that 
$\{\varepsilon_0\}(\gamma ,i)\! =\! u^X_f[1,i,f(\gamma ,i)]\!\in\! BC_1$ for each 
${(\gamma ,i)\!\in\!\omega^\omega\times\omega}$. We define a recursive map 
$g\! :\!\omega^\omega\!\rightarrow\!\omega^\omega$ by $g(\gamma )\! :=\! 0^\frown 
S^{\omega^\omega ,\omega ,\omega^\omega}_{\Boraone}(\varepsilon_0,\gamma )$. If 
$(\gamma ,i)\!\in\! BC_1\!\times\!\omega$, then $\{[g(\gamma )]^*\}(i)\! =\! 
\{S^{\omega^\omega ,\omega ,\omega^\omega}_{\Boraone}(\varepsilon_0,\gamma )\}(i)\! =\! 
\{\varepsilon_0\}(\gamma ,i)\! =\! u^X_f[1,i,f(\gamma ,i)]$ is defined,  
$\rho^{X}\Big(\{[g(\gamma )]^*\}(i)\Big)\!\in\!\borone$ since it is a finite union of clopen sets, and 
$\neg\rho^{X}\Big(\{[g(\gamma )]^*\}(i)\Big)\! =\!\rho^{X}[u^X_c(\gamma ,i)]$. Moreover,
$$\begin{array}{ll}\rho^X(\gamma )\!\!\! 
& \! =\bigcup_{i\in\omega}\ N[X,\gamma^*(i)]
\! =\!\displaystyle{\bigcup^{\bullet}_{i\in\omega}}\ N[X,\gamma^*(i)]\!
\setminus\!\Big(\bigcup_{j<i}\ N[X,\gamma^*(j)]\Big)\cr 
& \! =\!\displaystyle{\bigcup^{\bullet}_{i\in\omega}}\ \neg
\Big[\rho^X\Big( u^X_{\neg N}[\gamma^*(i)]\Big)\cup
\bigcup_{j<i}\ \rho^X\Big( u_{N}[\gamma^*(j)]\Big)\Big]\cr 
& \! =\displaystyle{\bigcup^{\bullet}_{i\in\omega}}\ \neg\rho^{X}\Big( u^X_f[1,i,f(\gamma ,i)]\Big)
\! =\!\displaystyle{\bigcup^{\bullet}_{i\in\omega}}\ \neg\rho^{X}\Big(\{[g(\gamma )]^*\}(i)\Big).
\end{array}$$
$\bullet$ For the general case, assume that $\gamma\!\in\! BC_{\xi}$, with $\xi\! \geq\! 2$. We set 
$B_j\! :=\!\rho^X[\{\gamma^*\}(j)]$, so that we can write $\rho^X(\gamma )\! =\!\bigcup_{j\in\omega}\ \neg B_j$. Note that 
$\{\gamma^*\}(j)\!\in\! BC_{\eta_{j}}$, where $1\!\leq\!\eta_{j}\! <\!\xi$. We set 
$$B_{j,i}\! :=\!\left\{\!\!\!\!\!\!\!\!\begin{array}{ll}
& \neg N\Big( X,[\{\gamma^*\}(j)]^*(i)\Big)~\ \hbox{\rm if}~\ [\{\gamma^*\}(j)](0)\! =\! 0\hbox{\rm ,}\cr & \cr 
& \rho^X\Big(\Big\{[\{\gamma^*\}(j)]^*\Big\}(i)\Big)~\ \hbox{\rm otherwise,}
\end{array}\right.$$
so that $B_j\! :=\!\bigcup_{i\in\omega}\neg B_{j,i}$.

 By Lemma 3.9, $B_{j,i}\! =\!\rho^X\Big[ u^X_{\neg N}\Big([\{\gamma^*\}(j)]^*(i)\Big)\Big]$ if $[\{\gamma^*\}(j)](0)\! =\! 0$, and  
$$\neg B_{j,i}\! =\!\left\{\!\!\!\!\!\!\!\!\begin{array}{ll}
& \rho^X\Big[ u_{N}\Big([\{\gamma^*\}(j)]^*(i)\Big)\Big]~\ \hbox{\rm if}~\ [\{\gamma^*\}(j)](0)\! =\! 0\hbox{\rm ,}\cr & \cr 
& \rho^X\Big[ u_{\neg}\Big(\Big\{[\{\gamma^*\}(j)]^*\Big\}(i)\Big)\Big]~\ \hbox{\rm otherwise,}
\end{array}\right.$$
by Lemma 3.1. Thus
$$\begin{array}{ll}\rho^X(\gamma )\!\!\!\!
& \! =\!\displaystyle{\bigcup_{k\in\omega}}\ \neg B_k\cr
& \! =\!\displaystyle{\bigcup^{\bullet}_{k\in\omega}}\ \bigcap_{j<k}\ B_j\!\setminus\! B_k\cr
& \! =\!\displaystyle{\bigcup^{\bullet}_{k\in\omega}}\ \bigcap_{j<k}\ 
\bigg(\ \bigcup_{i\in\omega}\ \neg B_{j,i}\ \bigg)\!\setminus\! B_k\cr
& \! =\!\displaystyle{\bigcup^{\bullet}_{k\in\omega}}\ \bigcap_{j<k}\ 
\bigg(\ \displaystyle{\bigcup^{\bullet}_{i\in\omega}}\ \bigcap_{l<i}\ B_{j,l}\setminus B_{j,i}\bigg)\!\setminus\! 
B_k\cr
& \! =\!\displaystyle{\bigcup^{\bullet}_{i\in\omega}}\ \bigg(\bigcap_{j<lh(i)}\ \bigcap_{l<(i)_j}\ B_{j,l}\setminus B_{j,(i)_j}\bigg)\!\setminus\! B_{lh(i)}\cr
& \! =\!\displaystyle{\bigcup^{\bullet}_{i\in\omega}}\ \neg
\bigg(B_{lh(i)}\cup\bigcup_{j<lh(i)}\ \bigcup_{l<(i)_j}\ B_{j,(i)_j}\cup\neg B_{j,l}\bigg).
\end{array}$$
Note that the code for $B_{i,j}$ is a partial recursive function of $\gamma ,i$ and $j$. Using Lemma 3.5, this shows the existence of a partial function 
$f^X\! :\!\omega^\omega\!\times\!\omega\!\rightarrow\!\omega^\omega$, recursive on its domain, such that $f^X(\gamma ,i)$ is in $\bigcup_{1\leq\eta <\xi}\ BC_{\eta}$ and  
$\rho^X(\gamma )\! =\!\displaystyle{\bigcup^{\bullet}_{i\in\omega}}\ \neg\rho^X[f^X(\gamma ,i)]$ for each 
$\gamma\!\in\! BC_{\xi}$ with $\xi\! \geq\! 2$. There is $\varepsilon_1\!\in\!\Boraone$ such that  
$f^X(\gamma ,i)\! =\! \{\varepsilon_1\}(\gamma ,i)\! =\!
\{ S^{\omega^\omega ,\omega ,\omega^\omega}_{\Boraone} (\varepsilon_1,\gamma )\}(i)$ 
if $f^X(\gamma ,i)$ is defined. We define 
$h\! :\!\omega^\omega\!\rightarrow\!\omega^\omega$ by the formula $h(\gamma )\! :=\! 
1^\frown S^{\omega^\omega ,\omega ,\omega^\omega}_{\Boraone} 
(\varepsilon_1,\gamma )$. The map $h$ is recursive, $h(\gamma )\!\in\! BC_{\xi}$ and 
$$\rho^X(\gamma )\! =\!\displaystyle{\bigcup^{\bullet}_{i\in\omega}}\ 
\neg\rho^X\Big(\{[h(\gamma )]^*\}(i)\Big)$$ 
if $\gamma\!\in\! BC_{\xi}$ and $\xi\! \geq\! 2$.\bigskip

\noindent $\bullet$ It remains to set $u^X_d(\gamma )\! :=\! g(\gamma)$ if $\gamma (0)\! =\! 0$, 
$h(\gamma )$ otherwise.\hfill{$\square$}\bigskip

 Now we will show that the hyperarithmetical hierarchy makes sense, i.e., the existence of sets of 
arbitrary complexity under $\omega_1^{CK}$. The intuition is quite simple: we take universal sets. 
But we have to check that this is effective.

\vfill\eject

\noindent\bf Notation.\rm\ Recall that if $\alpha\!\in\!\omega^\omega$, then 
$\leq_{\alpha}:=\!\{ (m,n)\!\in\!\omega^{2}\mid\alpha (<\! m,n\! >)\! =\! 1\}$ (see 4A 
in \cite{Moschovakis}), and $<_{\alpha}:=\!\{ (m,n)\!\in\!\omega^{2}\mid\alpha (<\! m,n\! >)\! =\! 1
\ \ \hbox{\rm and}\ \ \alpha (<\! n,m\! >)\!\not=\! 1\}$. The first relation is used to 
define the set $WO\! :=\!\left\{\alpha\!\in\!\omega^\omega\mid\ \leq_{\alpha}\ \hbox{\rm is\ a
\ wellordering\ on\ its\ domain}\ \{ n\!\in\!\omega\mid n\!\leq_\alpha n\}\right\}$, which is used to define   
$$\omega^{CK}_{1} :=\hbox{\rm sup}\{\vert\alpha\vert\mid\alpha\!\in\! WO\cap
\Boraone\}\hbox{\rm ,}$$ 
where $\vert\alpha\vert$ is the order type of $\leq_{\alpha}$. The ordinal $\omega^{CK}_{1}$ is the first non recursive ordinal. If $\alpha\!\in\!\omega^\omega$ and 
$p\!\in\!\omega$, then we define $\alpha_{\vert p}\!\in\! 2^\omega\!\subseteq\!\omega^\omega$ by   
$$\alpha_{\vert p}(q)\! =\! 1\ \Leftrightarrow\ \hbox{\rm Seq}(q)~\ \hbox{\rm and}~\ \hbox{\rm lh}(q)\! =\! 2~\ 
\hbox{\rm and}~\ \alpha (q)\! =\! 1~\ \hbox{\rm and}~\ \forall i\!\in\! 2~\ (q)_i\! <\!_\alpha\ p.$$
If $\alpha\!\in\! WO$, then $\alpha_{\vert p}\!\in\! WO$ and $\leq_{\alpha_{\vert p}}$ is the restriction of 
$\leq_{\alpha}$ to the strict $\leq_{\alpha}$-predecessors of $p$. The next lemma expresses the fact that one can find cofinal sequences of ordinals recursively.

\begin {lem} There is a partial function 
$\eta\! :\!\omega^\omega\!\times\!\omega\!\rightarrow\!\omega^\omega$, recursive on its domain, 
defined if $\alpha\!\in\! WO$ and $\vert\alpha_{\vert p}\vert\!\geq\! 1$, such that 
$\vert\alpha_{\vert p}\vert\! =\!\hbox{\it sup}\!\uparrow\! 
\{\vert\alpha_{\vert\eta(\alpha ,p)(n)}\vert\! +\! 1\mid n\!\in\!\omega\}$.\end{lem}

\noindent\bf Proof.\rm\ This is an application of Kleene's Recursion Theorem. We define a partial function $g\! :\!\omega^\omega\!\times\!\omega\!\rightarrow\!\omega$ by $g(\alpha ,p)\! :=
\!\hbox{\rm min}\{ m\!\in\!\omega\mid m<_{\alpha}p\}$ if it exists. Note that $g$ is recursive on 
its domain and defined on $D\! :=\!\{ (\alpha ,p)\!\in\! WO\!\times\!\omega\mid
\vert\alpha_{\vert p}\vert\!\geq\! 1\}$. We define a map 
$h\! :\!\omega^\omega\!\times\!\omega^{3}\!\rightarrow\!\omega$ by 
$$h(\alpha ,p,n,m)\! :=\!\left\{\!\!\!\!\!\!\!\!
\begin{array}{ll}
& n\ \ \hbox{\rm if}\ \ m<_{\alpha}n<_{\alpha}p\hbox{\rm ,}\cr & \cr 
& m\ \ \hbox{\rm otherwise.}
\end{array}
\right.$$
Note that $h$ is recursive. This allows us to define a partial function 
$\psi\! :\! (\omega^\omega )^{2}\!\times\!\omega^{2}\!\rightarrow\!\omega$ by:
$$\psi (\varepsilon ,\alpha ,p,n)\! :=\!\left\{\!\!\!\!\!\!\!\!
\begin{array}{ll}
& g(\alpha ,p)\ \ \hbox{\rm if}\ \ n\! =\! 0\hbox{\rm ,}\cr & \cr
& h[\alpha ,p,n\! -\! 1,
\{\varepsilon\}^{\omega^\omega\times\omega^{2},\omega}_{\Boraone}
(\alpha ,p,n\! -\! 1)]\ \ \hbox{\rm if}\ \ n\!\geq\! 1.
\end{array}\right.$$
Note that $\psi$ is recursive on its domain, so that there is 
$\varepsilon^{*}\!\in\!\Boraone$ such that 
$\{\varepsilon^{*}\}(\alpha ,p,n)\! =\!\psi (\varepsilon^{*},\alpha ,p,n)$ 
if $\psi (\varepsilon^{*},\alpha ,p,n)$ is defined. Now it is clear 
that $\psi (\varepsilon^{*},\alpha ,p,n)$ is defined if $(\alpha ,p)\!\in\! D$, by 
induction on $n$, and that $\vert\alpha_{\vert p}\vert\! =\!\hbox{\rm sup}
\!\uparrow\! \{\vert\alpha_{\vert \{\varepsilon^{*}\}(\alpha ,p,n)}\vert
\! +\! 1\mid n\!\in\!\omega\}$. We put 
$\eta (\alpha ,p)(n)\! :=\!\{\varepsilon^{*}\}(\alpha ,p,n)$ if 
$\psi (\varepsilon^{*},\alpha ,p,n)$ is defined. Clearly, $\eta$ is defined on $D$ 
and suitable.\hfill{$\square$}\bigskip 

\noindent\bf Notation.\rm\ In the next lemma we identify 
$(\omega^\omega )^\omega$ with $\omega^\omega$, using the formula 
$\Big( (\delta_{q} )_{q\in\omega}\Big)_{n}\! =\!\delta_{n}$. Let 
$\alpha\!\in\! WO$, $\gamma_{0}\!\in\!\omega^\omega$, and 
$u\! :\!\omega^\omega\!\rightarrow\!\omega^\omega$ a map. Using Lemma 3.11 we can 
define, by induction on $p$ (with respect to the wellordering $\leq_{\alpha}$), 
and if $\vert\alpha_{\vert p}\vert\!\geq\! 1$,  
$\gamma_{\vert\alpha_{\vert p}\vert}\! :=\! 
u[ (\gamma_{\vert\alpha_{\vert\eta(\alpha ,p)(n)}\vert})_{n\in\omega}]$. The next lemma expresses the fact that $\gamma_{\vert\alpha_{\vert p}\vert}$ is recursive if the datas are recursive. 

\begin {lem} Let $\theta\! <\!\omega^{CK}_{1}$, $\alpha\!\in\! WO\cap\Boraone$ with 
$\theta\! +\! 1\! =\!\vert\alpha\vert$, $\gamma_{0}\!\in\!\Boraone$,  
$u\! :\!\omega^\omega\!\rightarrow\!\omega^\omega$ a recursive map, and 
$p\!\in\!\omega$ such that $p\!\leq_\alpha\! p$. Then 
$\gamma_{\vert\alpha_{\vert p}\vert}$ is $\Boraone$.\end{lem}

\noindent\bf Proof.\rm\ Once again, this is an application of Kleene's Recursion Theorem. Fix 
$p_{0}\!\in\!\omega$ with $\vert\alpha_{\vert p_0}\vert\! =\! 0$. Using Lemma 3.11, we define a partial function $f\! :\!\omega^\omega\!\times\!\omega\!\rightarrow\!\omega^\omega$ by 
$f(\varepsilon ,p)\! :=\! u\Big[\Big(\{\varepsilon\}[\eta (\alpha ,p)(n)]\Big)_{n\in\omega}\Big]$. 
Note that $f$ is recursive on its domain. We define a partial function 
$\psi\! :\!\omega^\omega\!\times\!\omega\!\rightarrow\!\omega^\omega$ by
$$\psi (\varepsilon ,p)\! :=\!\left\{\!\!\!\!\!\!\!\!
\begin{array}{ll}
& \gamma_{0}\ \ \hbox{\rm if}\ \ p\! =\! p_{0}\hbox{\rm ,}\cr & \cr
& f(\varepsilon ,p)\ \ \hbox{\rm if}\ \ p\!\not=\! p_{0}.
\end{array}\right.$$
As $\psi$ is recursive on its domain, there is $\varepsilon^{*}\!\in\!\Boraone$ 
with $\{\varepsilon^{*}\}(p)\! =\!\psi(\varepsilon^{*},p)$ if 
$\psi(\varepsilon^{*},p)$ is defined. It remains to see that $\psi(\varepsilon^{*},p)$ is defined 
and equal to $\gamma_{\vert\alpha_{\vert p}\vert}$ if $p\leq_{\alpha}p$. We argue by induction on p (with respect to the wellordering $\leq_{\alpha}$). If $p\! =\! p_0$, then 
$\psi(\varepsilon^{*},p)\! =\!\gamma_0\! =\!\gamma_{\vert\alpha_{\vert p_0}\vert}\! =\!\gamma_{\vert\alpha_{\vert p}\vert}$. Assume now that $\vert\alpha_{\vert p}\vert\!\geq\! 1$, and that 
the statement is proved for $q$ satisfying $\vert\alpha_{\vert q}\vert\! <\!\vert\alpha_{\vert p}\vert$. Then 
$\psi[\varepsilon^{*},\eta(\alpha ,p)(n)]$ is defined and equal to 
$\gamma_{\vert\alpha_{\vert \eta(\alpha ,p)(n)}\vert}$ for each $n\!\in\!\omega$. It is also equal to 
$\{\varepsilon^*\}[\eta(\alpha ,p)(n)]$. Thus $f(\varepsilon^*,p)$ is defined and equal to 
$u[ (\gamma_{\vert\alpha_{\vert\eta(\alpha ,p)(n)}\vert})_{n\in\omega}]\! =\!
\gamma_{\vert\alpha_{\vert p}\vert}$.\hfill{$\square$}\bigskip

\noindent\bf Notation.\rm\ In the next lemma, $\alpha\!\in\! WO$ and we study the formula building universal sets for the additive Borel classes. We set 
$\eta_{\alpha ,p}\! :=\!\vert\alpha_{\vert p}\vert$, and 
$\eta_{\alpha ,p,n}\! :=\!\vert\alpha_{\vert \eta(\alpha ,p)(n)}\vert$ 
if $\vert\alpha_{\vert p}\vert\!\geq\! 1$.

\begin {lem} There is $u\! :\!\omega^\omega\!\rightarrow\!\omega^\omega$ 
recursive such that $(\gamma )_{n}\!\in\! BC_{1+\eta_{\alpha ,p,n}}$ for each 
integer $n$ implies that $u(\gamma )\!\in\! BC_{1+\eta_{\alpha ,p}}$ and  
$(\beta ,\delta )\!\in\!\rho^{(2^\omega )^2}[u(\gamma )]\Leftrightarrow 
\exists n\!\in\!\omega ~\ [(\beta )_{n},\delta ]\!\notin\!
\rho^{(2^\omega )^2}[(\gamma )_{n}]$.\end{lem}

\noindent\bf Proof.\rm\ First note that there is $\varepsilon_0\!\in\!\Boraone$ with 
$\{\varepsilon_0\}^{\omega\times (2^\omega )^2,(2^\omega )^2}(n,\beta ,\delta )
\! =\! [(\beta )_{n},\delta ]$ for each $(n,\beta ,\delta )$ in $\omega\!\times\! (2^\omega )^2$. 
Similarly, using Lemmas 3.2 and 3.3, we see that there is 
$\varepsilon_1\!\in\!\Boraone$ such that, for each $(\gamma ,n)$ in 
$\omega^\omega\!\times\!\omega$, $\{\varepsilon_1\}(\gamma ,n)\! =\! 
u^\omega_s\Big( u^{\omega\times (2^\omega )^2,(2^\omega )^2}_r[(\gamma )_n,\varepsilon_0],n\Big)$. We put $u(\gamma )\! :=\! 1^\frown 
S^{\omega^\omega ,\omega ,\omega^\omega}_{\Boraone}(\varepsilon_1,\gamma )$, so that 
$u$ is a recursive map. Moreover, $\{[u(\gamma )]^*\}(n)\! =\! \{\varepsilon_1\}(\gamma ,n)\! =\! 
u^\omega_s\Big( u^{\omega\times (2^\omega )^2,(2^\omega )^2}_r[(\gamma )_n,\varepsilon_0],n\Big)$ is defined and in $BC_{1+\eta_{\alpha ,p,n}}$, 
so that $u(\gamma )\!\in\! BC_{1+\eta_{\alpha ,p}}$. Finally,
$$\begin{array}{ll}(\beta ,\delta )\!\in\!\rho^{(2^\omega )^2}[u(\gamma )]\!\!\!\! 
& \Leftrightarrow\exists  n\!\in\!\omega ~\ (\beta ,\delta )\!\notin\!\rho^{(2^\omega )^2}
[\{[u(\gamma )]^*\}(n)]\cr 
& \Leftrightarrow\exists  n\!\in\!\omega ~\ (\beta ,\delta )\!\notin\!\rho^{(2^\omega )^2}
\Big[u^\omega_s\Big( u^{\omega\times (2^\omega )^2,(2^\omega )^2}_r[(\gamma )_n,\varepsilon_0],n\Big)\Big]\cr 
& \Leftrightarrow\exists  n\!\in\!\omega ~\ (n,\beta ,\delta )\!\notin\!\rho^{\omega\times (2^\omega )^2}
\Big( u^{\omega\times (2^\omega )^2,(2^\omega )^2}_r[(\gamma )_n,\varepsilon_0]\Big)\cr 
& \Leftrightarrow\exists  n\!\in\!\omega ~\ [(\beta )_{n},\delta ]\!\notin\!\rho^{(2^\omega )^2}[(\gamma )_n].
\end{array}$$
This finishes the proof.\hfill{$\square$}

\begin{thm} Let $1\!\leq\!\xi\! <\!\omega^{CK}_{1}$, and ${\it\Gamma}$ be one of the classes 
$\Boraxi$, $\Bormxi$. Then there is 
$B_{\xi}\!\in\! {\it\Gamma}(2^\omega)\!\setminus\!\check {\bf\Gamma}$.\end{thm}

\noindent\bf Proof.\rm\ Assume first that ${\it\Gamma}\! =\!\Boraxi$. As in 22.3 in \cite{Kechris} we set 
$$(\beta ,\delta )\!\in\! {\mathcal U}^{2^\omega}_{\boraone}\ \ \Leftrightarrow
\ \ \exists  k\!\in\!\omega\ \ \beta (k)\! =\! 0\ \ \hbox{\rm and}\ \ \delta\!\in\! N[2^\omega ,k]\hbox{\rm ,}$$
so that ${\mathcal U}^{2^\omega}_{\boraone}\!\in\!\Boraone [(2^\omega )^2]$ is universal for 
$\boraone (2^\omega )$. We define a recursive bijection $\psi\! :\! 2^\omega\!\rightarrow\! (2^\omega)^2$ by $\psi_i (\gamma )(k)\! :=\!\gamma (2k+i)$, for $i\!\in\! 2$. We set 
$B_1\! :=\!\psi^{-1}({\mathcal U}^{2^\omega}_{\boraone})$, so that $B_1$ is $\Boraone$. As in 22.4 in \cite{Kechris}, we see that ${\mathcal U}^{2^\omega}_{\boraone}\!\notin\!\bormone$. Thus 
$B_1\!\notin\!\bormone$ since $\psi$ is a homeomorphism.\bigskip 

 So we may assume that $\xi\!\geq\! 2$ and we will generalize this. Write $\xi\! =\! 1\! +\!\theta$, with 
$1\!\leq\!\theta\! <\!\omega_1^{\hbox{\rm CK}}$. Let $\alpha\!\in\! WO\cap\Boraone$ with 
$\theta\! +\! 1\! =\!\vert\alpha\vert$. Using the previous notation, we get 
 $\eta_{\alpha ,p}\! =\!\hbox{\rm sup}\!\uparrow\! \{\eta_{\alpha ,p,n}\! +\! 1\mid n\!\in\!\omega\}$ if 
 $\eta_{\alpha ,p}\!\geq\! 1$, by Lemma 3.11. As in 22.3 in \cite{Kechris} again we inductively define, if 
 $\eta_{\alpha ,p}\!\geq\! 1$,  
$$(\beta ,\delta )\!\in\! {\mathcal U}^{2^\omega}_{\borapeap}\Leftrightarrow\ \ \exists  n\!\in\!\omega\ \ 
[(\beta )_{n},\delta ]\!\notin\! {\mathcal U}^{2^\omega}_{\borapeapn}\hbox{\rm ,}$$
so that ${\mathcal U}^{2^\omega}_{\borapeap}$ is universal for $\borapeap (2^\omega )$.\bigskip

 Note the existence of $q\!\in\!\omega$ with $\eta_{\alpha ,q}\! =\!\theta$. As before we put $B_{\xi}\! :=\!\psi^{-1}({\mathcal U}^{2^\omega}_{\boraxi})$, so that $B_{\xi}$ is not $\bormxi$. By Lemma 3.3, it remains to see that ${\mathcal U}^{2^\omega}_{\borapeap}$ is $\Borapeap$. By Theorem 3.7 there is 
$\gamma_0\!\in\!\Boraone\cap BC_1$ such that 
${\mathcal U}^{2^\omega}_{\boraone}\! =\!\rho^{(2^\omega )^2}(\gamma_0)$. Lemma 3.13 gives $u$  recursive. We can apply Lemma 3.12, so that $\gamma_{\vert\alpha_{\vert p}\vert}\!\in\!\Boraone$ is defined for each $p$ with $p\!\leq_\alpha\! p$. By induction we see that $\gamma_{\eta_{\alpha ,p}}\!\in\! BC_{1+\eta_{\alpha ,p}}$, by Lemma 3.13. Moreover, 
$$(\beta ,\delta )\!\in\!\rho^{(2^\omega )^2}(\gamma_{\eta_{\alpha ,p}})
\Leftrightarrow
(\beta ,\delta )\!\in\!\rho^{(2^\omega )^2}\Big(u[(\gamma_{\eta_{\alpha ,p,n}})_{n\in\omega}]\Big)\Leftrightarrow 
\exists n\!\in\!\omega ~\ [(\beta )_{n},\delta ]\!\notin\!\rho^{(2^\omega )^2}(\gamma_{\eta_{\alpha ,p,n}}).$$
This inductively shows that 
$\rho^{(2^\omega )^2}(\gamma_{\eta_{\alpha ,p}})\! =\! {\mathcal U}^{2^\omega}_{\borapeap}$. Thus 
${\mathcal U}^{2^\omega}_{\borapeap}$ is $\Borapeap$.\bigskip

 Assume now that ${\it\Gamma}\! =\!\Bormxi$. The previous facts give 
$B_{\xi}\!\in\!\Boraxi (2^\omega )\!\setminus\!\bormxi$. But it is clear that  
$A_{\xi}\! :=\!\neg B_{\xi}$ is in $\Bormxi (2^\omega )\!\setminus\!\boraxi$.\hfill{$\square$}\bigskip 

\noindent\bf Remark.\rm\ We can define, for $\beta\!\in\! 2^\omega$, 
$\omega^{\beta}_{1} :=\hbox{\rm sup}\{\vert\alpha\vert\mid\alpha\!\in\! WO\cap
\Boraone (\beta )\}$. If $X$ is a recursively presented Polish space, then we can define 
$\Boraxi (\beta )(X)\! =\{\rho^{X}(\gamma )\mid\gamma\!\in\!\Boraone (\beta )\cap BC_{\xi}\}$, 
$\Bormxi (\beta )\! :=\!\check\Boraxi (\beta )$ and also $\Borxi (\beta )\! :=\!\Boraxi (\beta )\cap\Bormxi (\beta )$. One can check that this definition of $\Boraone (\beta )$ is equivalent to the one we gave in section 3. The previous proof shows the existence of $B_{\xi}\!\in\!\Boraxi (\beta )(2^\omega)\!\setminus\!\bormxi$, for $1\!\leq\!\xi\! <\!\omega^{\beta}_{1}$. Indeed, the only things to change in the proof are the following. In Lemma 3.12, $\theta\! <\!\omega_1^\beta$, $\alpha\!\in\!\Boraone (\beta )$, $f$ and $\psi$ become 
$\Boraone (\beta )$ on their domain by 3D.7, 3G.1 and 3G.2 in \cite{Moschovakis}. Then we can apply 7A.2 in \cite{Moschovakis} to get $\varepsilon^*$. The conclusion becomes 
$\gamma_{\vert\alpha_{\vert p}\vert}\!\in\!\Boraone (\beta )$. The result follows.

\section{$\!\!\!\!\!$ Effective versions of Kuratowski's theorem.}

\noindent\bf Notation.\rm\ Let $\xi\! <\!\omega_{1}$. Then $\xi\! -\! 1$ will 
denote the predecessor of $\xi$ if it exists, $\xi$ otherwise. We also define $\xi^-\! :=\!\xi\! -\! 1$ if 
$\xi\!\geq\! 3$, $\xi$ otherwise.

\begin{thm} Let $a\!\in\! 2$. There is a partial function 
$F^{a}\! :\!\omega^\omega\!\rightarrow\! (\omega^\omega )^{3}$, recursive on its domain, such that\bigskip\smallskip 

\noindent (a) For each $1\!\leq\!\xi\!\leq\! 2$ and for each $\gamma\!\in\! BC_{\xi}$, coding 
$B\! :=\!\neg\rho^{2^\omega}(\gamma )\!\in\!\bormxi$, $F^0(\gamma )$ is defined and\bigskip 

\noindent (1) $F^0_{0}(\gamma )\!\in\! BC_{1}$ (codes 
$C\! :=\!\neg\rho^{\omega^\omega}[F^0_{0}(\gamma )]\!\in\!\bormone$).\bigskip

\noindent (2) $f\! :=\!\{ F^0_{1}(\gamma )\}^{\omega^\omega ,2^\omega}_{\vert C}$ 
defines a continuous bijection from $C$ onto $B$.\bigskip

\noindent (3) $F^0_{2}(\gamma )\!\in\! BC_{1}$ codes an open set computing a partial function 
$g\! :\! 2^\omega\!\rightarrow\!\omega^\omega$, defined and continuous on $B$, which coincides with 
$f^{-1}$.

\vfill\eject

\noindent (b) For each $1\!\leq\!\xi\! <\!\omega_{1}$ and for each $\gamma\!\in\! BC_{\xi}$, coding 
$B\! :=\!\neg\rho^{2^\omega}(\gamma )\!\in\!\bormxi$, $F^1(\gamma )$ is defined and\bigskip 

\noindent (1) $F^1_{0}(\gamma )\!\in\! BC_{1}$ (codes 
$C\! :=\!\neg\rho^{\omega^\omega}[F^1_{0}(\gamma )]\!\in\!\bormone$).\bigskip

\noindent (2) $f\! :=\!\{ F^1_{1}(\gamma )\}^{\omega^\omega ,2^\omega}_{\vert C}$ 
defines a continuous bijection from $C$ onto $B$.\bigskip

\noindent (3) $F^1_{2}(\gamma )\!\in\! BC_{\xi^-}$ codes a $\boraxim$ set computing a partial function 
$g\! :\! 2^\omega\!\rightarrow\!\omega^\omega$, defined and $\boraxim$-measurable on $B$, which coincides with $f^{-1}$.\end{thm}

\noindent\bf Proof.\rm\ Let us look at the case where $\xi\! =\! 1$ first. We define $\mu\! :\!\omega\!\rightarrow\!\omega$ by 
$$\mu (k)\! :=\!\left\{\!\!\!\!\!\!\!\!
\begin{array}{ll}
& 0~\ \hbox{\rm if}~\ \Big( (k)_{1}\Big) _{1}\! =\! 0\hbox{\rm ,}\cr 
& \hbox{\rm min}\{ l\!\in\!\omega\mid\frac{1}{l+1}\! 
<\!\frac{\Big( (k)_{1}\Big) _{1}}{\Big( (k)_{1}\Big) _{2}+1}\}~\ \hbox{\rm otherwise.}
\end{array}\right.$$ 
Clearly $\mu$ is recursive. Let us recall, for each $k\!\in\!\omega$, the 
definition of the basic neighborhood: 
$$N(\omega^\omega ,k)\! :=\!\left\{\!\!\!\!\!\!\!\!\begin{array}{ll}
& \emptyset\ \ \hbox{\rm if}\ \ \Big( (k)_{1}\Big) _{1}\! =\! 0\hbox{\rm ,}\cr & \cr 
& \Big\{\ \delta\!\in\!\omega^\omega\mid\ \forall j\! <\!\mu (k)\ \ 
\delta (j)\! =\! \bigg( \Big( (k)_{1}\Big) _{0}\bigg)_{j}\ \Big\}\ 
\hbox{\rm otherwise.}
\end{array}\right.$$
In 3A.2 in \cite{Moschovakis} the recursive map $sg\! :\!\omega\!\rightarrow\!\omega$ is defined by $sg(n)\! :=\! 0$ if 
$n\! =\! 0$, $1$ otherwise. The recursive presentation of $2^\omega$ ensures that 
$$N(2^\omega,k)\! :=\!\left\{\!\!\!\!\!\!\!\!\begin{array}{ll}
& \emptyset\ \ \hbox{\rm if}\ \ \Big( (k)_{1}\Big) _{1}\! =\! 0\hbox{\rm ,}\cr 
& \Big\{\ \alpha\!\in\! 2^\omega\mid\ \forall j\! <\!\mu (k)\ \ 
\alpha (j)\! =\! sg\bigg[\bigg( \Big( (k)_{1}\Big) _{0}\bigg)_{j}\bigg]\ \Big\}\ 
\hbox{\rm otherwise.}
\end{array}\right.$$
We view $2^\omega$ as a subset of $\omega^\omega$. We denote by 
$\hbox{\rm Id}_{2^\omega}$ the partial function defined on 
$2^\omega\!\subseteq\!\omega^\omega$, with values in $2^\omega$, by 
$\hbox{\rm Id}_{2^\omega}(\alpha )\! :=\!\alpha$. It is recursive on $2^\omega$, 
since the relation ``$\alpha\!\in\! N(2^\omega,k)$" is 
$\Boraone (\omega^\omega\!\times\!\omega )$ on $2^\omega\!\times\!\omega$. Thus there is 
$\delta_{0}\!\in\!\Boraone$ with $\{\delta_{0}\}^{\omega^\omega ,2^\omega}
(\alpha )\! =\! \hbox{\rm Id}_{2^\omega}(\alpha )$ for each $\alpha\!\in\! 2^\omega$. 
By Lemma 3.3 we have 
$u_{r}^{\omega^\omega ,2^\omega}(\gamma ,\delta_{0})\!\in\! BC_{1}$ and 
$\alpha\!\in\!\rho^{\omega^\omega}[u^{\omega^\omega,2^\omega}_{r}
(\gamma ,\delta_{0})]\Leftrightarrow\alpha\!\in\!\rho^{2^\omega}(\gamma )$ if 
$\gamma\!\in\! BC_{1}$ and $\alpha\!\in\! 2^\omega$. As 
$2^\omega\!\in\!\Bormone (\omega^\omega)$, there is 
$\gamma_{0}\!\in\!\Boraone\cap BC_{1}$ with 
$2^\omega\! =\!\neg\rho^{\omega^\omega}(\gamma_{0})$, by Theorem 3.7. 
We define a recursive map $f\! :\!\omega^\omega\!\rightarrow\!\omega^\omega$ by 
$\Big(f(\gamma )\Big)_{i}\! :=\!\gamma_{0}$ if $i\! =\! 0$, 
$u^{\omega^\omega,2^\omega}_{r}(\gamma ,\delta_{0})$ otherwise.\bigskip

 If $\gamma\!\in\! BC_{1}$, then using Lemma 3.5 we set 
$F^{a}_{0}(\gamma )\! :=\! u_{f}^{\omega^\omega}[1,1,f(\gamma )]$, so that 
$F^{a}_{0}(\gamma )\!\in\! BC_{1}$ and also 
$2^\omega\!\setminus\!\rho^{2^\omega}(\gamma )\! =\!\neg\rho^{\omega^\omega}[F^{a}_{0}(\gamma )]$ since 
$$\rho^{\omega^\omega}[F^{a}_0(\gamma )]\! =\!\bigcup_{i\leq 1}\  \rho^{\omega^\omega}
\Big[\Big( f(\gamma )\Big)_i\Big]\! =\!\rho^{\omega^\omega}(\gamma_0)\cup\rho^{\omega^\omega}
[u_r^{\omega^\omega ,2^\omega}(\gamma ,\delta_0)]
\! =\!\omega^\omega\!\setminus\! 2^\omega\cup\rho^{2^\omega}(\gamma ).$$
Thus $B\! =\! 2^\omega\!\setminus\!\rho^{2^\omega}(\gamma )\!\in\!\bormone (2^\omega )$, and $C\! =\! B$. 
We set $F^{a}_1(\gamma )\! :=\!\delta_0$ if $\gamma\!\in\! BC_1$, so that condition (2) is fullfilled.

\vfill\eject

 We define $P\!\subseteq\! 2^\omega\!\times\!\omega$ by\medskip

\centerline{$(\alpha ,k)\!\in\! P\ \ \Leftrightarrow\ \ \alpha\!\in\! 
N(2^\omega ,k)\ \ \hbox{\rm and}\ \ \bigg[\forall j\! <\!\mu (k)
\ \ \bigg( \Big( (k)_1\Big)_0\bigg)_j\! <\! 2\bigg].$}\medskip

\noindent As $P$ is $\Boraone$, there is $\varepsilon_0\!\in\!\Boraone\cap BC_1$ 
with $P\! =\!\rho^{2^\omega\times\omega}(\varepsilon_0)$, by Theorem 3.7. We put 
$F^{a}_2(\gamma )\! :=\!\varepsilon_0$ if $\gamma\!\in\! BC_1$, so that $F^{a}_2(\gamma )$ codes $P$ computing the canonical injection from $2^\omega$ into $\omega^\omega$ since if 
$\alpha\!\in\! 2^\omega$, then we have $\alpha\!\in\! N(\omega^\omega ,k)\Leftrightarrow P(\alpha ,k)$. So we are done if $\gamma\!\in\! BC_1$.\bigskip

\noindent $\bullet$ For the general case, we give the classical scheme of the 
construction before getting into the effective details, to make things easier to 
understand. So let $B\!\in\!\bormxi$. There is 
$(B_{i})_{i\in\omega}\!\subseteq\!\bigcup_{1\leq\eta <\xi}\ \bormeta$ such that  
$B\! =\!\bigcap_{i\in\omega}\ \neg B_{i}$. Using Lemma 3.10 
we will find $(B_{i,j})_{i,j\in\omega}\!\subseteq\!\bigcup_{1\leq\eta <\xi}\ \boreta$ with 
$\neg B_{i}\! =\!\displaystyle{\bigcup^{\bullet}_{j\in\omega}}\ B_{i,j}$. We will 
argue by induction on $\xi$, so that we will get 
$C_{i,j}\!\in\!\bormone (\omega^\omega )$,  
$f_{i,j}\! :\! C_{i,j}\!\rightarrow\! B_{i,j}$, and $g_{i,j}\! :=\! f_{i,j}^{-1}$. 
The objects we are looking for will be the following:
$$C\! :=\!\Big\{\delta\!\in\!\omega^\omega\mid\ \forall i\!\in\!\omega\ \ 
[(\delta )_{i}]^{*}\!\in\! C_{i,(\delta )_{i}(0)}\ \ \hbox{\rm and}\ \ 
f_{i,(\delta )_{i}(0)}\Big( [(\delta )_{i}]^{*}\Big)\! =\! 
f_{0,(\delta )_{0}(0)}\Big( [(\delta )_{0}]^{*}\Big)\Big\}\hbox{\rm ,}$$
$f(\delta )\! :=\! f_{0,(\delta )_{0}(0)}\Big( [(\delta )_{0}]^{*}\Big)$. 
To define $g$, we define $h\! :\! B\!\rightarrow\!\omega^\omega$ by 
${h(\alpha )(i)\! :=\!\hbox{\rm min}\{ j\!\in\!\omega\mid\alpha\!\in\! B_{i,j}\}}$. Note that $h(\alpha )(i)$ is also the unique integer $j$ satisfying $\alpha\!\in\! B_{i,j}$. We will have 
$\Big( g(\alpha )\Big)_{i}\! :=\! h(\alpha )(i)^\frown g_{i,h(\alpha )(i)}(\alpha )$.\bigskip

\noindent $\bullet$ We set 
$$\begin{array}{ll} (\xi ,\gamma )\!\in\! Q 
& \!\!\Leftrightarrow\ \xi\!\geq\! 2\ \hbox{\rm and}\ \gamma\!\in\! BC_{\xi}\hbox{\rm ,}\cr & \cr
(\xi ,\gamma ,\varepsilon )\!\in\! Q^+ 
& \!\!\Leftrightarrow\ (\xi ,\gamma )\!\in\! Q\ \hbox{\rm and}\ 
\{\varepsilon\}_{2}^{\omega^\omega ,(\omega^\omega )^{3}}(\delta )\ \ 
\hbox{\rm is\ defined\ and\ in}\ BC\ \hbox{\rm for\ each}\ \ \delta\!\in\!\bigcup_{1\leq\eta <\xi}\ BC_{\eta}
\hbox{\rm ,}\cr & \cr
(\xi ,\gamma, \varepsilon,\alpha )\!\in\! Q^{++}
& \!\!\Leftrightarrow\ (\xi ,\gamma ,\varepsilon )\!\in\! Q^+\ \hbox{\rm and}\ \alpha\!\in\! B.
\end{array}$$
Assume that $(\xi ,\gamma )\!\in\! Q$ and $\gamma$ codes $B$, so that 
$\{\gamma^{*}\}(i)$ is defined for each integer $i$, and in $BC_{\eta_i}$ for some $1\!\leq\!\eta_i\! <\!\xi$. Using Lemma 3.10, we set $\gamma_{i,j}\! :=\!\Big\{\Big( u^{2^\omega}_d[\{\gamma^{*}\}(i)]\Big)^*\Big\}(j)$ 
for each $j$. Note that $\gamma_{i,j}$ is recursive in $(\gamma ,i,j)$, $\gamma_{i,j}\!\in\! BC_1$ if $\{\gamma^{*}\}(i)\!\in\! BC_1$, and 
$\gamma_{i,j}\!\in\!\bigcup_{1\leq\eta<\eta_i}\ BC_{\eta}$ if $\{\gamma^{*}\}(i)\!\in\! BC^*$. We also have 
$B_{i,j}\! =\! 2^\omega\!\setminus\!\rho^{2^\omega}(\gamma_{i,j})$.\bigskip

 The map $F^{a}$ will be obtained by Kleene's Recursion Theorem, so that, for some suitable 
$\varepsilon^{a}$, we will have $F^{a}(\gamma )\! =\!\varphi^{a}(\varepsilon^{a}, \gamma )\! =\!
\{\varepsilon^{a}\}^{\omega^\omega ,(\omega^\omega )^{3}}(\gamma )$. In order to describe $\neg C$, we define $R\!\in\!\Boraone [(\omega^\omega )^{3}]$ as follows:
$$\begin{array}{ll}(\varepsilon ,\gamma ,\delta )\!\in\! R\!\!
& \Leftrightarrow\ \exists i\!\in\!\omega\ \Big[\exists j\!\in\!\omega\ \ 
[(\delta )_{i}]^{*}\!\in\! 
N\Big(\omega^\omega ,[\{\varepsilon\}_{0}^{\omega^\omega ,(\omega^\omega )^{3}}
(\gamma_{i,(\delta )_{i}(0)})]^{*}(j)\Big)\Big]\ \ \hbox{\rm or}\cr & \cr 
& \!\!\!\!\Big[\Big\{\{\varepsilon\}_{1}^{\omega^\omega ,(\omega^\omega )^{3}}
(\gamma_{i,(\delta )_{i}(0)})\Big\}^{\omega^\omega ,2^\omega}
\Big( [(\delta )_{i}]^{*}\Big)\!\not=\!
\Big\{\{\varepsilon\}_{1}^{\omega^\omega ,(\omega^\omega )^{3}}
(\gamma_{0,(\delta )_{0}(0)})\Big\}^{\omega^\omega ,2^\omega}
\Big( [(\delta )_{0}]^{*}\Big)\Big].
\end{array}$$
By 3C.4 and 3C.5 in \cite{Moschovakis}, there is 
$R^{*}\!\in\!\Borone [(\omega^\omega )^{2}\!\times\!\omega^{2}]$ with 
$$(\varepsilon ,\gamma ,\delta )\!\in\! R\ \Leftrightarrow\ \exists i\!\in\!\omega\ 
\Big(\delta\!\in\! N[\omega^\omega ,(i)_{0}]\ \ \hbox{\rm and}\ \ 
[\varepsilon ,\gamma ,(i)_{0},(i)_{1}]\!\in\! R^*\Big)$$
(the idea is that an open subset of $(\omega^\omega )^3$ is a countable union of clopen sets)

\vfill\eject
 
 We define a map $\psi_{0}\! :\! (\omega^\omega )^{2}\!\rightarrow\!
\omega^\omega$ by  
$$\psi_{0}(\varepsilon ,\gamma )(i)\! :=\!\left\{\!\!\!\!\!\!\!\!\begin{array}{ll} 
& (i)_{0}\ \ \hbox{\rm if}\ \ [\varepsilon ,\gamma ,(i)_{0},(i)_{1}]\!\in\! R^*
\hbox{\rm ,}\cr & \cr
& 0\ \ \hbox{\rm otherwise.}
\end{array}\right.$$
Clearly, $\psi_{0}$ is recursive and 
$(\varepsilon ,\gamma ,\delta )\!\in\! R\Leftrightarrow\exists i\!\in\!\omega\ \ 
\delta\!\in\! N[\omega^\omega ,\psi_{0}(\varepsilon ,\gamma )(i)]$. We define a recursive map $\varphi_{0}$ by ${\varphi_{0}(\varepsilon ,\gamma )\! :=\! 0^\frown\psi_{0}
(\varepsilon ,\gamma )}$. Note that $\varphi_{0}(\varepsilon ,\gamma )\!\in\! BC_{1}$ (we will have 
$F^{a}_{0}(\gamma )\! =\!\varphi_{0}(\varepsilon^{a} ,\gamma )$, for $\varepsilon^{a}$ 
suitable, if $\gamma\!\in\! BC^{*}$).\bigskip

\noindent $\bullet$ We define a partial function 
$\psi_{1}\! :\! (\omega^\omega )^{3}\!\rightarrow\! 2^\omega$ by 
$$\psi_{1}(\varepsilon ,\gamma ,\delta )\! :=\!
\Big\{\{\varepsilon\}_{1}^{\omega^\omega ,(\omega^\omega )^{3}}
(\gamma_{0,(\delta )_{0}(0)})\Big\}^{\omega^\omega ,2^\omega}
\Big( [(\delta )_{0}]^{*}\Big).$$
As $\psi_{1}$ is recursive on its domain, there is 
$\varepsilon_{1}\!\in\!\Boraone$ such that 
$\psi_{1}(\varepsilon ,\gamma ,\delta )\! =\!
\{\varepsilon_{1}\}^{(\omega^\omega )^{3},2^\omega}
(\varepsilon ,\gamma ,\delta )$ if 
$\psi_{1}(\varepsilon ,\gamma ,\delta )$ is defined. We put 
$\varphi_{1}(\varepsilon ,\gamma )\! :=\! 
S^{(\omega^\omega )^{2},\omega^\omega ,2^\omega}_{\Boraone}
(\varepsilon_{1},\varepsilon ,\gamma )$, so that 
$\psi_{1}(\varepsilon ,\gamma ,\delta )$ is equal to 
$\{\varphi_{1}(\varepsilon ,\gamma )\}^{\omega^\omega ,2^\omega}(\delta )$ when it 
is defined. Note that $\varphi_{1}$ is a total recursive map.\bigskip

\noindent $\bullet$ Now we have to describe $\varphi^{a}_{2}(\varepsilon ,\gamma )$ 
coding a set computing $g$. By the proof of 3C.3 in \cite{Moschovakis} there are recursive maps 
$g'\! :\!\omega\!\rightarrow\!\omega$ and 
$h'\! :\!\omega^{2}\!\rightarrow\!\omega$ such that, for each 
$(\delta ,j,k)\!\in\!\omega^\omega\!\times\!\omega^{2}$, 
$$\delta (j)\! =\! k\ \Leftrightarrow\ \exists i\!\in\!\omega\ \ 
[\delta\!\in\! N(\omega^\omega ,i)\ \ \hbox{\rm and}\ \ 
j\! <\! g'(i)\ \hbox{\rm and}\ h'(i,j)\! =\! k].$$
We set $k_{j}\! :=\!\bigg( \Big( (k)_{1}\Big) _{0}\bigg)_{j}$. We have, for $\alpha\!\in\! B$,  
$$\!\!\!\!\!\!\!\begin{array}{ll}
& \ \ \ g(\alpha )\!\in\! N(\omega^\omega ,k)\cr & \cr  
& \Leftrightarrow \Big( (k)_{1}\Big)_1\!\not=\! 0\ \ \hbox{\rm and}\ \ \forall j\! <\!\mu (k)\ \ 
g(\alpha )(j)\! =\!k_{j}\cr & \cr 
& \Leftrightarrow \Big( (k)_{1}\Big)_1\!\not=\! 0\ \ \hbox{\rm and}\ \ \forall j\! <\!\mu (k)\ \ 
[g(\alpha )]_{\bf (j)_{0}}[{\bf (j)_{1}}]\! =\!k_{j}\cr & \cr  
& \Leftrightarrow \Big( (k)_{1}\Big)_1\!\not=\! 0\ \ \hbox{\rm and}\ \ \forall j\! <\!\mu (k)\ \ 
\bigg[\Big( {\bf (j)_{1}}\! =\! 0\ \ \hbox{\rm and}\ \ h(\alpha )[{\bf (j)_{0}}]\! =\! k_{j}
\Big)\ \hbox{\rm or}\cr & \cr 
& \ \ \ \ \   
\Big( {\bf (j)_{1}}\! >\! 0\ \ \hbox{\rm and}\ \ g_{{\bf (j)_{0}},h(\alpha )[{\bf (j)_{0}}]}
(\alpha )[{\bf (j)_{1}}\! -\! 1]\! =\! k_{j}\Big)\bigg]\cr & \cr  
& \Leftrightarrow \Big( (k)_{1}\Big)_1\!\not=\! 0\ \ \hbox{\rm and}\ \ \forall j\! <\!\mu (k)\ \ 
\bigg[\Big( {\bf (j)_{1}}\! =\! 0\ \ \hbox{\rm and}\ \ 
\alpha\!\in\! B_{{\bf (j)_{0}},k_{j}}\Big)\ \ \hbox{\rm or}\cr & \cr  
& \ \ \ 
\Big( {\bf (j)_{1}}\! >\! 0\ \ \hbox{\rm and}\ \ \exists i\!\in\!\omega\ 
\Big[\ g_{{\bf (j)_{0}},h(\alpha )[{\bf (j)_{0}}]}(\alpha )\!\in\! N(\omega^\omega ,i)\ 
\hbox{\rm and}\ {\bf (j)_{1}}\!\leq\! g'(i)\ \hbox{\rm and}\ 
h'[i,{\bf (j)_{1}}\! -\! 1]\! =\! k_{j}\ \Big]\Big)\bigg].
\end{array}$$

 But $\alpha\!\in\! B_{{\bf (j)_{0}},k_{j}}\ 
\Leftrightarrow\ \alpha\!\notin\!\rho^{2^\omega}(\gamma_{{\bf (j)_{0}},k_{j}})\ 
\Leftrightarrow\ \exists l\!\in\!\omega\  [k_{j}\! =\! l\ \ 
\hbox{\rm and}\ \ \alpha\!\notin\!\rho^{2^\omega}(\gamma_{{\bf (j)_{0}},l})]$. There is $\delta_{1}\!\in\!\Boraone$ such that $\{\delta_{1}\}^{2^\omega\times\omega ,2^\omega}(\alpha ,k)\! =\!\alpha$ if 
$(\alpha ,k)\!\in\! 2^\omega\!\times\!\omega$. We will code the relation 
``$\tilde R_0(\alpha ,k)\Leftrightarrow\alpha\!\in\! B_{{\bf (j)_{0}},l}$", via a partial function 
$\tilde g_{0}\! :\!\omega^\omega\!\times\!\omega^2\!\rightarrow\!\omega^\omega$.\bigskip

 If $\{\gamma^*\}[{\bf (j)_0}](0)\! =\! 0$, then by Lemmas 3.10 and 3.3 we get 
$$\begin{array}{ll}\alpha\!\notin\!\rho^{2^\omega}(\gamma_{{\bf (j)_{0}},l})\!\!\!\!
& \Leftrightarrow \alpha\!\in\!\rho^{2^\omega}\Big[u^{2^\omega}_c\Big(\{\gamma^*\}[{\bf (j)_0}],l\Big)\Big]\cr
& \cr
& \Leftrightarrow (\alpha ,k)\!\in\!\rho^{2^\omega\times\omega}
\bigg( u^{2^\omega\times\omega ,2^\omega}_{r}\Big[u^{2^\omega}_c\Big(\{\gamma^*\}[{\bf (j)_0}],l\Big),\delta_{1}\Big]\bigg).
\end{array}$$

 If $\{\gamma^*\}[{\bf (j)_0}](0)\!\not=\! 0$, then by Lemmas 3.3 and 3.1 we get 
$$\begin{array}{ll}\alpha\!\notin\!\rho^{2^\omega}(\gamma_{{\bf (j)_{0}},l})\!\!\!\!
& \Leftrightarrow (\alpha ,k)\!\notin\!\rho^{2^\omega\times\omega}
[u^{2^\omega\times\omega ,2^\omega}_{r}(\gamma_{{\bf (j)_{0}},l},\delta_{1})]\cr 
& \cr
& \Leftrightarrow (\alpha ,k)\!\in\!\rho^{2^\omega\times\omega}\Big( u_{\neg}
[u^{2^\omega\times\omega ,2^\omega}_{r}(\gamma_{{\bf (j)_{0}},l},\delta_{1})]\Big).
\end{array}$$
This shows the existence of a partial function 
$\tilde g_{0}\! :\!\omega^\omega\!\times\!\omega^2\!\rightarrow\!\omega^\omega$, recursive on its domain, such that $\tilde g_{0}(\gamma ,j,l)$ is defined if $(\xi ,\gamma )\!\in\! Q$. In this case,  
$\tilde g_{0}(\gamma ,j,l)\!\in\! BC_{1}$ if $\{\gamma^*\}[{\bf (j)_0}]\!\in\! BC_1$, 
$\tilde g_{0}(\gamma ,j,l)\!\in\! BC_{\eta_{{\bf (j)_{0}}}}$ if $\{\gamma^*\}[{\bf (j)_0}]\!\in\! BC^*$, and  
$$\alpha\!\in\! B_{{\bf (j)_{0}},l}\ \Leftrightarrow\ 
\alpha\!\notin\!\rho^{2^\omega}(\gamma_{{\bf (j)_{0}},l})\ \Leftrightarrow\ 
(\alpha ,k)\!\in\!\rho^{2^\omega\times\omega}[\tilde g_{0}(\gamma ,j,l)].$$
Similarly, we now code the relation 
``$R_0(\alpha ,k)\Leftrightarrow\alpha\!\in\! B_{{\bf (j)_{0}},k_j}$", via a partial function 
$g_{0}\! :\!\omega^\omega\!\times\!\omega\!\rightarrow\!\omega^\omega$. Choose 
$\gamma_{1}\!\in\!\Boraone\cap BC_{1}$ such that 
$(\alpha ,k,j,l)\!\in\!\rho^{2^\omega\times\omega^{3}}(\gamma_{1})\ 
\Leftrightarrow\ k_{j}\! =\! l$. Using Lemma 3.2 we see that 
$u^{\omega^{2}}_{s}(\gamma_{1},j,l)\!\in\! BC_{1}$ and 
$(\alpha ,k)\!\in\!\rho^{2^\omega\times\omega}[u^{\omega^{2}}_{s}(\gamma_{1},j,l)]
\ \Leftrightarrow\ k_{j}\! =\! l$, for each 
$(\alpha ,k,j,l)\!\in\! 2^\omega\!\times\!\omega^{3}$. Using Lemmas 3.5 and 3.6.(b), we get the existence of a partial function 
$g_{0}\! :\!\omega^\omega\!\times\!\omega\!\rightarrow\!\omega^\omega$, recursive on its domain, such that $g_{0}(\gamma ,j)$ is defined if $(\xi ,\gamma )\!\in\! Q$. In this case, 
$g_{0}(\gamma ,j)\!\in\! BC_{1}$ if $\{\gamma^*\}[{\bf (j)_0}]\!\in\! BC_1$, and 
$g_{0}(\gamma ,j)\!\in\! BC_{\eta_{{\bf (j)_{0}}}}$ if $\{\gamma^*\}[{\bf (j)_0}]\!\in\! BC^*$. If moreover 
$\alpha\!\in\! B$, then $\alpha\!\in\! B_{{\bf (j)_{0}},k_{j}}\ \Leftrightarrow\ 
(\alpha ,k)\!\in\!\rho^{2^\omega\times\omega}[g_{0}(\gamma ,j)]$.\bigskip

\noindent - Similarly, we now deal with the end of the computation of the relation 
``$g(\alpha )\!\in\! N(\omega^\omega ,k)$" above. We will have
$$\begin{array}{ll}
g_{{\bf (j)_{0}},h(\alpha )[{\bf (j)_{0}}]}(\alpha )\!\in\! N(\omega^\omega ,i)\!\!\!
& \Leftrightarrow\ (\alpha ,i)\!\in\!\rho^{2^\omega\times\omega}
[\{\varepsilon\}_{2}^{\omega^\omega ,(\omega^\omega )^{3}}(\gamma_{{\bf (j)_{0}},h(\alpha )
[{\bf (j)_0}]})]\cr & \cr 
& \Leftrightarrow\ \exists l\!\in\!\omega\ \Big[\ (\alpha ,i)\!\in\!
\rho^{2^\omega\times\omega}
[\{\varepsilon\}_{2}^{\omega^\omega ,(\omega^\omega )^{3}}(\gamma_{{\bf (j)_{0}},l})]\ 
\hbox{\rm and}\ h(\alpha )[{\bf (j)_0}]\! =\! l\ \Big]\cr & \cr 
& \Leftrightarrow\ \exists l\!\in\!\omega\ \Big[\ (\alpha ,i)\!\in\!\rho^{2^\omega\times\omega}
[\{\varepsilon\}_{2}^{\omega^\omega ,(\omega^\omega )^{3}}(\gamma_{{\bf (j)_{0}},l})]\ 
\hbox{\rm and}\ \alpha\!\in\! B_{{\bf (j)_0},l}\ \Big]
\end{array}$$
if $(\xi ,\gamma, \varepsilon ,\alpha )\!\in\! Q^{++}$. If we apply Lemmas 3.5 and 3.6.(a), then we obtain the existence of a partial function $g^0_{1}\! :\! (\omega^\omega )^2\!\times\!\omega^2\!\rightarrow\!\omega^\omega$, recursive on its domain, such that $g^0_{1}(\varepsilon ,\gamma ,j,l)$ is defined and in $BC$ if $(\xi ,\gamma ,\varepsilon )\!\in\! Q^+$, in which case  
$(\alpha ,k)\!\in\!\rho^{2^\omega\times\omega}[g^0_{1}(\varepsilon ,\gamma ,j,l)]$ is equivalent to 
$$\exists i\!\in\!\omega\ \Big[\ (\alpha ,i)\!\in\!\rho^{2^\omega\times\omega}
[\{\varepsilon\}_{2}^{\omega^\omega ,(\omega^\omega )^{3}}(\gamma_{{\bf (j)_{0}},l})]\ \ 
\hbox{\rm and}\ \ {\bf (j)_{1}}\!\leq\! g'(i)\ \ \hbox{\rm and}\ \ h'[i,{\bf (j)_{1}}\! -\! 1]\! =\! k_{j}\ \Big]\hbox{\rm .}$$

 If $(\xi ,\gamma, \varepsilon ,\alpha )\!\in\! Q^{++}$, then  
$$\begin{array}{ll}
& \ \ \ \ \ \ \exists i\!\in\!\omega\ \Big[\ g_{{\bf (j)_{0}},h(\alpha )[{\bf (j)_{0}}]}(\alpha )\!\in\! N(\omega^\omega ,i)\ 
\ \hbox{\rm and}\ \ {\bf (j)_{1}}\!\leq\! g'(i)\ \ \hbox{\rm and}\ \ h'[i,{\bf (j)_{1}}\! -\! 1]\! =\! k_{j}\ \Big]\cr & \cr
& \Leftrightarrow\ \exists l\!\in\!\omega\ \Big[\ (\alpha ,k)\!\in\!\rho^{2^\omega\times\omega}
[g^0_{1}(\varepsilon ,\gamma ,j,l)]\ \ \hbox{\rm and}\ \ \alpha\!\in\! B_{{\bf (j)_{0}},l}\ \Big].
\end{array}$$

 But $g^0_{1}(\varepsilon ,\gamma ,j,l)$ could be in $BC_1$ for some $l$'s, and in 
$BC^*$ for some others. This may happen if $\xi\!\geq\! 3$. This is a problem since we want to apply Lemma 3.6.(b). We will solve this problem with Lemma 3.4. 
We define a partial function $g^1_{1}\! :\! (\omega^\omega )^2\!\times\!\omega^2\!\rightarrow\!\omega^\omega$ by 
$g^1_{1}(\varepsilon ,\gamma ,j,l)\! :=\! u^{2^\omega\times\omega}_*[g^0_{1}(\varepsilon ,\gamma ,j,l)]$. As $\rho^{2^\omega\times\omega}[g^1_{1}(\varepsilon ,\gamma ,j,l)]\! =\!\rho^{2^\omega\times\omega}[g^0_{1}(\varepsilon ,\gamma ,j,l)]$ if 
$(\xi ,\gamma ,\varepsilon )\!\in\! Q^+$, it satisfies the previous properties of $g^0_1$.\bigskip

\noindent - Lemmas 3.5 and 3.6.(b) imply the existence of a partial function 
$g^1_{2}\! :\! (\omega^\omega )^2\!\times\!\omega\!\rightarrow\!\omega^\omega$, 
recursive on its domain, such that $g^1_{2}(\varepsilon ,\gamma ,j)$ is defined and in $BC$ if 
$(\xi ,\gamma, \varepsilon )\!\in\! Q^+$. If moreover $(\xi ,\gamma, \varepsilon ,\alpha)\!\in\! Q^{++}$, then 
$(\alpha ,k)\!\in\!\rho^{2^\omega\times\omega}[g^1_{2}(\varepsilon ,\gamma ,j)]$ is equivalent to 
$$\exists i\!\in\!\omega\ \Big[\ g_{{\bf (j)_{0}},h(\alpha )[{\bf (j)_{0}}]}(\alpha )\!\in\! N(\omega^\omega ,i)\ 
\ \hbox{\rm and}\ \ {\bf (j)_{1}}\!\leq\! g'(i)\ \ \hbox{\rm and}\ \ h'[i,{\bf (j)_{1}}\! -\! 1]\! =\! k_{j}\ \Big].$$
We also define a partial function $g^0_{2}\! :\! (\omega^\omega )^2\!\times\!\omega\!\rightarrow\!\omega^\omega$. It is defined relatively to $g^0_1$ exactly like $g^1_2$ was defined relatively to $g^1_1$. It will satisfy the previous properties of $g^1_2$ if $\xi\! =\! 2$, and we will have, for 
$(\xi ,\gamma, \varepsilon ,\alpha)\!\in\! Q^{++}$ and $a\!\in\! 2$,
$$\begin{array}{ll}
g(\alpha )\!\in\! N(\omega^\omega ,k)\ \Leftrightarrow\ \ \Big( (k)_{1}\Big)_1\!\not=\! 0\ \ \hbox{\rm and}\ \ 
\forall j\! <\!\mu (k)\!\!\!
& \Big[\Big( {\bf (j)_1}\! =\! 0\ \ \hbox{\rm and}\ \ (\alpha ,k)\in\rho^{2^\omega\times\omega}
[g_{0}(\gamma ,j)]\Big)\ \hbox{\rm or}\cr & \cr
& \Big( {\bf (j)_1}\!  >\! 0\ \ \hbox{\rm and}\ \ (\alpha ,k)\in\rho^{2^\omega\times\omega}
[g^{a}_{2}(\varepsilon ,\gamma ,j)]\Big)\Big].
\end{array}$$
- We define a partial function 
$g^1_3\! :\! (\omega^\omega )^2\!\times\!\omega\!\rightarrow\!\omega^\omega$ by
$$g^1_3(\varepsilon ,\gamma ,j)\! :=\!\left\{\!\!\!\!\!\!\!\!
\begin{array}{ll}
& u^{2^\omega\times\omega}_*[g_{0}(\gamma ,j)]\ \ \hbox{\rm if}\ \ {\bf (j)_1}\! =\! 0\hbox{\rm ,}\cr & \cr 
& g^1_{2}(\varepsilon ,\gamma ,j)\ \ \hbox{\rm if}\ \ {\bf (j)_1}\!  >\! 0.
\end{array}\right.$$
Note that $g^1_{3}$ is recursive on its domain, and $g^1_3(\varepsilon ,\gamma ,j)$ is defined and in $BC$ if $(\xi ,\gamma, \varepsilon )\!\in\! Q^+$. We also define a partial function 
$g^0_3\! :\! (\omega^\omega )^2\!\times\!\omega\!\rightarrow\!\omega^\omega$. It is defined relatively to $g^0_2$, like $g^1_3$ was defined relatively to $g^1_2$, except that 
$g^0_3(\varepsilon ,\gamma ,j)\! :=\! g_{0}(\gamma ,j)$ if $ {\bf (j)_1}\! =\! 0$. The function $g^0_3$ will satisfy the previous properties of $g^1_3$ if $\xi\! =\! 2$.\bigskip

\noindent - By Lemma 3.5, we get the existence of a partial function 
$g^1_4\! :\! (\omega^\omega )^2\!\times\!\omega\!\rightarrow\!\omega^\omega$, recursive on its domain, such that $g^1_4(\varepsilon ,\gamma ,m)$ is defined and in $BC$ if 
$(\xi ,\gamma, \varepsilon )\!\in\! Q^+$ and, if moreover $(\xi ,\gamma, \varepsilon ,\alpha )\!\in\! Q^{++}$, then 
$$(\alpha ,k)\!\in\!\rho^{2^\omega\times\omega}[g^1_{4}(\varepsilon ,\gamma ,m)]\ \ \Leftrightarrow\ \ 
\Big( (k)_{1}\Big)_1\!\not=\! 0\ \ \hbox{\rm and}\ \ \mu (k)\! =\! m\ \ \hbox{\rm and}\ \ 
\forall j\! <\!m\ \ (\alpha ,k)\!\in\!\rho^{2^\omega\times\omega}[g^1_{3}(\varepsilon ,\gamma ,j)].$$
Thus $(\xi ,\gamma, \varepsilon ,\alpha )\!\in\! Q^{++}$ will imply that 
$$\begin{array}{ll}
g(\alpha )\!\in\! N(\omega^\omega ,k)\!\!
& \Leftrightarrow\ \Big( (k)_{1}\Big)_1\!\not=\! 0\ \ \hbox{\rm and}\ \ \forall j\! <\!\mu (k)\ \ (\alpha ,k)\!\in\!
\rho^{2^\omega\times\omega}[g^1_{3}(\varepsilon ,\gamma ,j)]\cr & \cr 
& \Leftrightarrow\ \exists m\!\in\!\omega\ \Big[\Big( (k)_{1}\Big)_1\!\not=\! 0\ \ \hbox{\rm and}\ \ 
\mu (k)\! =\! m\ \ \hbox{\rm and}\ \ \forall j\! <\!m\ \ (\alpha ,k)\!\in\!
\rho^{2^\omega\times\omega}[g^1_{3}(\varepsilon ,\gamma ,j)]\Big]\cr & \cr 
& \Leftrightarrow\ \exists m\!\in\!\omega\ \ 
(\alpha ,k)\!\in\!\rho^{2^\omega\times\omega}[g^1_{4}(\varepsilon ,\gamma ,m)].
\end{array}$$

 We also define a partial function 
$g^0_4\! :\! (\omega^\omega )^2\!\times\!\omega\!\rightarrow\!\omega^\omega$. It is defined relatively to $g^0_3$ exactly like $g^1_4$ was defined relatively to $g^1_3$. It will satisfy the previous properties of $g^1_4$ if $\xi\! =\! 2$.\bigskip

\noindent - By Lemma 3.6.(b), we get the existence of a partial function 
$\varphi^1_2\! :\! (\omega^\omega )^2\!\rightarrow\!\omega^\omega$, 
recursive on its domain, such that $\varphi^1_2(\varepsilon ,\gamma )$ is defined and in $BC$ if 
$(\xi ,\gamma, \varepsilon )\!\in\! Q^+$, and
$$g(\alpha )\!\in\! N(\omega^\omega ,k)\ \Leftrightarrow\ 
(\alpha ,k)\!\in\!\rho^{2^\omega\times\omega}[\varphi^1_2(\varepsilon ,\gamma )]$$
if $(\xi ,\gamma, \varepsilon ,\alpha)\!\in\! Q^{++}$. We also define a partial function 
$\varphi^0_2\! :\! (\omega^\omega )^2\!\rightarrow\!\omega^\omega$. It is defined relatively to $g^0_4$ exactly the way $\varphi^1_2$ was defined relatively to $g^1_4$. It will satisfy the previous properties of 
$\varphi^1_2$ if $\xi\! =\! 2$.\bigskip

\noindent $\bullet$ Now we can define a partial function  
$\varphi^{a}\! :\! (\omega^\omega )^2\!\rightarrow\! (\omega^\omega )^{3}$ by 
$$\varphi^{a}(\varepsilon ,\gamma )\! :=\!\left\{\!\!\!\!\!\!\!\!
\begin{array}{ll}
& [F^{a}_{0}(\gamma ),F^{a}_{1}(\gamma ),F^{a}_{2}(\gamma )]\ \ \hbox{\rm if}\ \ 
\gamma (0)\! =\! 0\hbox{\rm ,}\cr & \cr 
& [\varphi_{0} (\varepsilon ,\gamma ),\varphi_{1} (\varepsilon ,\gamma ),\varphi^{a}_{2} 
(\varepsilon ,\gamma )]\ \ \hbox{\rm if}\ \ \gamma (0)\!\not=\! 0.
\end{array}\right.$$
As $\varphi^{a}$ is recursive on its domain, by Kleene's Recursion Theorem there is 
$\varepsilon^{a}\!\in\!\Boraone$ such that 
$$\{\varepsilon^{a}\}^{\omega^\omega ,(\omega^\omega )^{3}}(\gamma )\! =\!
\varphi^{a} (\varepsilon^{a},\gamma )$$ 
if $\varphi^{a}(\varepsilon^{a},\gamma )$ is defined. We define a partial 
function $F^{a}\! :\!\omega^\omega\!\rightarrow\! (\omega^\omega )^{3}$ by 
$F^{a} (\gamma )\! :=\!\{\varepsilon^{a}\}^{\omega^\omega ,(\omega^\omega )^{3}}
(\gamma )$, so that $F^{a}$ is recursive on its domain. We already checked that $F^{a}(\gamma )$ is suitable if $\gamma\!\in\! BC_{1}$.\bigskip 

 So assume that $2\!\leq\!\xi\! <\!\omega_{1}$, and $\gamma\!\in\! BC_{\xi}$ codes 
$B\! :=\!\neg\rho^{2^\omega}(\gamma )\!\in\!\bormxi$. We will prove that 
$F^{a}(\gamma )$ is defined and fullfills the required properties by induction on 
$\xi$.\bigskip

 Note that $\{\varepsilon^{a}\}_2^{\omega^\omega ,(\omega^\omega )^{3}}(\delta )$ is defined and in $BC$ for each $\delta\!\in\!\bigcup_{1\leq\eta <\xi}\ BC_{\eta}$, by induction assumption. This implies that $(\xi ,\gamma,\varepsilon^{a})\!\in\! Q^+$, $\varphi^{a}_{2} (\varepsilon^{a},\gamma )$ and $F^{a}(\gamma )$ are defined, and 
$$F^{a}(\gamma )\! =\!\{\varepsilon^{a}\}^{\omega^\omega ,(\omega^\omega )^{3}}(\gamma )
\! =\!\varphi^{a}(\varepsilon^{a},\gamma )\! =\! [\varphi_{0} (\varepsilon^{a} ,\gamma ),
\varphi_{1} (\varepsilon^{a} ,\gamma ),\varphi^{a}_{2} (\varepsilon^{a},\gamma )].$$
(1) Note that $F^{a}_{0}(\gamma )\!\in\! BC_{1}$ since 
$F^{a}_{0}(\gamma )\! =\!\varphi_{0}(\varepsilon^{a},\gamma )\! =\! 0^\frown\psi_{0}(\varepsilon^{a},
\gamma )$. Moreover, with the previous notation, we get $\delta\!\notin\!\rho^{\omega^\omega}
[F^{a}_0(\gamma )]\ \Leftrightarrow\ (\varepsilon^{a} ,\gamma ,\delta )\!\notin\! R\ \Leftrightarrow\ \delta
\!\in\! C$, by induction assumption.\bigskip

\noindent (2) We have $F^{a}_{1}(\gamma )\! 
=\!\varphi_{1}(\varepsilon^{a},\gamma )$, so that, by induction assumption, and for each 
$\delta\!\in\! C$, 
$$\begin{array}{ll}\{ F^{a}_{1}(\gamma )\}^{\omega^\omega ,2^\omega}(\delta )\!\!\!\!
&  =\psi_{1}(\varepsilon^{a},\gamma ,\delta )\cr & \cr
&  =\Big\{\{\varepsilon^{a}\}_{1}^{\omega^\omega ,(\omega^\omega )^{3}}
(\gamma_{0,(\delta )_{0}(0)})\Big\}^{\omega^\omega ,2^\omega}
\Big( [(\delta )_{0}]^{*}\Big)\cr & \cr
&  =\Big\{ F^{a}_{1}(\gamma_{0,(\delta )_{0}(0)})\Big\}^{\omega^\omega ,2^\omega}
\Big( [(\delta )_{0}]^{*}\Big)\cr & \cr
&  = f_{0,(\delta )_{0}(0)}\Big( [(\delta )_{0}]^{*}\Big)\! =\! f(\delta ).
\end{array}$$ 

 Clearly, $f$ is continuous. If $\delta\!\in\! C$ and $i\!\in\!\omega$, then 
$f(\delta )\! =\! f_{i,(\delta )_{i}(0)}\Big( [(\delta )_{i}]^{*}\Big)\!\in\! 
B_{i,(\delta )_{i}(0)}$, thus $f(\delta )\!\in\! B$. Let $\delta$, $\delta'\!\in\! C$ such that 
$\alpha\! :=\! f(\delta )\! =\! f(\delta' )$. Then 
$\alpha\! =\! f_{i,(\delta )_{i}(0)}\Big( [(\delta )_{i}]^{*}\Big)\!\in\! B_{i,(\delta )_{i}(0)}$, 
so that $\alpha$ is in $B_{i,(\delta )_{i}(0)}$ and $B_{i,(\delta' )_{i}(0)}$. This shows that 
$(\delta )_{i}(0)\! =\! (\delta' )_{i}(0)$. Thus $[(\delta )_{i}]^{*}\! 
=\! [(\delta' )_{i}]^{*}$ since $f_{i,(\delta )_{i}(0)}$ is one-to-one, 
$(\delta )_{i}\! =\! (\delta' )_{i}$, and 
$\delta\! =\!\delta'$. This shows that $f$ is one-to-one. If $\alpha\!\in\! B$ and 
$i\!\in\!\omega$, then there is a unique integer $j_{i}$ with $\alpha\!\in\! B_{i,j_{i}}$. There 
is $\delta^{i}\!\in\! C_{i,j_{i}}$ with $\alpha\! =\! f_{i,j_{i}}(\delta^{i})$. Put 
$(\delta )_{i}\! :=\! {j_{i}}^\frown \delta^{i}$. Then $\delta\!\in\! C$ and 
$\alpha\! =\! f(\delta )$. This shows that $f$ is onto.\bigskip

\noindent (3) We have $F^{a}_{2}(\gamma )\! =\!\varphi^{a}_{2}(\varepsilon^{a},\gamma )$.\bigskip 

\noindent - If $\xi\! =\! 2$, then $\eta_i\! =\! 1$ for each $i$, $\gamma_{i,j}\!\in\! BC_1$ for each 
$(i,j)$. Thus\bigskip

 $\circ$ $g^0_{1}(\varepsilon^{a},\gamma ,j,l)\!\in\! BC_1$, by induction assumption, since 
$\{\varepsilon^{a}\}_{2}^{\omega^\omega ,(\omega^\omega )^{3}}(\gamma_{{\bf (j)_{0}},l})\!\in\! BC_1$. This implies that $g^1_{1}(\varepsilon^{1},\gamma ,j,l)\!\in\! BC_2$.\bigskip

 $\circ$ $g^1_{2}(\varepsilon^{1},\gamma ,j)\!\in\! BC_2$ and 
$g^0_{2}(\varepsilon^{a},\gamma ,j)\!\in\! BC_1$.\bigskip

 $\circ$ $g^1_{3}(\varepsilon^{1},\gamma ,j)\!\in\! BC_2$ and 
$g^0_{3}(\varepsilon^{a},\gamma ,j)\!\in\! BC_1$.\bigskip

 $\circ$ $g^1_{4}(\varepsilon^{1},\gamma ,m)\!\in\! BC_2$ and 
$g^0_{4}(\varepsilon^{a},\gamma ,m)\!\in\! BC_1$.\bigskip

 $\circ$ $\varphi^1_2(\varepsilon^1,\gamma )\!\in\! BC_2\! =\! BC_{\xi^-}$ and 
$\varphi^0_2(\varepsilon^{a},\gamma )\!\in\! BC_1$.\bigskip

\noindent - If $\xi\!\geq\! 3$, then\bigskip

 $\circ$ $g^1_{1}(\varepsilon^{1},\gamma ,j,l)\!\in\! BC_{\hbox{\rm max}(2,\eta_{{\bf (j)_{0}}})}$, by induction assumption. Indeed, if $\gamma_{{\bf (j)_{0}},l}\!\in\! BC_1$, then 
$g^0_{1}(\varepsilon^{1},\gamma ,j,l)$ is in $BC_1$ and $g^1_{1}(\varepsilon^{1},\gamma ,j,l)\!\in\! BC_2$. If $\gamma_{{\bf (j)_{0}},l}\!\in\! BC^*$, then $\gamma_{{\bf (j)_{0}},l}\!\in\! BC_{\eta_{{\bf (j)_{0}}}}$, and 
$g^0_{1}(\varepsilon^{1},\gamma ,j,l)$, $g^1_{1}(\varepsilon^{1},\gamma ,j,l)$ too.\bigskip

 $\circ$ $g^1_{2}(\varepsilon^{1},\gamma ,j)\!\in\! BC_{\hbox{\rm max}(2,\eta_{{\bf (j)_{0}}})}$.\bigskip

 $\circ$ $g^1_{3}(\varepsilon^{1},\gamma ,j)\!\in\! BC_{\hbox{\rm max}(2,\eta_{{\bf (j)_{0}}})}$.\bigskip

 $\circ$ $g^1_{4}(\varepsilon^{1},\gamma ,m)\!\in\! BC_{\hbox{\rm max}_{j<m}(2,\eta_{{\bf (j)_{0}}})}\!\subseteq\! BC_{\xi -1}$.\bigskip

 $\circ$ $\varphi^1_2(\varepsilon^1,\gamma )\!\in\! BC_{\xi -1}\! =\! BC_{\xi^-}$.\bigskip
 
\noindent Thus $F^{0}_{2}(\gamma )\! =\!\varphi^0_2(\varepsilon^{0},\gamma )\!\in\! BC_{1}$ if 
$\gamma\!\in\! BC_2$, and 
$F^{1}_{2}(\gamma )\! =\!\varphi^1_2(\varepsilon^{1},\gamma )\!\in\! BC_{\xi^-}$. And  $\rho^{2^\omega\times\omega}[F^{a}_2(\gamma )]$ computes $g$ on $B$. If $\alpha\!\in\! B$ and $i\!\in\!\omega$, then $\Big[\Big( g(\alpha )\Big)_{i}\Big]^{*}\! =\! 
g_{i,h(\alpha )(i)}(\alpha )\!\in\! C_{i,h(\alpha )(i)}$ since $\alpha\!\in\! B_{i,h(\alpha )(i)}$. 
Thus 
$$f_{i,h(\alpha )(i)}[g_{i,h(\alpha )(i)}(\alpha )]\! =\!\alpha$$ 
and $g(\alpha )\!\in\! C$ since $\Big( g(\alpha )\Big)_{i}(0)\! =\! h(\alpha )(i)$. Moreover, 
$$f[g(\alpha )]\! =\! f_{0,( g(\alpha ))_{0}(0)}\bigg(\Big[\Big( g(\alpha )\Big)_{0}\Big]^{*}\bigg)
\! =\!\alpha .$$ 
If $\delta\!\in\! C$ and $i\!\in\!\omega$, then 
$\Big( g[f(\delta )]\Big)_{i}\! =\! h[f(\delta )](i)^\frown g_{i,h[f(\delta 
)](i)}[f(\delta )]\! =\! (\delta )_{i}(0)^\frown [(\delta )_{i}]^{*}\! =\! (\delta )_{i}$. 
Therefore $g[f(\delta )]\! =\!\delta$. This shows that $g$ coincides with $f^{-1}$.
\hfill{$\square$}

\vfill\eject

\noindent\bf Proof of Theorem 1.5.\rm\ Let $\gamma\!\in\!\Boraone\cap BC_{\xi +1}$ with 
$B\! =\!\neg\rho^{2^\omega}(\gamma )$.  By Theorem 4.1,\bigskip 

\noindent $\bullet$ If $\xi\! =\! 1$, then $F^0(\gamma )\!\in\!\Boraone$, 
$C\!\in\!\Bormone$, $f$ is a partial recursive function on $C$, and $g$ is a partial 
$\Boraxi$-recusive function on $B$.\bigskip

\noindent $\bullet$ If $\xi\!\geq\! 2$, then $F^1(\gamma )\!\in\!\Boraone$ and the same conclusion holds.\hfill{$\square$}\bigskip

 We also have a $\Borel$ version of Theorem 1.5:

\begin{thm} Let $\xi\!\geq\! 1$ be a countable ordinal, and $B\!\in\!\bormxipo (2^\omega )\cap\Borel$. 
Then there is $C\!\in\!\bormone\cap\Borel (\omega^\omega )$, a $\Borel$-recursive function 
$f\! :\!\omega^\omega\!\rightarrow\! 2^\omega$, and a $\Borel$-recursive function 
$g\! :\! 2^\omega\!\rightarrow\!\omega^\omega$, such that $f_{\vert C}$ defines a continuous 
bijection from $C$ onto $B$, $g_{\vert B}$ is $\boraxi$-measurable, and $g_{\vert B}$ 
coincides with $(f_{\vert C})^{-1}$.\end{thm}

\noindent\bf Proof.\rm\ We set $\boraxi (\Borel)(X)\! =\!
\{\rho^{X}(\gamma )\mid\gamma\!\in\!\Borel\cap BC_{\xi}\}$ if 
$1\!\leq\!\xi\! <\!\omega_1$. In \cite{Louveau80}, it is essentially proved that 
$\boraxi (\Borel)\! =\!\boraxi\cap\Borel$. Actually, Louveau does not use 
the coding for Borel sets that  we use here, but he proves this specific result, with this coding, in  his notes \cite{Louveau}. 
So let $\gamma\!\in\!\Borel\cap BC_{\xi +1}$ with 
$B\! =\!\neg\rho^{2^\omega}(\gamma )$.  By Theorem 4.1,\bigskip 

\noindent $\bullet$ If $\xi\! =\! 1$, then $F^0(\gamma )\!\in\!\Borel$, 
$C\!\in\!\bormone\cap\Borel$, $f$ and $g$ are partial $\Borel$ functions on 
$\Borel$ sets, and can be extended to total $\Borel$ maps.\bigskip

\noindent $\bullet$ If $\xi\!\geq\! 2$, then $F^1(\gamma )\!\in\!\Borel$ and the same conclusion holds.\hfill{$\square$}

\section{$\!\!\!\!\!\!$ Proof of Theorem 1.3.}\indent
 
 The proof of Theorem 1.3.(2) is essentially identical to that of Theorem 1.3.(1), so 
it is enough to prove Theorem 1.3.(1) to get Theorem 1.3. In the sequel we will assume that 
$\xi\! <\!\omega_{1}^{CK}$, except where indicated. Let us indicate the specifications of the proof of Theorem 1.2 that we need. Theorem 3.14 gives 
$B\!\in\! {\it\Gamma}(2^\omega )\!\setminus\!\check {\bf\Gamma}$. As $B\!\in\!\Bormxipo$, 
Theorem 1.5 gives $C$, $f$ and $g$. Here again, the dictionary $A$ will be made of two pieces: we will have $A\! =\!\mu\cup\pi$ if $\xi\!\geq\! 3$.\bigskip

\noindent\bf Notation.\rm\ Recall that $Q\! :=\!\{ (s,t)\!\in\! 2^{<\omega}\!\times\! 
2^{<\omega}\mid \vert s\vert\! =\!\vert t\vert\}$. We will sometimes view $Q$ as 
$\tilde Q\!\in\!\Borone (\omega )$: 
$$\tilde Q\! :=\!\Big\{ m\!\in\!\omega\mid\hbox{\rm Seq}(m)\ \ \hbox{\rm and}\ \ 
\forall i\! <\!\hbox{\rm lh}(m)\ \ \Big[\hbox{\rm Seq}[(m)_i]\ \ \hbox{\rm and}\ \ 
\hbox{\rm lh}[(m)_i]\! =\! 2~\ \hbox{\rm and}~\ 
\forall j\!\in\! 2\ \ \Big( (m)_i\Big)_j\!\! <\! 2\Big]\Big\}.$$
Implicitely, we have used the bijection $I\! :\! Q\!\rightarrow\!\tilde Q$ 
defined by 
$$I(s,t)\! :=\!\Big< <s(0),t(0)>\ ,\ \ldots\ ,\ <s(\vert s\vert\! -\! 1),
t(\vert s\vert\! -\! 1)>\Big>.$$ 
Note that the map $\chi\! :\!\omega\!\rightarrow\!\omega$ defined by $\chi (r)\! :=\! I(q_{r})$ is a recursive injection with range $\tilde Q$. We define a recursive map $M\! :\!\omega\!\rightarrow\!\omega$ by 
$M(j)\! :=\! M_{j}\! :=\!\Sigma_{i<j}\ 4^{i+1}$.

\begin {lem} The sets $\mu^0$, $\mu^1$ and $\mu$ can be coded by recursive subsets of 
$\omega$.\end{lem}

\noindent {\bf Proof.} We define a recursive map 
$\hbox{\rm Exp}\! :\!\omega^2\!\rightarrow\!\omega$ coding the finite sequence 
$k^j$: 
$$\hbox{\rm Exp}(k,j)\! :=\! c~\ \Leftrightarrow\ 
\hbox{\rm Seq}(c)\ \ \hbox{\rm and}\ \ \hbox{\rm lh}(c)\! =\! j\ \ \hbox{\rm and}
\ \ \forall i\! <\! j\ \ (c)_i\! =\! k.$$
Using Exp, it is easy to build a recursive map 
$f\! :\!\omega^5\!\rightarrow\!\omega$ such that $f(N,l,m,P,R)$ codes the sequence 
${\bf 2}^N\ {^\frown}\ [\ {^\frown}_{i\leq l+1}\ \ m_i\ {\bf 2}^{P_i}\ {\bf 3}\ {\bf 2}^{R_i}\ ]$. Then we just have to use bounded quantifiers.\hfill{$\square$}

\vfill\eject

 Now we show that $\mu^\infty$ is ``simple".

\begin {lem} The set $\mu^\infty$ is $\Bormtwo (4^\omega )$.\end{lem}

\noindent {\bf Proof.} We have
$$\gamma\!\in\!\mu^\infty\ \Leftrightarrow\ \exists i\!\in\! 2\ \ 
\forall j\!\in\!\omega\ \ \exists k\!\in\!\omega\ \ \exists t\!\in\! 
(\mu^i)^{<\omega}\ \ \vert t\vert\!\geq\! j\ \ \hbox{\rm and}\ \ 
\gamma\restriction k\! =\! {^\frown}_{l<\vert t\vert}\ t(l).$$
This shows that $\mu^\infty\!\in\!\Bormtwo (4^\omega )$, by Lemma 5.1, since $t$ 
can be coded by an integer, and the restriction and concatenation maps are 
recursive.\hfill{$\square$}\bigskip

\noindent\bf Notation.\rm\ We define a partial function 
$c\! :\! 2^\omega\!\times\!\omega\!\rightarrow\! Q$ on $B\!\times\!\omega$ by 
$c(\alpha ,l)\! :=\! [g(\alpha ),\alpha ]\restriction l$.

\begin {lem} The set 
$E\! :=\!\{ (N,\alpha )\!\in\!\omega\!\times\! 2^\omega\mid\alpha\!\in\! E_N\}$ is in 
$\it\Gamma$.\end{lem}

\noindent {\bf Proof.} The map $h\! :\!\omega\!\times\! 2^\omega\!\rightarrow\! 2^\omega$ defined by 
$h(N,\alpha )\! :=\! q^{1}_{N}\alpha$ is clearly recursive. From this we deduce that $E$ is in 
$\it\Gamma$, using Lemmas 3.3 and 3.5.\hfill{$\square$}\bigskip

\noindent\bf Notation.\rm\ Now we code the maps $\varphi_{N,j}$. We set 
$\hbox{\rm Dom}\! :=\!\{ (N,j,\gamma )\!\in\!\omega^{2}\!\times\! 
4^\omega\mid N\!\leq\! M_{j}\ \hbox{\rm and}\ \gamma\!\in\! K_{N,j}\}$. We define 
a partial function 
$\tilde\varphi\! :\!\omega^{2}\!\times\! 4^\omega\!\rightarrow\! 2^\omega$ as 
follows: $\tilde\varphi (N,j,\gamma )$ is defined if 
$(N,j,\gamma )\!\in\!\hbox{\rm Dom}$, and its coordinates are the coordinates of 
$\gamma$ in $2$, in the same order as in $\gamma$ (we forget the ${\bf 2}$'s and the 
${\bf 3}$'s).\bigskip 

 In the next lemma we consider the set expressing the fact that ``$\pi^\infty$ will look like $B$ 
on $K_{N,j}$".

\begin {lem} The set ${\cal F}\! :=\!\{ (N,j,\gamma )\!\in\!\hbox{\rm Dom}\mid
\tilde\varphi (N,j,\gamma )\!\in\! E_{N}\}$ is in $\it\Gamma$ if $\xi\!\geq\! 2$.
\end{lem}

\noindent {\bf Proof.} We define a map $\psi\! :\! 2^\omega\!\times\!\omega^{2}\!\rightarrow\! 4^\omega$ by 
$$\psi (\alpha ,N,j)\! :=\!\left\{\!\!\!\!\!\!
\begin{array}{ll}
& {\bf 2}^{N}\ \alpha (0)\ {^\frown}\  [\ {^\frown}_{k\in\omega}\ \ {\bf 2}^{M(j+k+1)}\ {\bf 3}\ 
{\bf 2}^{M(j+k+1)}\ \alpha (k\! +\! 1)\ ]\ \ \hbox{\rm if}\ \ N\!\leq\! M(j)
\hbox{\rm ,}\cr  
& 0^\infty\ \ \hbox{\rm if}\ \ N\! >\! M(j).
\end{array}\right.$$ 
It is easy to see that $\psi$ is recursive. If $N\!\leq\! M_{j}$ and 
$\gamma\!\in\! 4^\omega$, then $\gamma\!\in\! K_{N,j}$ is equivalent to 
$$\forall i\!\in\!\omega\ \ 
[\ \psi (0^\infty ,N,j)(i)\! =\! 0\ \hbox{\rm and}\ \gamma (i)\!\in\! 2\ ]
\ \ \hbox{\rm or}\ \ 
[\ \psi (0^\infty ,N,j)(i)\!\not=\! 0\ \hbox{\rm and}\ \gamma (i)\! =\! 
\psi (0^\infty ,N,j)(i)\ ].$$
This shows that $\hbox{\rm Dom}\!\in\!\Bormone$. Then $\tilde\varphi$ is 
clearly recursive on $\hbox{\rm Dom}$. This shows that ${\cal F}$ is in $\it\Gamma$ if 
$\xi\!\geq\! 2$, by Lemmas 3.3, 3.5 and 5.3.\hfill{$\square$}\bigskip

 Now we describe the elements of $A^\infty\!\setminus\!\mu^\infty$.\bigskip
 
\noindent\bf Notation.\rm\ Recall that $P_{t,S,j}:=\left\{\ \gamma\!\in\! 4^\omega\mid t\ {\bf 2}^S\!\prec\!\gamma\ \ \hbox{\rm and}\ \ \gamma\! -\! t\ {\bf 2}^S\!\in\! K_{0,j}\ \right\}$. Note that the relation defined by 
``$\gamma\!\in\! P_{t,S,j}$" is $\Bormone$ in $\gamma ,t,S,j$. Let $(t,S,j)$ be suitable and 
$N\!\leq\!\hbox{\rm min}(M_{j},S)$ ($N\! =\! S$ if $t\! =\!\emptyset$). Note that 
$(N,j,\gamma\! -\! t\ {\bf 2}^{S-N})\!\in\! {\cal F}$ means that 
$\gamma\! -\! t\ {\bf 2}^{S-N}\!\in\!\pi^\infty\cap K_{N,j}$. This implies that 
$$A_{t,S,j,N}=\left\{\ \gamma\!\in\! P_{t,S,j}\mid (N,j,\gamma\! -\! t\ {\bf 2}^{S-N})\!\in\! {\cal F}\ \right\}.$$ 

\begin {lem} The set of $(\gamma ,t,S,j,N)\!\in\! 4^\omega\!\times\! (\{\emptyset\}\cup\mu )\!
\times\!\omega^3$ such that $(t,S,j)$ is suitable, $N\!\leq\!\hbox{\it min}[M(j),S]$, 
$N\! =\! S\ \hbox{\it if}\ t\! =\!\emptyset$ and $\gamma\!\in\! A_{t,S,j,N}$ can be coded by a set in 
${\it\Gamma} (4^\omega\!\times\!\omega^4)$ if $\xi\!\geq\! 2$.\end{lem}

\noindent {\bf Proof.} Apply Lemmas 3.3, 3.5 and 5.4.\hfill{$\square$}

\vfill\eject

 Let us specify a few facts about the definition of $\pi$.\bigskip

\noindent\bf Notation.\rm\  As $C$ is $\Bormone$ and $f$ is recursive on $C$, the graph $\hbox{\rm Gr}(f)$ of $f$ is a $\Bormone$ subset of $\omega^\omega\!\times\! 2^\omega$. As the identity from $2^\omega$, viewed as a subset of $\omega^\omega$,  into $2^\omega$ is a partial recursive function on $2^\omega$ (see the proof of Theorem 4.1), we can also say that $\hbox{\rm Gr}(f)$ is a $\Bormone$ subset of 
$\omega^\omega\!\times\!\omega^\omega$, by Lemma 3.3.  By 4A.1 in \cite{Moschovakis}, there is 
$R\!\in\!\Borone (\omega^2)$ such that $\alpha\! =\! f(\beta )\ \Leftrightarrow\ 
\forall k\!\in\!\omega\ \ (\overline{\beta\restriction k},\overline{\alpha\restriction k})\!\in\! R$ (recall that $\overline{t}$ is defined at the beginning of section 3).\bigskip

\noindent $\bullet$ We set $Q_{f}\! :=\{ (t,s)\!\in\! Q\mid (\overline{t},\overline{s})\!\in\! R\ \hbox{\rm and}\ 
t\!\not=\!\emptyset\ \hbox{\rm and}\ t(\vert t\vert\! -\! 1)\! =\! 1\}$. Note that $Q_f$ can easily be coded by a recursive subset of $\omega$.\bigskip 

\noindent $\bullet$ The definition of $\pi$ is the same as the one in section 2. Here again, $\pi$ can easily be coded by a recursive subset of $\omega$.\bigskip 

\noindent {\bf Proof of Theorem 1.3.(1).} We refer to the proof of Theorem 1.2. We put $A\! :=\!\mu\cup\pi$, so that $A$ can be 
coded by a $\Borone$ subset of $\omega$. We will prove that 
$A^\infty\!\in\! {\it\Gamma}\!\setminus\!\check {\bf\Gamma}$.\bigskip

\noindent $\bullet$ Here again we have $\varphi_{N,j}[\pi^\infty\cap K_{N,j}]\! =\! E_N$ if 
$N\!\leq\! M_j$. If $\gamma\!\in\!\pi^\infty\cap K_{N,j}$, then the only thing to notice is that 
$[\overline{\beta\restriction k},\overline{(q^1_N\alpha )\restriction k}]\!\in\! R$ for each $k\!\in\!\omega$.\bigskip 

\noindent $\bullet$ We also have 
$$A^\infty\! =\!\mu^\infty\cup
\bigcup_{(t,S,j)\ \hbox{\rm suitable}}
\bigcup_{
\begin{array}{ll}
& N\leq\hbox{\rm min}(M_j,S)\cr 
& \ N=S\ \hbox{\rm if}\ t=\emptyset
\end{array}}
\ A_{t,S,j,N}\hbox{\rm .}$$ 
As $\it\Gamma$ is uniformly closed under finite unions, the set of $(\gamma ,t,S,j)\!\in\! 4^\omega\!\times\! (\{\emptyset\}\cup\mu )\!\times\!\omega^2$ such that $(t,S,j)$ is suitable and $\gamma\!\in\! A_{t,S,j}$ can be coded by a set in 
${\it\Gamma} (4^\omega\!\times\!\omega^3)$ if $\xi\!\geq\! 2$, by Lemma 5.5.\bigskip

\noindent $\bullet$ By Lemmas 3.5, 3.6 and 5.2, we get $A^\infty\!\in\! {\it\Gamma} (4^\omega )$ if 
$\xi\!\geq\! 3$ and ${\it\Gamma}\! =\!\Boraxi$.\bigskip

\noindent $\bullet$ If $\xi\!\geq\! 3$ and ${\it\Gamma}\! =\!\Bormxi$, then we can write
$$A^\infty\! =\!\mu^\infty\!\setminus\!\left(\bigcup_{(t,S,j)\ \hbox{\rm suitable}}\ P_{t,S,j}
\right)\ \cup\bigcup_{(t,S,j)\ \hbox{\rm suitable}}\ A_{t,S,j}\cap P_{t,S,j}.$$
Thus
$$\neg A^\infty\! =\!\neg\left[ \mu^\infty\cup\left(\bigcup_{(t,S,j)\ \hbox{\rm suitable}}\ P_{t,S,j}
\right)\right]\ \cup\bigcup_{(t,S,j)\ \hbox{\rm suitable}}\ P_{t,S,j}\!\setminus\! A_{t,S,j}.$$
Here $\neg\left[ \mu^\infty\cup\left(\bigcup_{(t,S,j)\ \hbox{\rm suitable}}\ P_{t,S,j}
\right)\right]\!\in\!\Borthree (4^\omega )\!\subseteq\!\check {\it\Gamma}(4^\omega )$. By Lemma 3.6, 
$\bigcup_{(t,S,j)\ \hbox{\rm suitable}}\ P_{t,S,j}\!\setminus\! A_{t,S,j}$ is in 
$\check {\it\Gamma}(4^\omega )$, and by Lemma 3.5 $\neg A^\infty$ is in 
$\check {\it\Gamma}(4^\omega )$. Thus $A^\infty\!\in\! {\it\Gamma} (4^\omega )$.\bigskip

\noindent $\bullet$ If $1\!\leq\!\xi\!\leq\! 2$, then we argue as in the proof of Theorem 1.2.\hfill{$\square$}

\section{$\!\!\!\!\!\!$ On the complexity of some sets of dictionaries.}

 The proof of Theorem 1.3 has the following consequence on the complexity of the sets ${\bf\Sigma}_{\xi}$ and ${\bf\Pi}_{\xi}$ defined in the introduction. Recall that if $1\!\leq\!\xi\! <\!\omega_1$, then
$${\bf\Sigma}_{\xi}\! :=\!\{ A\!\subseteq\! 2^{<\omega}\mid A^\infty\!\in\!\boraxi\}\ \ \ \ \ \ 
\hbox{\rm and}\ \ \ \ \ \ {\bf\Pi}_{\xi}\! :=\!\{ A\!\subseteq\! 2^{<\omega}\mid A^\infty\!\in\!\bormxi\}.$$
\bf Notation.\rm\ We set
$${\bf\Sigma}'_{\xi}\! :=\!\{\gamma\!\in\! BC\mid\rho^{2^\omega}(\gamma )\!\in\!\boraxi\}\ \ \ \ \ \ 
\hbox{\rm and}\ \ \ \ \ \ {\bf\Pi}'_{\xi}\! :=\!\{\gamma\!\in\! BC\mid\rho^{2^\omega}(\gamma )\!\in\!\bormxi\}.$$

\begin {cor} Let $3\!\leq\!\xi\! <\!\omega_1$. Then there is 
$\varphi\! :\!\omega^\omega\!\rightarrow\! 2^{2^{<\omega}}$ continuous with 
${\bf\Sigma}'_{\xi}\! =\! BC\cap\varphi^{-1}({\bf\Sigma}_{\xi})$ and 
${\bf\Pi}'_{\xi}\! =\! BC\cap\varphi^{-1}({\bf\Pi}_{\xi})$.\end{cor}

So ${\bf\Sigma}_{\xi}$ (resp., ${\bf\Pi}_{\xi}$) is more complicated than the set of Borel codes for 
$\boraxi$ (resp., $\bormxi$) sets, on $BC$, if $\xi\!\geq\! 3$.\bigskip

\noindent {\bf Proof.} Theorem 4.1 gives a partial function $F^1$. Recall that $F^1_1(\gamma )$ codes a continuous bijection defined on a closed subset of 
$\omega^\omega$ if $\gamma\!\in\! BC$. We now express the fact that its graph is a closed subset of 
$\omega^\omega\!\times\!\omega^\omega$ (see the notation after Lemma 5.5). In Theorem 4.1, the complement of $\rho^{2^\omega}(\gamma )$ is involved. This leads us to use the map $u_\neg$ given 
by Lemma 3.1. There is $P\!\in\!\Bormone [(\omega^\omega )^3]$ such that 
$$(\gamma ,\beta ,\alpha )\!\in\! P\ \Leftrightarrow\ \alpha\!\in\! 2^\omega\ \ \hbox{\rm and}\ \ 
\beta\!\notin\!\rho^{\omega^\omega}\Big( F^1_0[u_\neg(\gamma )]\Big)\ \ \hbox{\rm and}\ \ 
\alpha\! =\!\{ F^1_1[u_\neg(\gamma )]\}^{\omega^\omega ,2^\omega}(\beta )$$
if $\gamma\!\in\! BC$. By 4A.1 in \cite{Moschovakis} there is 
$\tilde R\!\in\!\Borone (\omega^\omega\!\times\!\omega^2)$ such that 
$$(\gamma ,\beta ,\alpha )\!\in\! P\ \Leftrightarrow\ \forall k\!\in\!\omega\ \
 [\gamma ,\overline{\beta\restriction k},\overline{\alpha\restriction k}]\!\in\!\tilde R$$
(see the notation after Lemma 5.5).\bigskip

 We say that $(t,s)\!\in\!\tilde Q_f$ if $(t,s)\!\in\! Q$, 
$[\gamma ,\overline{t},\overline{s}]\!\in\!\tilde R$, $t\!\not=\!\emptyset$ and $t(\vert t\vert\! -\! 1)\! =\! 1$  
(we use again the definition of $Q_f$ after Lemma 5.5, but here it is uniform in $\gamma$). Now we define $\pi$ as we did in section 2, with ``$q_{p_l}\!\in\!\tilde Q_f$" instead of ``$q_{p_l}\!\in\! Q_f$". After a coding of $4^{<\omega}$ with $\omega$, we can define a recursive map 
$\tilde\varphi\! :\!\omega^\omega\!\rightarrow\! 2^{\omega}$ coding $\mu\cup\pi\!\subseteq\! 4^{<\omega}$ 
(we will identify $\tilde\varphi (\gamma )$ with $\mu\cup\pi$, identifying $\omega$ with $4^{<\omega}$; the notation $\tilde\varphi$ instead of $\varphi$ is for $\omega$ in the range of $\tilde\varphi$ instead of 
$2^{<\omega}$ in the range of $\varphi$).\bigskip 

 Now let $\gamma\!\in\! BC$. Then $u_\neg (\gamma )\!\in\! BC$, $F^1[u_\neg (\gamma )]$ is defined, 
$f\! :\!\neg\rho^{\omega^\omega}\Big(F^1_0[u_\neg (\gamma )]\Big)\!\rightarrow\!\rho^{2^\omega}(\gamma )$ is a bijection. The proof of Theorem 1.3.(1) shows that $[\tilde\varphi (\gamma )]^\infty$ is $\boraxi$ (resp., $\bormxi$) if $\rho^{2^\omega}(\gamma)$ is $\boraxi$ (resp., $\bormxi$), when $\xi\!\geq\! 3$. It also shows that $\varphi_{0,0}\Big([\tilde\varphi (\gamma )]^\infty\cap K_{0,0}\Big)\! =\!\rho^{2^\omega}(\gamma)$, so that 
$\rho^{2^\omega}(\gamma)$ is $\boraxi$ (resp., $\bormxi$) if $[\tilde\varphi (\gamma )]^\infty$ is $\boraxi$ (resp., $\bormxi$), when $\xi\!\geq\! 3$.\hfill{$\square$}

\begin {cor} Let $B\!\in\!\borel [(2^\omega )^2]$ and $3\!\leq\!\xi\! <\!\omega_1$. Then there is $\psi\! :\! 2^\omega\!\rightarrow\! 2^{2^{<\omega}}$ continuous such that\smallskip

\noindent (a) ${\bf\Sigma}^B_{\xi}\! :=\!\{\alpha\!\in\! 2^\omega\mid B_\alpha\!\in\!\boraxi\}\! =\!
\psi^{-1}({\bf\Sigma}_{\xi})$.\smallskip

\noindent (b) ${\bf\Pi}^B_{\xi}\! :=\!\{\alpha\!\in\! 2^\omega\mid B_\alpha\!\in\!\bormxi\}\! =\!\psi^{-1}({\bf\Pi}_{\xi})$.\end{cor}

\noindent {\bf Proof.} (a) Let $\gamma_0\!\in\! BC$ such that 
$B\! =\!\rho^{(2^\omega )^2} (\gamma_0 )$. By Lemma 3.2, we get 
$\rho^{(2^\omega )^2} (\gamma_0 )_\alpha\! =\!
\rho^{2^\omega} [u^{2^\omega}_s(\gamma_0 ,\alpha )]$ for each 
$\alpha\!\in\! 2^\omega$. So we just have to set 
$\psi (\alpha )\! :=\!\varphi [u^{2^\omega}_s(\gamma_0 ,\alpha )]$, using Corollary 6.1.\bigskip

\noindent (b) The proof is similar.\hfill{$\square$}

\begin {thm} (Saint Raymond) Let $1\!\leq\!\xi\! <\!\omega_1$. Then there is 
$B\!\in\!\borel [(2^\omega )^2]$ such that ${\bf\Sigma}^B_{\xi}$ is $\ca$-complete. Similarly, there is 
$B\!\in\!\borel [(2^\omega )^2]$ such that ${\bf\Pi}^B_{\xi}$ is $\ca$-complete.\end{thm}

\noindent {\bf Proof.} Let $P\!\subseteq\! 2^\omega$ be a $\ca$-complete set, 
$G\!\in\!\bormtwo [(2^\omega )^2]$ such that $\neg P$ is the first projection of $G$, 
$X$ in $\borel (2^\omega)\!\setminus\!\boraxi$, and 
$B\! :=\!\{ (\alpha ,\beta )\!\in\! (2^\omega )^2\mid [\alpha ,(\beta )_0,(\beta )_1]\!\in\! G\!\times\! X\}$. Then $B$ is clearly Borel. If $\alpha\!\in\! P$, then $B_\alpha\! =\!\emptyset\!\in\!\boraxi$, so 
$\alpha\!\in\! {\bf\Sigma}^B_{\xi}$. If $\alpha\!\notin\! P$, let $\beta_0\!\in\! 2^\omega$ such that 
$(\alpha ,\beta_0)\!\in\! G$, and $f\! :\! 2^\omega\!\rightarrow\! 2^\omega$ defined by 
$f(\gamma )\! :=\ <\!\beta_0,\gamma\! >$. Then $B_\alpha\! =\!\{\beta\!\in\! 2^\omega\mid 
[\alpha ,(\beta )_0,(\beta )_1]\!\in\! G\!\times\! X\}\!\notin\!\boraxi$ since 
$X\! =\! f^{-1}(B_\alpha )\!\notin\!\boraxi$. Thus $\alpha\!\notin\! {\bf\Sigma}^B_{\xi}$. We proved that ${\bf\Sigma}^B_{\xi}\! =\! P$ is $\ca$-complete. We argue similarly for ${\bf\Pi}^B_{\xi}$.\hfill{$\square$}\bigskip

\noindent\bf Remarks.\rm\ (a) We actually proved that if $\xi\!\geq\! 3$ and $P\!\in\!\ca (2^\omega)$, then there is $M\!\in\!\bormxi [(2^\omega )^2]$ such that 
$P\! =\! {\bf\Sigma}^M_{\xi}$. Similarly, there is $A\!\in\!\boraxi [(2^\omega )^2]$ such that 
$P\! =\! {\bf\Pi}^A_{\xi}$.\bigskip

\noindent (b) This proof also shows that if $P\!\in\! {\bf\Pi}^1_2 (2^\omega)$, then there is 
$M\!\in\!\ca [(2^\omega )^2]$ such that $P\! =\! {\bf\Sigma}^M_{\xi}$. Similarly, there is 
$A\!\in\!\ca [(2^\omega )^2]$ such that $P\! =\! {\bf\Pi}^A_{\xi}$.\bigskip

\noindent (c) This proof also shows that if $P\!\in\! {\bf\Pi}^1_2 (2^\omega)$, then there is 
$C\!\in\!\ca [(2^\omega )^2]$ such that $P\! =\!\{\alpha\!\in\! 2^\omega\mid C_\alpha\!\in\!\borel\}$.

\begin {cor} Let $3\!\leq\!\xi\! <\!\omega_1$. Then ${\bf\Sigma}_{\xi}$ and ${\bf\Pi}_{\xi}$ are 
$\ca$-hard (and also ${\bf\Sigma}^1_2(2^{2^{<\omega}})\!\setminus\!\ana$).\end{cor}

\noindent {\bf Proof.} We just have to apply Theorem 6.3 and Corollary 6.2.\hfill{$\square$}\bigskip

\noindent\bf Remark.\rm\ Recall that if $X$ is a recursively presented Polish space and 
$\beta\!\in\! 2^\omega$, then 
$$\Ana (\beta )(X)\! :=\!\{ Q_\beta\mid Q\!\in\!\Ana (2^\omega\!\times\! X)\}\mbox{,}$$
$\Ca (\beta )\! :=\!\check\Ana (\beta )$ and $\Borel (\beta )\! :=\!\Ana (\beta )\cap\Ca (\beta )$. In \cite{Lecomte05}, the following sets are introduced: 
$$\begin{array}{ll}
{\it\Sigma}_{\xi} & \!\!\!\! :=\!\{ A\!\subseteq\! 2^{<\omega}\mid A^\infty\!\in\!\boraxi\cap\Borel (A)\}\hbox{\rm ,}\cr & \cr
{\it\Pi}_{\xi} & \!\!\!\! :=\!\{ A\!\subseteq\! 2^{<\omega}\mid A^\infty\!\in\!\bormxi\cap\Borel (A)\}.
\end{array}$$
It is proved in \cite{Lecomte05} that they are $\ca\!\setminus\!\borel$ if $\xi\!\geq\! 2$. Under the axiom of 
$\ana$-determinacy, this implies that they are $\ca$-complete. Here we can say more: they are $\ca$-complete if $\xi\!\geq\! 3$, without any axiom of determinacy. Indeed, fix a $\ca$-complete set $\Pi\!\subseteq\! 2^\omega$. The proof of Theorem 6.3 gives $B\!\in\!\borel [(2^\omega )^2]$ such that $B_\alpha\!=\!\emptyset$ if 
$\alpha\!\in\!\Pi$, and $B_\alpha\!\notin\!\boraxi$ if $\alpha\!\notin\!\Pi$. Now the proof of Corollary 6.2 gives $\gamma_0$. If $\alpha\!\in\!\Pi$, then 
$\rho^{2^\omega} [u^{2^\omega}_s(\gamma_0 ,\alpha )]\! =\!\emptyset$, and the proof of Theorem 1.3.(1) shows that 
$$[\psi (\alpha )]^\infty\! =\! (\varphi [u^{2^\omega}_s(\gamma_0 ,\alpha )])^\infty\! =\!\mu^\infty\!\in\!\Bormtwo\!\subseteq\!\Borel .$$ 
Thus $\psi (\alpha )\!\in\! {\it\Sigma}_{\xi}$ if $\alpha\!\in\!\Pi$. If $\alpha\!\notin\!\Pi$, then 
$\psi (\alpha )\!\notin\! {\bf\Sigma}_{\xi}$, thus $\psi (\alpha )\!\notin\! {\it\Sigma}_{\xi}$. Therefore 
$\Pi\! =\!\psi^{-1}({\it\Sigma}_{\xi})$ and ${\it\Sigma}_{\xi}$ is $\ca$-hard. As ${\it\Sigma}_{\xi}$ is $\ca$, it is $\ca$-complete. We argue similarly for ${\it\Pi}_{\xi}$.

\begin{defi} Let $\bf\Gamma$ be a class, and 
${\cal U}^{2^\omega}_{{\bf\Gamma}}\!\subseteq\! (2^\omega )^2$ universal for ${\bf\Gamma}(2^\omega )$. We say that ${\cal U}^{2^\omega}_{{\bf\Gamma}}$ is a 
$good\ universal$ for ${\bf\Gamma}$ if for each set 
${\cal U}^{(2^\omega )^2}_{{\bf\Gamma}}\!\subseteq\! (2^\omega )^3$ which is universal for 
${\bf\Gamma}[(2^\omega )^2]$, there is 
$S\! :\! (2^\omega )^2\!\rightarrow\! 2^\omega$ continuous such that 
$[S(\alpha ,\beta ),\gamma ]\!\in\! {\cal U}^{2^\omega}_{{\bf\Gamma}}\Leftrightarrow
(\alpha ,\beta ,\gamma )\!\in\! {\cal U}^{(2^\omega )^2}_{{\bf\Gamma}}$ 
for each $(\alpha ,\beta ,\gamma )\!\in\! (2^\omega )^3$.\end{defi} 

\begin{prop}  Let $1\!\leq\!\xi\! <\!\omega_1$. Then there are good universals for $\boraxi$,  $\bormxi$, $\ana$ and $\ca$.\end{prop}

\noindent\bf Proof.\rm\ Let ${\bf\Gamma}$ be one of the classes of the statement, and 
${\cal V}^{(2^\omega )^2}_{{\bf\Gamma}}$ universal for 
${\bf\Gamma}[(2^\omega )^2]$. We define, for $\alpha\!\in\! 2^\omega$ and 
$\varepsilon\!\in\! 2$, $(\alpha )_\varepsilon\!\in\! 2^\omega$ by 
$(\alpha )_\varepsilon (n)\! :=\!\alpha (2n\! +\!\varepsilon)$. We set 
$${\cal U}^{2^\omega}_{{\bf\Gamma}}\! :=\!\{ (\alpha ,\beta )\!\in\! (2^\omega )^2\mid 
[(\alpha )_0,(\alpha )_1,\beta ]\!\in\! {\cal V}^{(2^\omega )^2}_{{\bf\Gamma}}\}.$$
It is clear that ${\cal U}^{2^\omega}_{{\bf\Gamma}}\!\in\! {\bf\Gamma}$, so that 
$\{ ({\cal U}^{2^\omega}_{{\bf\Gamma}})_\alpha\mid\alpha\!\in\! 2^\omega\}\!\subseteq\!
{\bf\Gamma}(2^\omega )$. Conversely, let $A\!\in\! {\bf\Gamma}(2^\omega )$. Then the set 
$E\! :=\!\{ (\gamma ,\beta )\!\in\! (2^\omega )^2\mid\beta\!\in\! A\}\!\in\! {\bf\Gamma}$, so  there is $\alpha\!\in\! 2^\omega$ such that $E\! =\! ({\cal V}^{(2^\omega )^2}_{{\bf\Gamma}})_\alpha$. We define $<\! .,.\! >:\! (2^\omega )^2\!\rightarrow\! 2^\omega$ by 
$<\!\alpha ,\beta\! >(2n)\! :=\!\alpha (n)$ and 
$<\!\alpha ,\beta\! >(2n\! +\! 1)\! :=\!\beta (n)$. We get 
$A\! =\! ({\cal U}^{2^\omega}_{{\bf\Gamma}})_{<\alpha ,0^\infty >}$. We proved that 
${\cal U}^{2^\omega}_{{\bf\Gamma}}$ is universal for 
${\bf\Gamma}(2^\omega )$.\bigskip 

 Now let ${\cal U}^{(2^\omega )^2}_{{\bf\Gamma}}$ be universal for 
${\bf\Gamma}[(2^\omega )^2]$, and 
$$F\! :=\!\{ (\beta ,\gamma )\!\in\! (2^\omega )^2\mid 
[(\beta )_0,(\beta )_1,\gamma ]\!\in\! {\cal U}^{(2^\omega )^2}_{{\bf\Gamma}}\}.$$
As $F\!\in\! {\bf\Gamma}[(2^\omega )^2]$, there is $\alpha_0\!\in\! 2^\omega$ such that 
$F\! =\! ({\cal V}^{(2^\omega )^2}_{{\bf\Gamma}})_{\alpha_0}$. We get 
$$\begin{array}{ll}(\alpha ,\beta ,\gamma )\!\in\! {\cal U}^{(2^\omega )^2}_{{\bf\Gamma}}\!\!\!\!
& \Leftrightarrow (<\!\alpha ,\beta\! >,\gamma )\!\in\! F\cr & \cr
& \Leftrightarrow (\alpha_0,<\!\alpha ,\beta\! >,\gamma )\!\in\! 
{\cal V}^{(2^\omega )^2}_{{\bf\Gamma}}\cr & \cr
& \Leftrightarrow\Big(\Big<\alpha_0,<\!\alpha ,\beta\! >\!\Big>,\gamma\Big)\!\in\! 
{\cal U}^{2^\omega}_{{\bf\Gamma}}
\end{array}$$
So we just have to set $S(\alpha ,\beta )\! :=\!
\Big<\alpha_0,<\!\alpha ,\beta\! >\!\Big>$.\hfill{$\square$}

\begin{lem}  We consider the good universal ${\cal U}^{2^\omega}_{\ana}$ for $\ana$ given by Proposition 6.6. Then there is a continuous map $c\! :\! 2^{2^{<\omega}}\!\rightarrow\! 2^\omega$ such that 
$A^\infty\! =\! ({\cal U}^{2^\omega}_{\ana})_{c(A)}$ for each $A\!\in\! {\cal P}(2^{<\omega})\!\equiv\! 2^{2^{<\omega}}$.\end{lem}

\noindent\bf Proof.\rm\ Recall that ${\cal U}^{2^\omega}_{\boraone}\!\subseteq\! (2^\omega )^2$ is universal for $\boraone (2^\omega )$ and defined in the proof of Theorem 3.14 as follows:
$$(\gamma ,\alpha )\!\in\! {\mathcal U}^{2^\omega}_{\boraone}\ \ \Leftrightarrow
\ \ \exists  m\!\in\!\omega\ \ \gamma (m)\! =\! 0\ \ \hbox{\rm and}\ \ \alpha\!\in\! N[2^\omega ,m].$$
Similarly, we can define ${\cal U}^{2^\omega\times\omega^\omega}_{\boraone}\!\subseteq\! (2^\omega )^2\!\times\!\omega^\omega$, universal for $\boraone (2^\omega\!\times\!\omega^\omega )$:
$$(\gamma ,\alpha ,\beta )\!\in\! {\mathcal U}^{2^\omega\times\omega^\omega}_{\boraone}\ \ \Leftrightarrow
\ \ \exists  m\!\in\!\omega\ \ \gamma (m)\! =\! 0\ \ \hbox{\rm and}\ \ (\alpha ,\beta )\!\in\! N[2^\omega\!\times\!\omega^\omega ,m].$$
Using this, we can define ${\cal V}^{2^\omega}_{\ana}\!\subseteq\! (2^\omega )^2$, universal for $\ana (2^\omega )$:
$$(\gamma ,\alpha )\!\in\! {\mathcal V}^{2^\omega}_{\ana}\ \ \Leftrightarrow
\ \ \exists  \beta\!\in\!\omega^\omega\ \ (\gamma ,\alpha ,\beta )\!\notin\! {\mathcal U}^{2^\omega\times\omega^\omega}_{\boraone}.$$

\vfill\eject

By \cite{Lecomte05} there is a continuous map 
$\pi\! :\! 2^\omega\!\times\!\omega^\omega\!\times\!\omega\!\rightarrow\! 2^{<\omega}$ 
such that 
$$\alpha\!\in\! A^\infty ~\Leftrightarrow ~\exists\beta\!\in\!\omega^\omega ~~\forall n\!\in\!\omega ~~
[\beta (n\! +\! 1)\! >\! 0\ ~\hbox{\rm and}~\ \pi (\alpha ,\beta ,n)\!\in\! A]\hbox{\rm ,}$$
for each $\alpha\!\in\! 2^\omega$ and $A\!\subseteq\! 2^{<\omega}$. We define 
$R\!\in\!\boraone (2^\omega\!\times\!\omega^\omega\!\times\! 2^{2^{<\omega}})$ by 
$$(\alpha ,\beta ,A)\!\in\! R\ \Leftrightarrow\ \exists n\!\in\!\omega\ \ 
[\beta (n\! +\! 1)\! =\! 0\ ~\hbox{\rm or}~\ \pi (\alpha ,\beta ,n)\!\notin\! A].$$

 By 3C.5 in \cite{Moschovakis}, there is $R^*\!\subseteq\!\omega$ such that 
$$(\alpha ,\beta ,A)\!\in\! R\ \Leftrightarrow\ \exists m\!\in\!\omega\ \ \Big(\alpha\!\in\! N[2^\omega ,(m)_1]\ \ 
\hbox{\rm and}\ \ \beta\!\in\! N[\omega^\omega ,(m)_2]\ \ \hbox{\rm and}\ \ 
A\!\in\! N[2^{2^{<\omega}},(m)_3]\ \ \hbox{\rm and}\ \ m\!\in\! R^*\Big).$$
We define $d\! :\! 2^{2^{<\omega}}\!\rightarrow\! 2^\omega$ by 
$d(A)(m)\! =\! 0\ \Leftrightarrow\ A\!\in\! N[2^{2^{<\omega}},(m)_3]\ \ \hbox{\rm and}\ \ m\!\in\! R^*$. If 
$A\!\subseteq\! 2^{<\omega}$, then
$$\begin{array}{ll}
\alpha\!\in\! ({\cal V}^{2^\omega}_{\ana})_{d(A)} 
& \Leftrightarrow\ \exists\beta\!\in\!\omega\ \ [d(A),\alpha ,\beta ]\!\notin\! 
{\mathcal U}^{2^\omega\times\omega^\omega}_{\boraone}\cr 
& \Leftrightarrow\ \exists\beta\!\in\!\omega\ \ \neg\Big(\exists  m\!\in\!\omega\ \ d(A)(m)\! =\! 0\ \ 
\hbox{\rm and}\ \ (\alpha ,\beta )\!\in\! N[2^\omega\!\times\!\omega^\omega ,m]\Big)\cr 
& \Leftrightarrow\ \exists\beta\!\in\!\omega\ \ \neg\Big(\exists  m\!\in\!\omega\ \ d(A)(m)\! =\! 0\ \ 
\hbox{\rm and}\ \ \alpha\!\in\! N[2^\omega ,(m)_1]\ \ \hbox{\rm and}\ \ 
\beta\!\in\! N[\omega^\omega ,(m)_2]\Big)\cr 
& \Leftrightarrow\ \exists\beta\!\in\!\omega\ \ (\alpha ,\beta ,A)\!\notin\! R\cr 
& \Leftrightarrow\ \alpha\!\in\! A^\infty
\end{array}$$
As ${\cal V}^{2^\omega}_{\ana}\!\in\!\ana [(2^\omega)^2]$, there is $\alpha_0\!\in\! 2^\omega$ such that 
${\cal V}^{2^\omega}_{\ana}\! =\! ({\cal U}^{(2^\omega )^2}_{\ana})_{\alpha_0}$. As 
${\cal U}^{2^\omega}_{\ana}$ is a good universal, we get $S$ continuous, and 
$({\cal V}^{2^\omega}_{\ana})_{d(A)}\! =\! ({\cal U}^{2^\omega}_{\ana})_{S[\alpha_0,d(A)]}$. So we just have to set $c(A)\! :=\! S[\alpha_0,d(A)]$.\hfill{$\square$}\bigskip

 Recall that ${\cal U}({\bf\Gamma},{\bf\Gamma'})\! :=\!\{\alpha\!\in\! 2^\omega\mid ({\cal U}^{2^\omega}_{\bf\Gamma})_\alpha\!\in\! {\bf\Gamma'}\}$ and 
${\bf\Delta}\! :=\!\{ A\!\subseteq\! 2^{<\omega}\mid A^\infty\!\in\!\borel\}$.
 
\begin {cor} Let $3\!\leq\!\xi\! <\!\omega_1$. We consider the good universals given by Proposition 6.6.\medskip 

\noindent (a) The set ${\cal U}(\bormxi ,\boraxi )$ is $\ca$-complete, 
${\cal U}(\bormxi ,\boraxi )\leq_W {\bf\Sigma}_{\xi} <_W {\cal U}(\ana ,\boraxi )$, and the set 
${\cal U}(\ana ,\boraxi )$ is ${\bf\Pi}^1_2$-hard and 
${\bf\Sigma}^1_3\!\setminus\! {\bf\Sigma}^1_2$.\medskip

\noindent (b) The set ${\cal U}(\boraxi ,\bormxi )$ is $\ca$-complete, 
${\cal U}(\boraxi ,\bormxi )\leq_W {\bf\Pi}_{\xi} <_W {\cal U}(\ana ,\bormxi )$, and the set 
${\cal U}(\ana ,\bormxi )$ is ${\bf\Pi}^1_2$-hard and 
${\bf\Sigma}^1_3\!\setminus\! {\bf\Sigma}^1_2$.\medskip 

\noindent (c) ${\bf\Delta} <_W {\cal U}(\ana ,\borel )$, and the set ${\cal U}(\ana ,\borel )$ is ${\bf\Pi}^1_2$-hard and ${\bf\Sigma}^1_3\!\setminus\! {\bf\Sigma}^1_2$. Moreover, the set 
${\cal U}(\ana ,\Borel )$ is ${\bf\Pi}^1_2$-complete.\end{cor}

\noindent {\bf Proof.} (a) By Theorem 6.3 and Remark (a) just after, there is 
$M\!\in\!\bormxi [(2^\omega )^2]$ such that ${\bf\Sigma}^M_{\xi}$ is $\ca$-complete. Fix 
$\alpha_0\!\in\! 2^\omega$ with $M\! =\! ({\cal U}^{(2^\omega )^2}_{\bormxi})_{\alpha_0}$. We define 
$f\! :\! 2^\omega\!\rightarrow\! 2^\omega$ by $f(\alpha )\! :=\! S(\alpha_0,\alpha )$, where $S$ is provided by the fact that ${\cal U}^{2^\omega}_{\bormxi}$ is a good universal. Then we get 
${\bf\Sigma}^M_{\xi}\! =\! f^{-1}({\cal U}(\bormxi ,\boraxi ))$, which proves that 
${\cal U}(\bormxi ,\boraxi )$ is $\ca$-hard. By \cite{Louveau80} (or 35.H in \cite{Kechris}), ${\cal U}(\bormxi ,\boraxi )$ is 
$\ca$, so it is $\ca$-complete.\bigskip

 By Corollary 6.2, we get ${\cal U}(\bormxi ,\boraxi )\leq_W {\bf\Sigma}_{\xi}$ since 
$${\cal U}(\bormxi ,\boraxi )\!=\! {\bf\Sigma}^{{\cal U}^{2^\omega}_{\bormxi}}_{\xi}.$$

\vfill\eject

By Lemma 6.7 we get ${\bf\Sigma}_{\xi}\leq_W {\cal U}(\ana ,\boraxi )$. Remark (b) 
after Theorem 6.3 gives $\Sigma\! :=\!\neg A\!\in\!\ana [(2^\omega )^2]$ such that ${\bf\Sigma}^{\Sigma}_{\xi}$ is 
${\bf\Pi}^1_2$-complete. The beginning of the proof shows that ${\cal U}(\ana ,\boraxi )$ is ${\bf\Pi}^1_2$-hard. 
In particular, ${\cal U}(\ana ,\boraxi )\!\notin\! {\bf\Sigma}^1_2$, and 
${\bf\Sigma}_{\xi} <_W {\cal U}(\ana ,\boraxi )$ since ${\bf\Sigma}_{\xi}\!\in\! {\bf\Sigma}^1_2$. 
Finally, ${\cal U}(\ana ,\boraxi )$ is ${\bf\Sigma}^1_3$ since 
$$\alpha\!\in\! {\cal U}(\ana ,\boraxi )\Leftrightarrow\exists\beta\!\in\! 2^\omega\ 
({\cal U}^{2^\omega}_{\ana})_\alpha\! =\! ({\cal U}^{2^\omega}_{\boraxi})_\beta .$$
(b) The proof is very similar to that of (a).\bigskip

\noindent (c) The proof of the first sentence is very similar to that of (a), using Remark (c) after Theorem 6.3. 
This proof shows that ${\cal U}(\ana ,\Borel )$ is ${\bf\Pi}^1_2$-hard. It remains to see that 
${\cal U}(\ana ,\Borel )$ is ${\bf\Pi}^1_2$. Recall the existence of $\Ca$ sets 
$W^{2^\omega}\!\subseteq\!\omega$, $C^{2^\omega}\!\subseteq\!\omega\!\times\! 2^\omega$ with 
${\Borel (2^\omega )\! =\!\{ C^{2^\omega}_{n}\mid n\!\in\! W^{2^\omega}\}}$ and 
$$\{(n,\alpha )\!\in\!\omega\!\times\!\! 2^\omega\mid n\!\in\! W^{2^\omega}\mbox{ and }\alpha\!\notin\! C_{n}^{2^\omega}\}\!\in\!\Ca (\omega\!\times\! 2^\omega)$$ 
(see Theorem 3.3.1 in \cite{HKL}). This implies that 
$$\alpha\!\in\! {\cal U}(\ana ,\Borel )\Leftrightarrow\exists n\!\in\! W^{2^\omega}\ 
({\cal U}^{2^\omega}_{\ana})_\alpha\! =\! C^{2^\omega}_{n}.$$
Thus ${\cal U}(\ana ,\Borel )$ is ${\bf\Pi}^1_2$, and ${\bf\Pi}^1_2$-complete.\hfill{$\square$}

\hs {\bf Acknowledgements.}~We would like to thank J. Saint Raymond for allowing us to present here one of his results (see Theorem 6.3), and the referee for valuable 
comments on a first version of this paper.

\end{document}